\documentclass[a4paper,twoside,11pt]{article}

\usepackage{amsmath}
\usepackage{amssymb}
\usepackage[utf8x]{inputenc} 
\usepackage{accents}
\usepackage{mathtools}
\usepackage{amsthm}
\usepackage[english]{babel}
\usepackage{newpxtext}
\usepackage[euler-digits]{eulervm}

\DeclareFontFamily{U}{mathx}{\hyphenchar\font45}
\DeclareFontShape{U}{mathx}{m}{n}{
	<5> <6> <7> <8> <9> <10>
	<10.95> <12> <14.4> <17.28> <20.74> <24.88>
	mathx10
}{}
\DeclareSymbolFont{mathx}{U}{mathx}{m}{n}
\DeclareFontSubstitution{U}{mathx}{m}{n}
\DeclareMathAccent{\widecheck}{0}{mathx}{"71}
\DeclareMathAccent{\wideparen}{0}{mathx}{"75}

\usepackage{fullpage}
\usepackage{fancyhdr}

\usepackage[dvipsnames,svgnames,x11names,hyperref]{xcolor}
\usepackage{letltxmacro}
\usepackage{xparse}
\usepackage{relsize}
\usepackage{mathrsfs}
\usepackage[bbgreekl]{mathbbol}
\DeclareSymbolFontAlphabet{\mathbbm}{bbold}
\DeclareSymbolFontAlphabet{\mathbb}{AMSb}%
\usepackage{euscript} 
\usepackage{mleftright}
\usepackage{dirtytalk}
\usepackage{scalerel} 
\usepackage{comment}
\usepackage{soul}
\usepackage[disable]{todonotes}

\usepackage{capt-of}

\setul{0pt}{0.4pt} 

\AtBeginDocument{%
	\LetLtxMacro\autoreforig\autoref
	\RenewDocumentCommand{\autoref}{som}{%
		\IfValueT{#2}{[}%
		\IfBooleanTF{#1}{%
			\autoreforig*{#3}%
		}{%
			\autoreforig{#3}%
		}%
		\IfValueT{#2}{,\space#2]}%
	}
}

\makeatletter
\renewcommand*\env@matrix[1][\arraystretch]{%
	\edef\arraystretch{#1}%
	\hskip -\arraycolsep
	\let\@ifnextchar\new@ifnextchar
	\array{*\c@MaxMatrixCols c}}
\makeatother

\usepackage{commath}
\usepackage{xfrac}
\usepackage{array}
\usepackage{enumitem}
\usepackage{stmaryrd}
\usepackage{tensor}
\usepackage{wasysym}

\newcommand\independent{\protect\mathpalette{\protect\independenT}{\perp}}
\def\independenT#1#2{\mathrel{\rlap{$#1#2$}\mkern2mu{#1#2}}}

\usepackage[scr=boondoxo]{mathalfa}

\newcommand{\edif}{\emph{\dif}}

\newcommand{\bbR}{\mathbb R}
\newcommand{\bbX}{\mathbb X}
\newcommand{\bbY}{\mathbb Y}
\usepackage[euler]{textgreek}

\newcommand{\scrg}{\mathscr g}

\newcommand{\Xphi}{{^\varphi} \! X}

\newcommand{\HFi}{{^{\mathscr F}} \! H}

\newcommand{\bfH}{\boldsymbol{H}}
\newcommand{\bfK}{\boldsymbol{K}}
\newcommand{\oH}{\bfH}

\newcommand{\bfX}{\mathbf{X}}
\newcommand{\bfY}{\mathbf{Y}}
\newcommand{\bfZ}{\mathbf{Z}}

\newcommand{\uX}{\bfX}

\newcommand{\Ptr}{/\hspace{-2.5pt}/}
\newcommand{\Aptr}{{\setminus\hspace{-2.5pt}\setminus}}
\newcommand{\Dev}{\text{\setul{2.35pt}{0.45pt}\ul{$\bigcirc$}}}
\newcommand{\Adev}{\text{\setul{-8.275pt}{0.45pt}\ul{$\bigcirc$}}}

\newcommand{\ptr}[4]{\Ptr^{#2}_{#3} \hspace{-1.5pt}(#1)_{#4}}
\newcommand{\aptr}[4]{{\Aptr^{#2}_{#3}} \hspace{-1.5pt}(#1)_{#4}}
\newcommand{\dev}[4]{\Dev^{#2}_{#3} \hspace{-1.5pt}(#1)_{#4}}
\newcommand{\adev}[4]{\Adev^{#2}_{#3} \hspace{-1.5pt}(#1)_{#4}}

\newcommand{\de}[3]{\Dev^{#2} \hspace{-2pt}(#1)_{#3}}
\newcommand{\ade}[3]{\Adev^{#2} \hspace{-3.5pt}(#1)_{#3}}

\newsavebox{\overlongequation}

\usepackage{graphicx}
\makeatletter
\newcommand*\bigcdot{\mathpalette\bigcdot@{.5}}
\newcommand*\bigcdot@[2]{\mathbin{\vcenter{\hbox{\scalebox{#2}{$\m@th#1\bullet$}}}}}
\makeatother

\newenvironment{talign*}
{\let\displaystyle\textstyle\csname align*\endcsname}
{\endalign}

\usepackage{tikz}
\usepackage{graphicx}
\usepackage{wrapfig}
\usepackage{caption}
\usepackage[percent]{overpic}
\usepackage{tikz-cd}

\usepackage{aliascnt}
\usepackage[pdftex, ocgcolorlinks, colorlinks = true,
linkcolor = ProcessBlue,
citecolor = LimeGreen,
urlcolor = RoyalPurple
]{hyperref}

\numberwithin{equation}{section}

\addto\extrasenglish{
	
} 

\addto\extrasenglish{
	
} 

\usepackage[refpage, intoc]{nomencl}

\makenomenclature

\newcommand{\mynewtheorem}[2]{
	\newaliascnt{#1}{dummy}
	\newtheorem{#1}[#1]{#2}
	\aliascntresetthe{#1}
	\expandafter\def\csname #1autorefname\endcsname{#2}
}

\theoremstyle{plain}
\mynewtheorem{thm}{Theorem}
\mynewtheorem{thmdef}{Theorem/Definition}
\mynewtheorem{prop}{Proposition}
\mynewtheorem{cor}{Corollary}
\mynewtheorem{lem}{Lemma}
\mynewtheorem{conj}{Conjecture}
\theoremstyle{definition}
\mynewtheorem{defn}{Definition}
\mynewtheorem{expl}{Example}
\mynewtheorem{prob}{Problem}
\mynewtheorem{ques}{Question}
\mynewtheorem{cond}{Condition}
\mynewtheorem{conv}{Convention}
\mynewtheorem{defthm}{Definition/Theorem}
\theoremstyle{remark}
\mynewtheorem{rem}{Remark}

\renewcommand{\dif}{\text{d}}

\let\amsamp=& 

\mathchardef\mhyphen="2D

\mleftright

\usepackage{authblk}

\title{Non-Geometric Rough Paths on Manifolds}
\date{\today}
\author[2]{John Armstrong}
\author[1]{Damiano Brigo}
\author[1]{Thomas Cass}
\author[1]{Emilio Rossi Ferrucci\thanks{Corresponding author: \href{mailto:emilio.rossi-ferrucci16@imperial.ac.uk}{\nolinkurl{emilio.rossi-ferrucci16@imperial.ac.uk}}}}
\affil[1]{Dept.\ of Mathematics, Imperial College London}
\affil[2]{Dept.\ of Mathematics, King's College London}
\date{\today}

\usepackage{cutwin,kantlipsum,mwe}
\usepackage{regexpatch}
\opencutright

\makeatletter
\xpretocmd{\cutout}{\leavevmode\hrule \@height\z@ \@width\linewidth\relax}
\makeatother

\begin{document}

\maketitle

\begin{abstract}
We provide a theory of manifold-valued rough paths of bounded $3 > p$-variation, which we do not assume to be geometric. Rough paths are defined in charts, and coordinate-free (but connection-dependent) definitions of the rough integral of cotangent bundle-valued controlled paths, and of RDEs driven by a rough path valued in another manifold, are given. When the path is the realisation of semimartingale we recover the theory of It\^o integration and SDEs on manifolds \cite{E89}. We proceed to present the extrinsic counterparts to our local formulae, and show how these extend the work in \cite{CDL15} to the setting of non-geometric rough paths and controlled integrands more general than 1-forms. In the last section we turn to parallel transport and Cartan development: the lack of geometricity leads us to make the choice of a connection on the tangent bundle of the manifold $TM$, which figures in an It\^o correction term in the parallelism RDE; such connection, which is not needed in the geometric/Stratonovich setting, is required to satisfy properties which guarantee well-definedness, linearity, and optionally isometricity of parallel transport. We conclude by providing numerous examples, some accompanied by numerical simulations, which explore the additional subtleties introduced by our change in perspective.
\end{abstract}

\section*{Introduction}

The theory of rough paths, first introduced in \cite{L98}, has as its primary goal that of providing a rigorous mathematical framework for the study of differential equations driven by highly irregular inputs. The roughness of such signals renders the traditional definition of differentiation and integration inapplicable, and motivates the definition of \emph{rough path}, a path $X$ accompanied by functions, satisfying certain algebraic and analytic constraints, which postulate the values of its (otherwise undefined) iterated integrals. This concept leads to definitions of \emph{rough integration} against the rough path $\bfX$ and of \emph{rough differential equation} (RDE) driven by $\bfX$, which bear the important feature of being continuous in the signal $\bfX$, according to appropriately defined $p$-variation norms. Rough path theory applies to a wide variety of settings, including to the case in which $X$ is given by the realisation of a stochastic process, for which it constitutes a pathwise approach to stochastic integration, extending the classical stochastic analysis of semimartingales.

An important feature that a rough path can satisfy is that of being \emph{geometric}: this can be interpreted as the statement that it obeys the integration by parts and change of variable laws of first-order calculus, its irregularity notwithstanding. The theory of geometric rough paths has been the most studied \cite{FV10}, and applies to semimartingales through the use of the Stratonovich integral. Other notions of stochastic integration, however, cannot be modelled by geometric rough paths, the It\^o integral being the prime (but not the only \cite{ER03}) example. 

Since smooth manifolds are meant to provide a general setting for ordinary differential calculus to be carried out, it is natural to ask how \say{rougher} calculi can be defined in the curved setting. In the context of stochastic calculus, this question has led to a rich literature on Brownian motion on manifolds. More recently, it has been raised a number of times with regards to rough paths \cite{CLL12,DS17,CDL15,BL15,Bai19}. In all cases, however, only the case of geometric rough paths has been discussed. 

The main goal of this paper is to construct a theory of manifold-valued rough paths of bounded $p$-variation, with $p < 3$, which are not required to be geometric. The regularity assumption ensures that we may draw on the familiar setting of \cite{FH14} for vector space-valued rough paths; dropping this requirement would require the more complex algebraic tools of \cite{HK15}. Our theory includes defining rough integration and differential equations, both from the intrinsic and extrinsic points of view, and showing how the classical notions of parallel transport and Cartan development can be extended to the case of non-geometric rough paths.

Although the definition of the It\^o integral on manifolds has been known for decades, Stratonovich calculus has been preferred in the vast majority of the literature on stochastic differential geometry. Nevertheless, there are phenomena that are best captured by It\^o calculus, particularly those which relate to the martingale property. In this spirit, three of the authors recently showed how a concrete problem involving the approximation of SDEs with ones defined on submanifolds necessitates the use of It\^o notation, and that the result naturally provided by projecting the Stratonovich coefficients is suboptimal in general \cite{ABR18,ABR19}. The reason that Stratonovich integration and geometric rough paths are preferred in differential geometry is that they admit a simple coordinate-free description, as is also remarked on \cite[p.219]{L98}. An important point, however, that we wish to make in this paper is the following: an invariant theory of integration against non-geometric rough paths may also be given, albeit one that depends on the choice of a linear connection on the tangent bundle of the manifold. Although geometric rough path theory still retains the important property of being connection-invariant, all rough paths may be treated in a coordinate-free manner, since, while manifolds may not admit global coordinate systems, they always admit covariant derivatives. Overlooking this principle leads to the common misconception that It\^o calculus/non-geometric rough integration cannot be carried out on manifolds, even in cases where a connection is already independently and canonically specified, e.g.\ when the manifold is Riemannian. In much of stochastic differential geometry the focus is not on the stochastic integral per se, which is viewed as a tool to investigate laws of processes defined on Riemannian Wiener space: in this context it is certainly justifiable to only work with the Stratonovich integral. Our emphasis here, however, is on pathwise integration itself: for this reason we believe it to be of value to build up the theory in a way that is faithful to the choice of the calculus, as specified through the rough path $\bfX$.

This paper is organised as follows:\todo{modify appropriately once paper is finished} in \autoref{sec:backRps} we review the theory of vector space-valued rough paths of bounded $3>p$-variation, controlled rough integrations and RDEs, relying (with a few modifications and additions) on \cite{FH14}.

In \autoref{sec:backDG} we review the differential geometry necessary in the following sections.

In \autoref{sec:rough} we develop the theory at the heart of the paper: this entails defining rough paths on manifolds and their controlled integrands in a coordinate-free manner by using pushforwards and pullbacks through charts, showing how the choice of a linear connection gives rise to a definition of rough integral, and defining RDEs in a similar spirit. We follow the \say{transfer principle} philosophy \cite{E90} of replacing all instances of Euclidean spaces with smooth manifolds, which means that both the driving rough path and the solution are valued in (possibly different) manifolds. When we restrict our theory to semimartingales we recover the known framework for It\^o integration and stochastic differential equations (SDEs) on manifolds \cite{E89}.

In \autoref{sec:extrinsic} we switch from the local to the extrinsic framework, and show how our theory extends that of \cite{CDL15} to non-geometric integrators and controlled integrands more general than 1-forms. Our broader assumptions require us to make additional nondegeneracy requirements on the path $X$, which are not needed in the local setting. We also remark that in this section we are confining ourselves to the Riemannian case (with the metric being induced by an embedding), while in the rest of the paper we allow for general connections.

Finally, in \autoref{sec:par} we return to our local coordinate framework to carry out the constructions of parallel transport along rough paths and the resulting notion of Cartan development, or \say{rolling without slipping}, a cornerstone of stochastic differential geometry which yields a convenient way of moving back and forth between the linear and curved setting. Since we are dealing with parallel transport as a $TM$-valued RDE driven by an $M$-valued rough path $\bfX$, the lack of geometricity leads us to require the choice of a connection not just on the tangent bundle of $M$ but also of one on the tangent bundle of the manifold $TM$. The latter connection may not be chosen arbitrarily, and we identify criteria (formulated in terms of the former connection) that guarantee well-definedness, linearity, and, if $M$ is Riemannian, isometricity of parallel transport. Different choices of such connection give rise to different definitions of parallel transport and Cartan development, which are only detectable at a second-order level, and all collapse to the same RDE when the rough path is geometric. Though we develop the theory in the most general way possible, three examples for how a connection on $TM$ may be lifted to one on $TTM$ are drawn from the literature; a case not analysed until now concerns the Levi-Civita connection of the Sasaki metric, which results in parallel transport coinciding with Stratonovich parallel transport. We end by seizing the opportunity to explore a few additional topics in stochastic analysis on manifolds, such as Cartan development in the presence of torsion, with a pathwise emphasis.

We hope that the framework laid out in this paper may be used in the future to extend our understanding of manifold-valued rough paths, both deterministic and stochastic, and in \autoref{sec:concl} offer some ideas in this direction.

\section{Background on rough paths}\label{sec:backRps}
In this section we review the core theory of finite-dimensional vector space-valued (controlled) rough paths, and the corresponding notions of rough integrals and RDEs. We refer mainly to \cite{FH14}, with the caveat that we are in the setting of arbitrary control functions, as opposed to H\"older regularity. The former has has the advantage of being a parametrisation-invariant framework, and of allowing us to consider a larger class of paths (e.g.\ all semimartingales, and not just Brownian motion). Other authors have already been treating controlled rough paths in the setting of bounded $p$-variation \cite[\S 2.4]{CL16}. When a result in this first section is stated without proof, it is understood that the proof can be found in \cite[Ch.\ 1-10]{FH14}, possibly with trivial modifications needed to adapt the arguments to the case of arbitrary controls.\todo{I have made this substitution in the statements of \cite{FH14} without worrying too much about not introducing falsehoods. I've seen this done elsewhere, but it could get contested.} Many of the more quantitative aspects of rough paths are left out, as they will not be relevant for the transposition of the theory to manifolds. Since our vector spaces are finite-dimensional, and since we will rely on arbitrary charts to make the manifold-valued theory coordinate-free, we will use fixed coordinates to express all of our formulae.

\subsection{$\bbR^d$-valued rough paths}\label{subsec:rdRough}
Throughout this document $p$ will be a real number $\in [1,3)$; we will not exclude the case of $p \in [1,2)$ in which the theory reduces to Young integration, and remains valid with trivial adjustments. A \emph{control} on $[0,T]$ is a continuous function $\omega$ defined on the subdiagonal $\Delta_T \coloneqq \{(s,t) \in [0,T]^2 \mid s \leq t\}$, s.t.\ $\omega(t,t) = 0$ for $0 \leq t \leq T$ and $\omega(s,u) + \omega(u,t) \leq \omega(s,t)$ for $0 \leq s \leq u \leq t \leq T$. $\omega$ will denote a control throughout this document, and should be thought of as being a fixed property of the (rough) path which relates to its parametrisation; the main example is the \emph{H\"older control} $\omega(s,t) = t-s$. Given a path $X \colon [0,T] \to \bbR^d$ we will denote its increment $X_{st} \coloneqq X_t - X_s$. Let $\mathcal C^p_\omega([0,T],\bbR^d)$ denote the set of $\bbR^d$-valued continuous paths $X \colon [0,T] \to \bbR^d$ with
\begin{equation}
\sup_{0 \leq s < t \leq T} \frac{\abs{X_{st}}}{\omega(s,t)^{1/p}} < \infty
\end{equation}
For there to exist a control $\omega$ s.t.\ the above holds is equivalent to saying that $X$ is a path of bounded $p$-variation \cite[Proposition 5.10]{FV10}; if $\omega$ is the H\"older control we recover the definition of H\"older regularity. This kind of regularity is invariant under smooth maps:
\begin{lem}\label{lem:regSmooth}
	Let $\omega$ be a control, $X \in \mathcal C^p_\omega([0,T],\bbR^d)$ and $f \in C^\infty(\bbR^d,\bbR^e)$. Then $f(X) \in \mathcal C^p_\omega([0,T],\bbR^e)$.
	\begin{proof}
		We have
		\begin{equation}
		\begin{split}
		\sup_{0 \leq s < t \leq T} \frac{\abs[0]{f(X)_{st}}}{\omega(s,t)^{1/p}} &= \sup_{0 \leq s < t \leq T} \frac{\abs[0]{\partial_\gamma f(X_s)X^\gamma_{st} + O(\abs[0]{X_{st}}^2)}}{\omega(s,t)^{1/p}} \\
		&\leq \norm[0]{Df|_{X_{[0,T]}}}_\infty \sup_{0 \leq s < t \leq T} \frac{\abs[0]{X_{st}}}{\omega(s,t)^{1/p}} + \sup_{0 \leq s < t \leq T} \frac{\abs[0]{O(\omega(s,t)^2)}}{\omega(s,t)^{1/p}} \\
		&< \infty
		\end{split}
		\end{equation}
		In writing $\partial_\gamma f(X_s)X^\gamma_{st}$ we have used the Einstein summation convention, as shall be done throughout this paper. This concludes the proof.
	\end{proof}
\end{lem}

Recall that if $ X \in \mathcal C^p_\omega([0,T],\bbR^d)$, $H \in \mathcal C^q_\omega([0,T],\bbR^{e \times d})$ with $1/p + 1/q > 1$ (which happens, in particular, when $p = q \in [1,2)$) we may define the Young integral
\begin{equation}
\int_s^t H \dif X \coloneqq \lim_{n \to \infty} \sum_{[u,v] \in \pi_n} H_u X_{uv}
\end{equation}
where $(\pi_n)_n$ is a sequence of partitions on $[s,t]$ with vanishing step size; the resulting path $\int_0^\cdot H \dif X$ belongs to $\mathcal C^p_\omega([0,T],\bbR^d)$. When the regularity requirement is no longer satisfied the Riemann sums no longer converge, and the definition of integral will require $X$ and $H$ to carry additional structure.

\begin{defn}[Rough path]
	A $p$-\emph{rough path} controlled by $\omega$ on $[0,T]$, valued in $\bbR^d$ consists of a pair $\bfX = (X,\bbX)$ with $X \in \mathcal C^p_\omega([0,T],\bbR^d)$ (the \emph{trace}) and a continuous function $\bbX \colon \Delta_T \to (\bbR^d)^{\otimes 2} = \bbR^{d\times d}$ (the \emph{second order part}) satisfying the regularity condition
	\begin{equation}
	\sup_{0 \leq s < t \leq T} \frac{\abs{\bbX_{st}}}{\omega(s,t)^{2/p}} < \infty
	\end{equation}
	with the property that the \emph{Chen identity} holds: for all $0 \leq s \leq u \leq t \leq T$ and $\alpha,\beta = 1,\ldots,d$
	\begin{equation}
	\bbX_{st}^{\alpha\beta} = \bbX_{su}^{\alpha\beta} + X_{su}^\alpha X_{ut}^\beta + \bbX_{ut}^{\alpha\beta}
	\end{equation}
	We denote the set of all such $\bfX$ as $\mathscr C^p_\omega([0,T],\bbR^d)$ (note the difference in font with $\mathcal C$, used for simple paths). Its \emph{bracket} path is given by
	\begin{equation}\label{eq:bracket}
	[\bfX]_{st}^{\alpha\beta} \coloneqq X^\alpha_{st}X^\beta_{st} - (\bbX^{\alpha\beta}_{st} + \bbX^{\beta\alpha}_{st})
	\end{equation}
	These are indeed the increments of an element of $\mathcal C^{p/2}_\omega([0,T],(\bbR^d)^{\odot 2})$-valued path, where $\odot$ denotes symmetric tensor product. We will say that $\bfX$ is \emph{geometric} if $[\bfX] = 0$, and denote the set of these with $\mathscr G^p_\omega([0,T],\bbR^d)$.
\end{defn}
The idea is that $\bbX_{st}$ represents the value of the (otherwise undefined) integral 
\begin{equation}\label{eq:bbXint}
\int_s^t \int_s^u \dif X_r \otimes \dif X_u = \int_s^t X_{su} \otimes \dif X_u
\end{equation}
In this interpretation it is easily checked that the Chen relation is simply the statement that the integral $\int X_u \otimes \dif X_u$ is additive on consecutive time intervals, and the property of $\bfX$ of being geometric represents an integration by parts formula. Relaxing these two requirements to 
\begin{equation}
\bbX_{st}^{\alpha\beta} = \bbX_{su}^{\alpha\beta} + X_{su}^\alpha X_{ut}^\beta + \bbX_{ut}^{\alpha\beta} + \varepsilon_{st}^{\alpha\beta}\quad \text{(and} \quad X^\alpha_{st}X^\beta_{st} = \bbX^{\alpha\beta}_{st} + \bbX^{\beta\alpha}_{st} + \varepsilon_{st}^{\alpha\beta}\text{)} 
\end{equation}
for some function of two parameters $\varepsilon_{st} \in o(\omega(s,t))$ as $t \searrow s$ for all $s$ gives us the definition of \emph{almost (geometric) rough path} and space of these denoted with $\widetilde{\mathscr C}$ (and $\widetilde{\mathscr G}$ for geometric rough paths); this definition is motivated by the fact that the $\varepsilon_{st}$'s vanish in the limit of a sum over a sequence of partitions:
\begin{equation}\label{eq:varepsilon}
\sum_{[s,t] \in \pi} \varepsilon_{st} = \sum_{[s,t] \in \pi}  \frac{\varepsilon_{st}}{\omega(s,t)} \omega(s,t) \leq \omega(0,T) \sup_{[s,t] \in \pi}\frac{\varepsilon_{st}}{\omega(s,t)} \xrightarrow{\abs{\pi} \to 0 } 0 
\end{equation}
since $p < 3$ and $O(\omega(s,t)^{3/p}) \subseteq o(\omega(s,t))$. The same reasoning is also at the root of the following lemma \cite[Theorem 3.3.1]{L98} \cite[Proposition 3.5]{CDLL16}. We write $\approx$ for equality up to an $\varepsilon_{st} \in o(\omega(s,t))$ as $t \searrow s$.
\begin{lem} \label{lem:almost}
	\begin{enumerate}
		\item If $\bfX, \bfY \in \mathscr C^p_\omega([0,T],\bbR^d)$, $\bfX \approx \bfY \Rightarrow \bfX = \bfY$;
		\item Given $\widetilde\bfX \in \mathscr C^p_\omega([0,T],\bbR^d)$, there exists a unique $\bfX \in \mathscr C^p_\omega([0,T],\bbR^d)$ with $\bfX \approx \widetilde {\bfX}$, which is given by
		\begin{equation}
		\bfX_{st} = \lim_{n \to \infty} \bigotimes_{[u,v] \in \pi_n} \widetilde{\bfX}_{uv}  
		\end{equation}
		where $\pi_n$ is any sequence of partitions of $[s,t]$ with vanishing step size. Moreover, if $\widetilde \bfX \in \widetilde{\mathscr G}^p_\omega([0,T],\bbR^d)$, $\bfX \in \mathscr G^p_\omega([0,T],\bbR^d)$.
	\end{enumerate}
	Both statements also hold when restricted to the level of paths $\in \mathcal C^p_\omega([0,T],\bbR^d)$.
\end{lem}
If $\bfX$ and $\widetilde\bfX$ are related as in \autoref[2.]{lem:almost} we will say that latter is the rough path \emph{associated} to the former.

Given $\bfX \in \mathscr C^p_\omega([0,T],\bbR^d)$ we may associate a canonical element ${_\text{g}\hspace{-0.1em}}\bfX \in \mathscr G^p_\omega([0,T],\bbR^d)$, which we call its \emph{geometrisation}, with trace equal to that of $\bfX$ and
\begin{equation}
{_\text{g}\hspace{-0.1em}}\bbX^{\alpha\beta}_{st} \coloneqq \tfrac 12 (\bbX^{\alpha\beta}_{st} - \bbX^{\beta\alpha}_{st}) + \tfrac 12 X^\alpha_{st}X^\beta_{st}
\end{equation}
In other words, ${_\text{g}\hspace{-0.1em}}\bbX$ has the same antisymmetric part as $\bbX$ and symmetric part fixed by the trace and the geometricity condition, and it is easily checked that the Chen identity continues to hold.

\begin{rem}[Inhomogeneous degrees of regularity]\label{rem:inhom}
	Let $d = d_1 + d_2$, and assume that we actually have that $(X^\gamma)^{\gamma = d_1+1,\ldots,d} \in \mathcal C_\omega^{p/2}([0,T],\bbR^{d_2})$. Then, for rough paths with trace $X$ we will only require for the second order components $\bbX^{\alpha\beta}$ with $\alpha, \beta = 1,\ldots,d_1$ to be defined; we will denote the set of these $\mathscr C^{p,p/2}_\omega([0,T],\bbR^{d_1,d_2})$. This is because (according to the interpretation \autoref{eq:bbXint}) all other components would have to be defined as having regularity $O(\omega(s,t)^{3/p})$ and thus be negligible in all expressions involving them. Note that results in the literature which are not explicitly stated for inhomogeneous degrees of regularity may be inferred for such rough paths by simply defining the missing second order components by Young integration (although sharper results exist when keeping track of the exact regularities of the components). In other words, everything that we shall prove for $\bfX \in \mathscr C_\omega^p([0,T],\bbR^d)$ also holds when $\bfX \in \mathscr C^{p,p/2}_\omega([0,T],\bbR^{d_1,d_2})$. The most important cases of a rough path with inhomogeneous degrees of regularity are given by, for $\bfX \in \mathscr C_\omega^p([0,T],\bbR^d)$, $(\bfX, [\bfX])\in \mathscr C^{p,p/2}_\omega([0,T],\bbR^{d,d^2})$ and $({_\text{g}}\!\bfX, [\bfX])\in \mathscr G^{p,p/2}_\omega([0,T],\bbR^{d,d^2})$ (or $\bbR^{d,d(d+1)/2}$ if we view $[\bfX]$ as having only $d(d+1)/2$ components, with the rest determined by symmetry; however we will usually be summing over all $d^2$ components of the bracket). For more on rough paths with inhomogeneous degrees of smoothness (in the geometric context) see \cite{Gyu}.
\end{rem}
\subsection{Controlled paths and rough integration}\label{subsec:contr}
We proceed to define the objects which are, in some sense, dual to rough paths, and are original to \cite{Gub04}:
\begin{defn}[Controlled path] \label{def:controlled}
	Let $X \in \mathcal C^p_\omega([0,T], \bbR^d)$. An $\bbR^e$-valued, \emph{$X$-controlled path}, or element of $\mathscr D_X^p ([0,T],\bbR^e)$ is a pair $\bfH = (H,H')$, where $H \in \mathcal C^p_\omega([0,T], \bbR^e)$ (the \emph{trace}), $H' \in \mathcal C^p_\omega([0,T], \bbR^{e \times d})$ (the \emph{Gubinelli derivative} of $H$ w.r.t.\ $X$), and
	\begin{equation} \label{eq:resto}
	R^k_{st} \coloneqq H_{st}^k - H'^k_{\gamma;s} X^\gamma_{st}, \quad \sup_{0 \leq s < t \leq T} \frac{\abs[0]{R_{st}}}{\omega(s,t)^{2/p}} < \infty, 
	\end{equation}
\end{defn}
Here $\bbR^{e \times d}$ should be thought of as $\mathcal L(\bbR^d,\bbR^e)$ (where $\mathcal L$ means \say{linear maps}). We will identify $\bbR^n$-valued expressions with their coordinate expression throughout this paper, e.g.\ we will write $X = (X^\gamma)$, $\bfX = (X^\gamma,\bbX^{\alpha\beta})$, $\bfH = (H^k,H^k_{\gamma})$.  We will use $\approx_2$ as a shorthand for equality up to $O(\omega(s,t)^{2/p})$, i.e.\ \autoref{eq:resto} may be written as $H^k_{st} \approx_2 H'^k_\gamma X^\gamma_{st}$.

The following definition and theorem establishes that rough paths should be thought of as integrators, and their controlled paths as a class of admissible integrands.
\begin{defthm}[Rough integral]\label{def:roughInt}
	Let $\bfX \in \mathscr C^p_\omega([0,T],\bbR^d)$ and $\bfH \in \mathscr D_{X}(\bbR^{e \times d})$. We then define, for $0 \leq s \leq t \leq T$
	\begin{equation}\label{eq:compR}
	\int_s^t \bfH \dif \bfX \coloneqq \lim_{n \to \infty}\sum_{[u,v] \in \pi_n} H_{\gamma;u} X_{uv}^\gamma + H_{\alpha\beta;u}' \bbX_{uv}^{\alpha\beta}
	\end{equation}
	where $(\pi_n)_n$ is a sequence of partitions on $[s,t]$ with vanishing step size. \emph{This limit exists, is independent of such sequence and is obtained by applying \autoref[1., restricted to the path level]{lem:almost} to
		\begin{equation}\label{eq:approxInt}
		H_{\gamma;s} X_{st}^\gamma + H_{\alpha\beta;s}' \bbX_{st}^{\alpha\beta}
		\end{equation}
	}
\end{defthm}
Here $H_t$ is an $\bbR^{e\times d}$-valued path and $H'_t$ is a $\bbR^{e \times d \times d}$-valued path, with superscripts denoting $\bbR^e$-coordinates and subscripts denoting $\bbR^d$-coordinates; in $H'^k_{\alpha\beta}$ the coordinate of the Gubinelli derivative is $\alpha$, i.e.\ the controlled path property now reads $H^k_{\beta;st} - H'^k_{\alpha\beta;s}X^\alpha_{st} \in O(\omega(s,t)^{2/p})$. We will often refer to controlled paths with trace valued in $\bbR^{e \times d}$ as \emph{controlled integrands}. Clearly if $\widetilde \bfX \in \widetilde{\mathscr C}^p_\omega([0,T],\bbR^d)$ we may have substituted it for $\bfX$ in \autoref{eq:compR} and \autoref{eq:approxInt}. We will often omit the integration extrema: in this case identities are to be intended to hold when the integral is taken on any interval. Also notice that it is obvious from the definition that the integral is linear in the integrand and additive on consecutive time intervals.

The condition of $H$ admitting a Gubinelli derivative w.r.t.\ $X$ is a strong condition, and one can only expect it to be satisfied when $H$ bears a special relationship with $X$. One may also ask whether there are conditions on $X$ under which any Gubinelli derivative $H'$ is unique: this is not always true, since if $X$ is too regular inside $\mathcal C^p_\omega([0,T],\bbR^d)$ the regularity requirement on $H'$ becomes less stringent. A condition on $X$ that rules this out, and guarantees uniqueness of the Gubinelli derivative is given by \emph{true roughness} of $X$: this means that for all $s$ in a dense set of $[0,T]$ and for all $\gamma = 1,\ldots,d$
\begin{equation}
\limsup_{t \searrow s} \frac{\abs[0]{X_{st}^\gamma}}{\omega(s,t)^{2/p}} = + \infty
\end{equation}
It is satisfied, for instance, by a.a.\ sample paths of fractional Brownian motion with Hurst parameter $1/3 < H \leq 1/2$, when considered as elements of $\mathcal C^p$, $1/H < p < 3$.
\begin{thm}[Uniqueness of the Gubinelli derivative]\label{thm:gubUniq}
	Let $\bfX \in \mathscr C^p_\omega([0,T],\bbR^d)$ with trace $X$ truly rough, $(H,{^1\!}H'), (H, {^2\!}H') \in \mathscr D_X(\bbR^d)$. Then ${^1\!}H' = {^2\!}H'$.
\end{thm}
A corollary of this result is the uniqueness of the decomposition of the sum of a Young integral and a rough integral:
\begin{thm}[Doob-Meyer for rough paths]\label{thm:doobMey}
	Let $\bfX \in \mathscr C^p_\omega([0,T],\bbR^d)$ with trace $X$ truly rough, $Y \in \mathcal C^{p/2}([0,T],\bbR^d)$, ${^1\!}\bfH,{^2\!}\bfH \in \mathscr D_X(\bbR^{e \times d})$, ${^1\!}K,{^2\!}K \in \mathcal C^p([0,T], \bbR^{e \times d})$ then
	\begin{equation}
	\int {^1\!}\bfH \edif \bfX + \int {^1\!}K \edif Y = \int {^2\!}\bfH \edif \bfX + \int {^2\!}K \edif Y
	\end{equation}
	implies ${^1\!}\bfH = {^2\!}\bfH$ and ${^1\!}K = {^2\!}K$.
\end{thm}
In most cases, as for \autoref{ex:1form} below, the Gubinelli derivative is defined in a canonical manner, and is intended to be computed accordingly, regardless of whether uniqueness holds or not. 

\begin{expl}[Examples of canonically controlled paths]\label{ex:1form}
	\begin{enumerate}
		\item The simplest example of an $X$-controlled path is a smooth function $f \in C^\infty(\bbR^d,\bbR^e)$ applied to $X$: its Gubinelli derivative is given by $Df(X)$ (where $Df \in C^\infty(\bbR^d, \bbR^{e \times d})$ is the differential of $f$, with coordinates $\partial_\gamma f^k$) since  
		\begin{equation}
		f^k(X)_{st} - \partial_\gamma f^k(X_s)X^\gamma_{st} \in O(\abs[0]{X_{st}}^2) \subseteq O(\omega(s,t)^{2/p})
		\end{equation}
		by Taylor's theorem. The regularity requirements on $f(X)$ and $Df(X)$ are satisfied by \autoref{lem:regSmooth}. We call this $X$-controlled path $\boldsymbol{f}(X)$;
		
		\item Let $\bfX, \bfH$ be as in \autoref{def:roughInt}, then the rough integral $\int_0^\cdot \bfH \dif \bfX$ admits Gubinelli derivative $H$. We denote the resulting element of $\mathscr D_X(\bbR^e)$ simply by $\boldsymbol \int \bfH \dif \bfX$;
		
		\item Assume $\bfH \in \mathscr D_X(\bbR^e)$ and that $K \in \mathcal C^{p/2}_\omega([0,T],\bbR^e)$, then we may use $H'$ as the Gubinelli derivative of $H + K$ and we have that $(H + K,H') \in \mathscr D_X(\bbR^e)$.
	\end{enumerate}
\end{expl}

\begin{expl}[Difference of rough integrals against rough paths with common trace]\label{expl:difference12}
	Let ${_1 \!}\uX = (X,{_1 \!}\bbX), {_2 \!}\uX = (X,{_2 \!}\bbX) \in \mathscr C^p_\omega([0,T],\bbR^d)$, $\oH \in \mathscr D^p_X(\bbR^{e \times d})$. Then it is easy to verify that there must exist a path $D \in \mathcal C^{p/2}_\omega(\bbR^{d \times d})$ s.t.\ ${_2\!}\bbX_{st} = {_1\!} \bbX_{st} + D_{st}$, and it is easily deduced from the \autoref{eq:approxInt}
	\begin{equation}\label{eq:diffD}
	\int \oH \dif {_2\!}\uX = \int \oH \dif {_1\!}\uX + \int H' \dif D
	\end{equation}
	where the second integral on the right is intended in the sense of Young. An important special case is when ${_2 \!}\uX = {_\text{g} \!}\bfX$ for a rough path $\bfX$, in which case $D = \frac 12 [\bfX]$. Note that this identity also holds at the level of controlled paths, since the Gubinelli derivatives (taken according to \autoref[2.,3.]{ex:1form}) both coincide with $H$. We will often use the notation
	\begin{equation}
	\circ \dif \bfX \coloneqq \dif {_\text{g}\!}\bfX
	\end{equation}
	which is motivated by Stratonovich calculus (see \autoref{rem:stochrough} below).
\end{expl}

A controlled path may be transformed into a rough path in a canonical fashion:
\begin{defn}[Lift of a controlled path]
	Let $\bfX \in \mathscr C^p_\omega([0,T],\bbR^d)$, $\bfH \in \mathscr D_X(\bbR^e)$. Define $\uparrow_\bbX \hspace{-0.25em} \bfH$ to be the rough path associated to $\upharpoonleft_\bbX \hspace{-0.25em} \bfH$, defined as
	\begin{equation}
	(\upharpoonleft_\bbX \hspace{-0.25em} \bfH)_{st} \coloneqq (H_{st}^k, H^i_{\alpha;s}H^j_{\beta;s}\bbX^{\alpha\beta}_{st}) 
	\end{equation}
	which is easily verified to belong to $\widetilde{\mathscr{C}}^p_\omega([0,T],\bbR^e)$ (and $\widetilde{\mathscr{G}}^p_\omega([0,T],\bbR^e)$ if $\bfX$ is geometric).
\end{defn}
Note that by \autoref[2.]{lem:almost} this implies that if $\bfX$ is geometric, so is $\uparrow_\bbX \hspace{-0.25em} \bfH$.
\begin{expl}[Lifts of controlled paths] \label{expl:lifts}
	\begin{enumerate}
		\item Given a $f \in C^\infty(\bbR^d, \bbR^e)$ we define $f_*\bfX \coloneqq \uparrow_\bbX \hspace{-0.25em} \boldsymbol{f}(X)$ the \emph{pushforward} of $\bfX$ through $f$, and by Taylor's formula we have
		\begin{equation}
		(f_* \bfX)_{st} \approx (\partial_\gamma f^k(X_s)X^\gamma_{st} + \tfrac 12 \partial_{\alpha\beta} f^k(X_s)X^\alpha_{st}X^\beta_{st}, \partial_\alpha f^i \partial_\beta f^j(X_s)\bbX_{st}^{\alpha\beta})
		\end{equation}
		\item Rough integrals may be lifted to rough paths: if $\bfH$ is as in \autoref{def:roughInt} we abuse the notation once again by setting $\int_s^t \bfH \dif \bfX \coloneqq (\uparrow_\bbX \hspace{-0.25em}\int \bfH \dif \bfX)_{st}$ and we have
		\begin{equation}
		\int_s^t \bfH \dif \bfX \approx (H^k_{\gamma;s}X^\gamma_{st} + H'^k_{\alpha\beta;s}\bbX^{\alpha\beta}_{st}, H^i_{\alpha;s}H^j_{\beta;s}\bbX^{\alpha\beta}_{st})
		\end{equation}
		Note, however, that \autoref{eq:diffD} does not hold at the rough path level, since the lift on the LHS would be computed using ${_2\!}\bbX$, and the one on the RHS using ${_1\!}\bbX$, and this would affect the second order part of the rough integrals. For this reason we will mostly consider It\^o-Stratonovich type corrections only at the trace level.
	\end{enumerate}
\end{expl}
Whenever there is an ambiguity as to whether a function on $\Delta_T$ is a controlled or rough path we will rely on coordinate notation to distinguish these two possibilities, e.g.\ $(\int \bfH \dif \bfX)^k$ are the coordinates of the trace of the controlled/rough path, $(\int \bfH \dif \bfX)_{\gamma}^k = H^k_{\gamma}$ are the coordinates of the Gubinelli derivative of the controlled path and $(\int_s^t \bfH \dif \bfX)^{ij} \approx H^i_{\alpha;s}H^j_{\beta;s}\bbX^{\alpha\beta}_{st}$ those of the second order part of the rough path. We will often use coordinate notation inside the integral too, to track the action of the integrand on the integrator, e.g.\ $(\int \bfH \dif \bfX)^k \eqqcolon \int \bfH^k_\gamma \dif \bfX^\gamma$, with the understanding that we also need the second-order coordinates of $\bfX$ and $\bfH$ to compute this integral (this will help make more complicated expressions clearer). 

\begin{prop}[Operations on controlled paths]\label{prop:contrProp}
	Let $X \in \mathcal C^p_\omega([0,T],\bbR^d)$.
	\begin{description}
		\item[Change of controlling path.] Let $\bfH \in \mathscr D_X(\bbR^e)$, $\boldsymbol K \in \mathscr D_H(\bbR^f)$, then
		\begin{equation}
		\bfK * H' \coloneqq (K^c,K'^c_k H'^k_\gamma) \in \mathscr D_X(\bbR^f)
		\end{equation}
		In particular, if $\bfK = \boldsymbol f(H)$ for $f \in C^\infty(\bbR^e,\bbR^f)$ we denote this $f_* \bfH$ and call it the pushforward of $\bfH$ through $f$;
		\item[Leibniz rule.] Let $\bfH \in \mathscr D_X(\bbR^{f \times e})$ and $\bfK \in \mathscr D_X(\bbR^{g \times f})$, then
		\begin{equation}
		\bfK \cdot \bfH \coloneqq (K^r_c H^c_k, K'^r_{\gamma c} H^c_k + K^r_c H'^c_{\gamma k}) \in \mathscr D_X(\bbR^{g \times e})
		\end{equation}
		\item[Pullback.] Let $g \in C^\infty(\bbR^d, \bbR^e)$, $\bfH \in \mathscr D_{g(X)}(\bbR^{f \times e})$, then
		\begin{equation}
		\begin{split}
		g^*\bfH &\coloneqq (\bfH * \underbrace{Dg(X)}_{= g(X)'}) \cdot \boldsymbol{Dg}(X) \in \mathscr D_X(\bbR^{f \times d}) \\
		&= (H^c_k \partial_\gamma g^k(X), H'^c_{ij} \partial_\alpha g^i \partial_\beta g^j(X) + H^c_k \partial_{\alpha\beta} g^k(X) )
		\end{split}
		\end{equation}
	\end{description}
	\begin{proof}
		Clearly all three paths belong to $\mathcal C^p_\omega$. We need to check that \autoref{eq:resto} holds in all three cases. In the case of the change of controlling path we have
		\begin{equation}
		K^c_{st} - K'^c_{k;s} H'^k_{\gamma;s} X_{st}^\gamma \approx_2 K^c_{st} - K'^c_{k;s}H_{st}^k \approx_2 0
		\end{equation}
		
		As for the Leibniz rule, consider the matrix multiplication function
		\begin{equation}
		\mathscr m \colon \bbR^{g \times f} \times \bbR^{f \times e} \to \bbR^{g \times e}, \quad (z^r_c,y^c_k) \mapsto (z^r_cy^c_k)
		\end{equation}
		It is easily verified that $\bfK \cdot \bfH = \mathscr{m}_*(\bfH,\bfK)$, the pushforward of controlled paths being defined in the step above.
		
		The case of the pullback readily follows from its expression as a combination of the two above constructions.
	\end{proof}
\end{prop}

\begin{prop}[Compatibility]\label{prop:compat}
	Let $\bfX \in \mathscr C^p_\omega([0,T],\bbR^d)$ and $\bfH \in \mathscr D_X(\bbR^e)$:
	\begin{enumerate}
		\item Lifting is compatible with change of controlling path in the sense that, for $\bfK \in \mathscr D_H(\bbR^f)$ we have
		\begin{equation}
		\uparrow_\mathbb{H} \hspace{-0.25em} \bfK = \uparrow_\bbX \hspace{-0.25em} (\bfK * H')
		\end{equation}
		where $\mathbb H$ denotes the second order part of the rough path $\uparrow_\bbX \hspace{-0.25em} \bfH$. In particular, for $f \in C^\infty(\bbR^e,\bbR^f)$ pushforward of rough and controlled paths are related through lift by $f_*\hspace{-0.35em}\uparrow_\bbX \hspace{-0.25em} \bfH = \uparrow_\bbX \hspace{-0.25em} f_*\bfH$. Moreover, $f_*(g_*\bfX) = (f \circ g)_*\bfX$ for appropriately valued smooth maps $f,g$;
		\item Lifting is compatible with geometrisation in the sense that
		\begin{equation}
		{_\emph{g}\hspace{-0.1em}}(\uparrow_\bbX \hspace{-0.25em} \bfH) = \uparrow_{{_\emph{g}\hspace{-0.1em}}\bbX} \hspace{-0.25em} \bfH
		\end{equation}
		In particular, pushforward of rough paths and rough integration preserve geometricity;
		\item For appropriately value smooth maps $f,g$ and controlled integrands $\bfK$ we have
		\begin{equation}
		(f \circ g)^*\bfK = g^*(f^*\bfK)
		\end{equation}
		
	\end{enumerate}
	\begin{proof}
		As for the first claim, the two rough paths agree on the trace $K$ and second order part
		\begin{equation}
		(\uparrow_\mathbb{H} \hspace{-0.25em} \bfK)^{ab}_{st} \approx K'^a_{i;s} K'^b_{j;s} \mathbb{H}^{ij}_{st} \approx K'^a_{i;s} K'^b_{j;s} H'^i_{\alpha;s}H'^j_{\beta;s} \bbX^{\alpha\beta}_{st} = (\uparrow_\bbX \hspace{-0.25em} (\bfK * H'))^{ab} 
		\end{equation}
		Identity of the two rough paths therefore holds by \autoref[1.]{lem:almost}. Now, taking $\bfK = \boldsymbol{f}(H)$ and the definitions of pushforward this yields 
		\begin{equation}
		f_*\hspace{-0.35em}\uparrow_\bbX \hspace{-0.25em} \bfH = \uparrow_\mathbb{H} \hspace{-0.25em} (\boldsymbol f(H)) = \uparrow_\bbX \hspace{-0.25em} f_*\bfH
		\end{equation}
		Taking, furthermore, $\bfH = \boldsymbol{g}(X)$ we obtain
		\begin{equation}
		f_*(g_* \bfX) = f_*\hspace{-0.35em}\uparrow_\bbX \hspace{-0.25em} \boldsymbol{g}(X) = \uparrow_\bbX \hspace{-0.25em} (f_*\boldsymbol{g}(X)) = \uparrow_\bbX \hspace{-0.25em}(\boldsymbol{f \circ g}(X)) = (f \circ g)_*\bfX
		\end{equation}
		
		As for the second claim, the two rough paths have the same trace, and therefore the same symmetric part of the second order part, and antisymmetric part equal to half of
		\begin{equation}
		\begin{split}
		{_\text{g}\hspace{-0.1em}}(\uparrow_\bbX \hspace{-0.25em} \bfH)^{ij} - {_\text{g}\hspace{-0.1em}}(\uparrow_\bbX \hspace{-0.25em} \bfH)^{ji} &= (\uparrow_\bbX \hspace{-0.25em} \bfH)^{ij} - (\uparrow_\bbX \hspace{-0.25em} \bfH)^{ji} \\
		&\approx H'^i_{\alpha;s} H'^j_{\beta;s} (\bbX^{\alpha\beta}_{st} - \bbX^{\beta\alpha}_{st})\\
		&\approx (\uparrow_{{_\text{g}\hspace{-0.1em}}\bbX} \hspace{-0.25em} \bfH)^{ij} - (\uparrow_{{_\text{g}\hspace{-0.1em}}\bbX} \hspace{-0.25em} \bfH)^{ji}
		\end{split}
		\end{equation}
		Therefore ${_\text{g}\hspace{-0.1em}}(\uparrow_\bbX \hspace{-0.25em} \bfH) \approx \uparrow_{{_\text{g}\hspace{-0.1em}}\bbX} \hspace{-0.25em} \bfH$ and we conclude again by \autoref[1.]{lem:almost}.
		
		The final statement is verified using a similar comparison of the expressions in coordinates.
	\end{proof}
\end{prop}

The following is a rough path version of the It\^o lemma. Note how the formula simplifies to a first order chain rule in the case of $\bfX$ geometric. It is followed by the rough path-version of the Kunita-Watanabe identity, where the bracket path takes the role of quadratic covariation matrix.
\begin{thm}[It\^o lemma for rough paths]\label{thm:ito}
	Let $\bfX \in \mathscr C^p_\omega([0,T],\bbR^d)$ and $f \in C^\infty(\bbR^d, \bbR^e)$. Then
	\begin{equation}
	f(X) = f(X_0) + \int_0^\cdot \boldsymbol{Df}(X) \edif \bfX + \frac 12 \int_0^\cdot D^2f(X) \edif [\bfX]
	\end{equation}
	Moreover, the Gubinelli derivatives of the LHS and RHS, computed canonically according to \autoref{ex:1form} agree, thus giving rise to an identity in $\mathscr D_X(\bbR^e)$, and after applying $\uparrow_\bbX$, to one in $\mathscr C^p_\omega([0,T],\bbR^e)$ (with the term $f(X_0)$ only influencing the trace).
	\begin{proof}
		The path-level statement is proved in \cite[Proposition 5.6]{FH14}. The Gubinelli derivative of the LHS according to \autoref[1.]{ex:1form} is $Df(X)$, which coincides with the Gubinelli derivative of the LHS according to \autoref[2.,3.]{ex:1form} (since the bracket path, and thus the Young integral has higher regularity).
	\end{proof}
\end{thm}

Note how the integral of the exact 1-form $Df(X)$ does not require the whole of $\bfX$: this is because its Gubinelli derivative, $D^2f(X)$, is symmetric. Only the symmetric part of $\bbX$ is needed: the pair $(X,\odot \bbX)$ (with $\odot$ denoting the symmetrisation operator) is called a \emph{reduced rough path}. 

\begin{prop}[Kunita-Watanabe identity for rough paths]\label{prop:KW}
	Let $\bfX \in \mathscr C^p_\omega([0,T],\bbR^d)$, $\oH \in \mathscr D_X(\bbR^e)$. Then
	\begin{equation}
	[\uparrow_\bbX \hspace{-0.25em} \oH ]^{ij}_{st} \approx H_{st}^i H_{st}^j - H'^i_{\alpha;s} H'^j_{\beta;s}(\bbX_{st}^{\alpha\beta} + \bbX_{st}^{\beta\alpha})
	\end{equation}
	so in particular (if $e = f \times d$ in the second case below)
	\begin{equation}
	[f_* \uX ]^{ij}_{st} = \int_s^t \partial_\alpha f^i(X)\partial_\beta f^j(X) \edif [\uX]^{\alpha\beta},\quad \bigg[ \int \bfH \edif \bfX  \bigg]^{ij}_{st} = \int_s^t H^i_\alpha H^j_\beta \edif [\bfX]^{\alpha\beta}
	\end{equation}
	\begin{proof}
		The first claim is immediate from \autoref{eq:bracket}, \autoref{expl:lifts}. The bracket of a pushforward is computed as
		\begin{align}
		\begin{split}
		[f_*\uX]_{st}^{ij} &\approx f^i(X)_{st}f^j(X)_{st} -\partial_\alpha f^i(X_s)\partial_\beta f^j(X_s) (\bbX^{\alpha\beta}_{st} + \bbX^{\beta\alpha}_{st})\\
		&\approx \partial_\alpha f^i(X_s)\partial_\beta f^j(X_s) (X^\alpha_{st} X^\beta_{st} - (\bbX^{\alpha\beta}_{st} + \bbX^{\beta\alpha}_{st}))\\
		&= \partial_\alpha f^i(X_s)\partial_\beta f^j(X_s) [\uX]_{st}^{\alpha\beta}\\
		&\approx \int \partial_i f(X)\partial_j f(X) \dif [\uX]^{ij}
		\end{split}
		\end{align}
		and since the integral is additive on consecutive intervals we conclude that we have equality by uniqueness in \autoref[1.]{lem:almost}. As for the rough integral
		\begin{align}
		\begin{split}
		\bigg[ \int \bfH \dif \bfX  \bigg]_{st}^{ij} &\approx (H^i_{\gamma; s}X^\gamma_{st} + H'^i_{\alpha\beta;s} \bbX_{st}^{\alpha\beta})(H^j_{\gamma; s}X^\gamma_{st} + H'^j_{\alpha\beta;s} \bbX_{st}^{\alpha\beta}) - H^i_{\alpha;s}H^j_{\beta;s}\bbX_{st}^{\alpha\beta}  \\
		&\approx H^i_{\alpha;s}H^j_{\beta;s} (X_{st}^\alpha X_{st}^\beta - (\bbX^{\alpha\beta}_{st} +\bbX^{\beta\alpha}_{st}))  \\
		&= H^i_{\alpha;s}H^j_{\beta;s} [\uX]^{\alpha\beta}_{st} \\
		&\approx \int_s^t H^i_{\alpha}H^j_{\beta} \dif [\bfX]^{\alpha\beta}
		\end{split}
		\end{align}
		and conclude as before that equality holds.
	\end{proof}
\end{prop}

The fact that the rough integral can be canonically considered a rough path in its own right naturally leads to the question of associativity, which is answered in the affirmative:
\begin{thm}[Associativity of the rough integral]\label{thm:assoc}
	Let $\bfX \in \mathscr C^p_\omega([0,T],\bbR^d)$, $\bfH \in \mathscr D_X(\bbR^{e \times d})$, $\boldsymbol I \coloneqq \int \bfH \edif \bfX \in \mathscr D_X(\bbR^e)$, $\bfY \coloneqq \uparrow_\bbX \hspace{-0.25em} \boldsymbol I$, $\bfK \in \mathscr D_I(\bbR^{f \times e})$. Then
	\begin{equation}
	\bigg(\int \bfK \edif \bfY \bigg) * I' =  \int (\bfK * I') \cdot \bfH \edif \bfX \in \mathscr D_X(\bbR^f)
	\end{equation}
	As a result, the identity $\int \bfK \edif \bfY =  \int (\bfK * H) \cdot \bfH \edif \bfX$ also holds in $\mathscr C^p_\omega([0,T],\bbR^f)$.
	\begin{proof}
		At the level of the trace we have
		\begin{equation}
		\begin{split}
		\int_s^t \bfK^c_k \dif \bfY^k &\approx K^c_{k;s}Y^k_{st} + K'^c_{ij;s} \bbY^{ij}_{st} \\
		&\approx K^c_{k;s}(H^k_{\gamma;s}X^\gamma_{st} + H'^k_{\alpha\beta;s}\bbX^{\alpha\beta}_{st}) + K'^c_{ij;s} H^i_{\alpha;s} H^j_{\beta;s} \bbX^{\alpha\beta}_{st} \\
		&\approx \int_s^t (\bfK * H) \cdot \bfH \dif \bfX
		\end{split}
		\end{equation}
		which proves the identity of the traces, since $I' = H$. Their Gubinelli derivatives w.r.t.\ $X$, as computed according to \autoref[2.]{ex:1form} and \autoref{prop:contrProp} both coincide with $(K^c_k H^k_\gamma)$. Passing to the lift on this identity we have
		\begin{equation}
		\bigg\uparrow_\bbX  \int (\bfK * H) \cdot \bfH \dif \bfX  = \bigg\uparrow_\bbX  \int \bfK \dif \bfY * H  = \bigg\uparrow_\mathbb{\bbY} \int \bfK \dif \bfY  
		\end{equation}
		where we have used \autoref[1.]{prop:compat} in the second identity. This is the identity required in the second statement.
	\end{proof}
\end{thm}
The next proposition expresses the degree to which pushforward of rough paths and pullback of controlled paths fail to be adjoint operators under the rough integral pairing; in particular the adjunction does hold when the integrator is geometric or when $g$ is an affine map.
\begin{cor}\label{cor:pushpull}
	Let $\bfX, \bfH, g$ be as in \autoref[Pullback.]{prop:contrProp}. Then
	\begin{equation}
	\bigg(\int \bfH \edif (g_*\bfX)\bigg) * Dg(X) = \int g^* \bfH \edif \bfX + \frac 12 \int H \cdot D^2g(X) \edif [\uX]
	\end{equation}
	where, as usual, the identity holds in $\mathscr D_X(\bbR^e)$ according to \autoref{ex:1form} and thus in $\mathscr C^p_\omega([0,T],\bbR^e)$. 
	\begin{proof}
		Plugging in the the expression for $g_*\bfX$ given by \autoref{thm:ito} and applying \autoref{thm:assoc} we have
		\begin{equation}
		\begin{split}
		\bigg(\int \bfH \dif (g_*\bfX) \bigg) * Dg(X) &= \bigg(\int \bfH \dif ({\textstyle\int} \boldsymbol{Dg}(X) \dif \uX + {\textstyle\frac 12 \int D^2g(X) \dif [\bfX]}) \bigg) * Dg(X) \\
		&= \int (\bfH * Dg(X)) \cdot \boldsymbol{Dg}(X) \dif \bfX  + \frac 12 \int H \cdot D^2g(X) \dif [\bfX] \\
		&= \int g^* \bfH \dif \bfX + \frac 12 \int H \cdot D^2g(X) \dif [\uX]
		\end{split}
		\end{equation}
		As usual, the more regular Young integral only contributes to the trace of the $X$-controlled/rough paths in question.
	\end{proof}
\end{cor} 

\subsection{Rough differential equations}\label{subsec:RDEs}
We proceed to discuss a central theme of rough path theory: that of rough differential equations, or RDEs.
\begin{defn}\label{def:rde}
	Let $\bfX \in \mathscr C^p_\omega([0,T],\bbR^d)$, $F \in C^\infty(\bbR^{e+d},\bbR^{e \times d})$. A \emph{controlled solution} to the RDE
	\begin{equation}\label{eq:rde}
	\dif \bfY = F(Y,X) \dif \bfX, \quad Y_0 = y_0
	\end{equation}
	(which we will write in coordinates as $\dif \bfY^k = F^k_\gamma(Y,X) \dif \bfX^\gamma$ when we wish to emphasise the action of the field of linear maps $F$ on the driver $\bfX$) is an element $\bfY \in \mathscr D_X(\bbR^e)$ s.t.\
	\begin{equation}\label{eq:intRde}
	\bfY = y_0 + \int F_* (\bfY,\bfX) \dif \bfX \in \mathscr D_X(\bbR^e)
	\end{equation}
	where $F_*(\bfY,\bfX)$ is the pushforward of the $\bbR^{e+d}$-valued $X$-controlled path with trace $(Y,X)$ and Gubinelli derivative $(Y',\mathbbm 1)$. We will call $\uparrow_\bbX \hspace{-0.25em}\bfY$ (which we will denote again $\bfY$) a \emph{rough path solution} to \autoref{eq:rde}.
\end{defn}
We will sometimes write $\dif Y$ (without the bold font for $Y$) on the LHS of \autoref{eq:rde} when only referring to the trace level of the solution. Note that the definition of controlled solution implies the requirement $Y' = F(Y,X)$ and 
\begin{equation}
F_*(\bfY,\bfX) = (F^k_\gamma(Y,X),\partial_\alpha F^k_\beta(Y,X) + F_\alpha^h \partial_h F^k_\beta(Y,X))
\end{equation}
Since a solution of either type is entirely determined by its trace and $F$, $\bfX$ we will often just use the term \emph{solution} without specifying which type we intend.
\begin{rem}\label{rem:RDEredu}
	Usually only RDEs of the form $\dif \bfY = F(Y) \dif \bfX$ are considered. \autoref{eq:rde} can be considered as a special case of this by simply \say{doubling the variables}, i.e.\ considering the joint RDE
	\begin{equation}\label{eq:doubling}
	\dif \begin{pmatrix}
	\bfX \\ \bfY
	\end{pmatrix} = \begin{pmatrix}
	\mathbbm 1 \\ F(Y,X)
	\end{pmatrix} \dif \bfX
	\end{equation}
	We have chosen to consider RDEs that also depend on $X$ since this will become a compulsory requirement when $\bfX$ is manifold-valued; this framework is taken from \cite{E89}, where it is used in the context of manifolds-valued SDEs. A side benefit of introducing this dependence is that rough integrals may now be seen as a particular case of RDEs, namely when $F$ is independent of $Y$. In particular, thanks to \autoref{thm:ito}, controlled paths given by 1-forms may be viewed as particular cases of solutions to RDEs driven by the rough path $(\bfX, [\bfX])$.
\end{rem}

\begin{expl}[RDEs driven by rough paths with common trace]\label{expl:rdediff}
	Let ${_1\!}\bfX, {_2\!}\bfX,D$ be as in \autoref{expl:difference12}. Then we have the following identity of controlled solutions
	\begin{equation}
	\dif Y^k = F^k_\gamma(Y,X) \dif {_2\!}\bfX^\gamma \Longleftrightarrow \dif Y^k = F^k_\gamma(Y,X) \dif {_1\!}\bfX^k + (\partial_\alpha F^k_\beta + F_\alpha^h \partial_h F^k_\beta)(Y,X) \dif D^{\alpha\beta}
	\end{equation}
	Note this identity does not hold for rough path solutions, for the reason provided in \autoref[2.]{expl:lifts}. The second expression is an RDE driven by the rough path with trace $(X,D)$ and second order part $\mathbb X$ (since $D^{\alpha\beta}_{st} \in O(\omega(s,t)^{p/2})$ we do not need other components by \autoref{rem:inhom}, i.e.\ the integral against $D$ is intended in the sense of Young). This is particularly important when ${_2 \!}\uX = {_\text{g} \!}\bfX$, $D = \frac 12 [\bfX]$ for a rough path $\bfX$, as it informs us that every RDE may be rewritten as an RDE driven by the geometric rough path $({_\text{g}\!}\bfX,[\bfX])$.
\end{expl}

The following theorem is proved in \cite[Corollary 2.17, Theorem 4.2]{CDL15}, and its proof carries over to the case of $\bfX$ non-geometric (thanks to \autoref{expl:rdediff}) and with $F$ depending on $X$ (thanks to \autoref{rem:RDEredu}). We will say that $\bfY$ is a controlled/rough path solution \emph{up to time $S \leq T$} if it is a solution to \autoref{eq:rde} where the driving rough path is substituted with $\bfX|_{[0,R]} \in \mathscr C^p_\omega([0,R],\bbR^d)$, for all $R < S$. Note that, according to this terminology, a solution up to time $T$ is a not necessarily a solution on the whole of $[0,T]$ (the former may explode precisely at time $T$, while for the latter we have $Y_T = y_0 + \int_0^T F_*(\bfY, \bfX)\dif \bfX$): to distinguish the two we will call the latter a \emph{global solution}.
\begin{thm}[Local existence and uniqueness]\label{thm:localE}
	Precisely one of the following two possibility holds w.r.t.\ \autoref{eq:rde}
	\begin{enumerate}
		\item A global solution exists;
		\item There exists an $S \leq T$ and a solution up to time $S$, with $Y_{[0,S)}$ not contained in any compact set of $\bbR^e$.
	\end{enumerate}
	Moreover, in either case, the solution is unique on the interval on which it is defined.
\end{thm}
The following lemma further specifies that the exit time from an open neighbourhood is bounded from below, uniformly in the initial time and initial condition (ranging in a precompact neighbourhood) of an RDE with fixed driver $\bfX$. It can be found in \cite[Corollary 2.17]{CDL15}, and its proof carries over to the setting considered here once again by \autoref{expl:rdediff} and \autoref{rem:RDEredu} (and using the obvious fact that $X_{[0,T]}$ is compact).
\begin{lem}\label{lem:localUnifE}
	Let $U, V \subseteq \bbR^e$ be open with $V \supseteq \overline U$ compact. Then there exists a $\delta > 0$ s.t.\ for all $t_0 \in [0,T]$ and $y_0 \in U$ the unique solution to
	\begin{equation}
	\edif Y = F(Y,X) \edif \bfX, \quad Y_{t_0} = y_0
	\end{equation}
	is defined and satisfies $Y \in V$ on $[t_0,(t_0 + \delta) \wedge T]$.
\end{lem}

\begin{expl}[Operations on RDEs]\label{expl:RDEOps}
	\begin{enumerate}
		\item Let 
		\begin{equation}\label{eq:compRDEs}
		\dif \bfY^k = F^k_\gamma(Y,X) \dif \bfX^\gamma, \quad \dif \bfZ^c = G^c_k(Z,Y) \dif \bfY^k
		\end{equation}
		Then by \autoref{thm:assoc}
		\begin{equation}
		\bfZ = \int (G(Z,Y) * (Z,Y)') \cdot F(Y,X) \dif \bfX
		\end{equation}
		Where the Gubinelli derivative of $(Z,Y)$ w.r.t.\ $X$ is computed as
		\begin{equation}
		(Z,Y)' = (G(Z,Y)Y', Y') = (G(Y,Z)F(Y,X), F(Y,X))
		\end{equation}
		$\bfZ$ can thus be viewed as the solution, jointly with $\bfY$, to an RDE driven by $\bfX$, and we will simply write
		\begin{equation}
		\dif \bfZ^c = G^c_k(Z,Y) F^k_\gamma(Y,X) \dif \bfX
		\end{equation}
		\item If $\bfY$ is as in \autoref{eq:compRDEs} we may use \autoref{prop:KW} to show that its bracket solves
		\begin{equation}
		\dif[\bfY]^{ij} = F^i_\alpha F^j_\beta(Y,X) \dif [\bfX]^{\alpha\beta}
		\end{equation}
		This can be viewed as the solution, jointly with $\bfY$, to an RDE driven by $(\bfX,[\bfX])$, which by \autoref{expl:rdediff} can be transformed into one driven by $({_\text{g}\!}\bfX,[\bfX])$.
	\end{enumerate} 
	From now on we will manipulate RDEs formally as It\^o SDEs, in coordinates, with the understanding that their precise meaning can be justified as done above.
\end{expl}

Although this is not a paper on global existence, we will need the following lemma that guarantees it in an important special case.
\begin{lem}\label{lem:globE}
	Let $V$ be as in \autoref{def:rde} with
	\begin{equation}
	F^k_\gamma(y,x) = A^k_{\gamma h}(x)y^h + b^k_\gamma(x)
	\end{equation}
	for some $A \in C^\infty(\bbR^d,\bbR^{e \times e \times d})$, $b \in C^\infty(\bbR^d,\bbR^{e \times d})$. Then \autoref{eq:rde} admits a global solution.
	\begin{proof}
		First of all, observe that
		\begin{equation}
		(\partial_\alpha F^k_\beta + F_\alpha^h \partial_h F^k_\beta)(y,x) = (\partial_\alpha A^k_{\beta h} + A^l_{ \alpha h}A^k_{\beta l})(x)y^h + (\partial_\alpha b^k_\beta + b^h_\alpha A^k_{h \beta})(x)
		\end{equation}
		has the same form as $F$: by \autoref{expl:rdediff} we may therefore assume $\bfX$ is geometric. To prove the result it does not suffice to invoke the well-known existence of global solutions for linear RDEs, as $F$ is not linear in $x$. However, we may assume $A$ and $b$ to be bounded with bounded derivatives of all orders, since we only require the values of $A$ and $b$ on the compact set $X_{[0,T]}$ (and may thus multiply them by a mollifier on $\bbR^d$ that vanishes outside an open set containing $X_{[0,T]}$). Now \cite[Theorem 10.53]{FV10} may be applied to \autoref{eq:doubling}, with the only caveat that we have to replace $v$ with $\max_{\gamma}(\norm[0]{A_\gamma}_\infty, \norm[0]{b_\gamma}_\infty)$ in its proof (and in that of \cite[Lemma 10.52]{FV10}).\todo[inline,backgroundcolor=yellow]{This kind of referencing is obviously not ideal, but I can't find an off-the-shelf global $\exists$ result that can be applied to \autoref{eq:doubling} (e.g. Weidner's criterion fails). This is because none of them can take into account the fact that $A$ and $b$ are only evaluated on a compact. The alternative would be to reprove the whole of \cite{FV10} introducing dependencies on $x$ from the start, which I would rather not do.}
	\end{proof}
\end{lem}

\subsection{Stochastic rough paths}\label{subsec:stochRough}
Finally, we address the topic of stochastic processes lifted to rough paths. We denote with $\mathcal S(\Omega,[0,T],\bbR^d)$ the set of $\bbR^d$-valued continuous adapted semimartingales defined up to time $T$ on some stochastic setup $(\Omega, (\mathcal F_t)_{t \in [0,T]}, \mathbb P)$ satisfying the usual conditions. We may define $\bfX, \widehat \bfX \in \mathscr C^p([0,T],\bbR^d)$ a.s.\ by Stratonovich and It\^o integration respectively
\begin{equation}\label{eq:ItoStrat}
\bbX_{st}^{\alpha\beta} \coloneqq \int_s^t X_{su}^\alpha \circ \dif X_u^\beta, \quad \widehat \bbX_{st}^{\alpha\beta} \coloneqq \int_s^t X_{su}^\alpha \dif X_u^\beta, \quad \bbX_{st}^{\alpha\beta} = \widehat \bbX_{st}^{\alpha\beta} + \tfrac 12 [X]_{st}^{\alpha\beta}
\end{equation}
where $[X]$ denotes the quadratic covariation tensor of $X$.

\begin{rem}\label{rem:stochrough}
	We have $[\widehat{\bfX}] = [X]$ and $[\bfX] = 0$ a.s.\, so ${_\text{g}\!}\widehat\uX = \uX$. In general, rough path theory applied to semimartingales extends the usual stochastic calculus, i.e.\ It\^o/Stratonovich stochastic integrals agree a.s.\ with the path-by-path computed rough integrals w.r.t.\ the It\^o/Stratonovich lifts \cite[Proposition 5.1, Corollary 5.2]{FH14}, and the strong solution to an It\^o/Stratonovich SDE coincides a.s.\ with the path-by-path computed solution to the RDE driven by the It\^o/Stratonovich-enhanced rough path \cite[Theorem 9.1]{FH14} (these results are only shown for Brownian integrators, but may be extended to general continuous semimartingales, e.g.\ by reducing to the Brownian case by splitting the integrator into its bounded variation and local martingale parts and applying the Dubins-Schwarz theorem to the latter).
\end{rem}
The following is a statement made in the same spirit which will be important later on.
\begin{prop}\label{prop:sound}
	Let $X \in \mathcal S(\Omega,[0,T])$ and $f \in C^\infty(\bbR^d,\bbR^e)$. Then $f_*\uX$ and $f_* \widehat \uX$ coincide a.s.\ with the lifts of the semimartingale $f(X)$ computed respectively through Stratonovich and It\^o integration.\todo{I tried to show this in the more general case in which $f(X)$ is replaced with an a.s.\ $X$-controlled semim. and the pushforwards with $\uparrow_{\bbX}$ and $\uparrow_{\widehat\bbX}$, but wasn't able to. Not really essential but would have been nice to have. UPDATE: I think it can be done if we have the second Gubinelli derivative $H''$.}
	\begin{proof}
		We begin with the It\^o case. By the classical It\^o formula and \autoref{rem:stochrough} we have that, a.s.\
		\begin{align}
		\begin{split}
		\int_s^t f^i(X) \dif f^j(X) &= \int_s^t f^i(X) \partial_\gamma f^j(X) \dif X + \frac 12 \int_s^t f^i(X) \partial_{\alpha\beta}f(X) \dif [X]^{\alpha\beta} \\
		&= \int_s^t \boldsymbol f^i \boldsymbol \partial_\gamma \boldsymbol f^j(X) \dif \widehat \bfX^\gamma + \frac 12 \int_s^t f^i \partial_{\alpha\beta}f(X) \dif [\widehat \bfX]^{\alpha\beta} \\
		&\approx f^i \partial_\gamma f^j(X_s)X_{st}^\gamma + (\partial_\alpha f^i \partial_\beta f^j + f^i \partial_{\alpha\beta}f^j)(X_s)\widehat\bbX_{st}^{\alpha\beta} \\
		&\mathrel{\phantom{=}} + \tfrac 12 f^i \partial_{\alpha\beta} f^j(X_s) [\widehat\uX]_{st}^{\alpha\beta} \\
		&= f^i \partial_\gamma f^j(X_s)X_{st}^\gamma + \partial_\alpha f^i \partial_\beta f^j\widehat\bbX^{\alpha\beta}_{st} + \tfrac 12 f^i \partial_{\alpha\beta} f^j(X_s) X^\alpha_{st} X^\beta_{st} \\
		&\approx f^i(X_s) f^j(X)_{st} + \partial_\alpha f^i \partial_\beta f^j\widehat\bbX^{\alpha\beta}_{st} 
		\end{split}
		\end{align}
		Therefore a.s.\
		\begin{equation}
		\int_s^t f^i(X)_{su} \dif f^j(X_u) = \int_s^t f^i(X) \dif f^j(X) - f^i(X_s) f^j(X)_{st} \approx \partial_\alpha f^i \partial_\beta f^j\widehat\bbX^{\alpha\beta}_{st}
		\end{equation}
		and we conclude by \autoref[1.]{lem:almost}. The Stratonovich case is handled analogously, with the only difference that the first order change of variable formula holds, and that brackets vanish.
	\end{proof}
\end{prop}

For other examples of stochastic rough paths, which include lifts of Gaussian and Markov processes, we refer to \cite[Ch.\ III]{FV10}. Though these rough paths are mostly geometric, examples of non-geometric, non-semimartingale stochastic rough paths also exist in the literature \cite{QiXu}.

\section{Background on differential geometry}\label{sec:backDG}
In this section we review the differential geometry needed in the rest of this paper. We begin by recalling various equivalent notions of connections on manifolds, and proceed to specialise this study to the case where $M$ is a Riemannian submanifold of $\bbR^d$. We follow \cite{L97} and \cite{KobNom96} for classical differential geometry (with an occasional glance at \cite{Nak03} for expressions in local coordinates), and \cite[Ch.\ 8]{L97}, \cite{CDL15}, \cite{ABR18} and \cite{Dr04} for the extrinsic theory.

\subsection{Linear connections}\label{subsec:linearconn}
Let $M$ be a smooth $m$-dimensional manifold, and $\tau M \colon TM \to M$ its tangent bundle; throughout this paper we will identify fibre bundles with their projection. We will denote the tangent map of a smooth map of manifolds $f \colon M \to N$ by $\tau f \colon \tau M \to \tau N$ (a morphism in the category of vector bundles), with map of total spaces $Tf \colon TM \to TN$ (a smooth map), and by $T_xf$ its restriction to the tangent space $T_xM$. In this subsection we review equivalent notions of a connection on a manifold. Given a smooth fibre bundle $\pi \colon E \to M$ we denote with $\Gamma \pi$ its $C^\infty M$-module of sections and $E_A \coloneqq \pi^{-1}(A)$ for $A \subseteq M$, $E_x \coloneqq E_{\{x\}}$ for $x \in M$.
\begin{defn}[Covariant derivative]
	A \emph{linear connection}, or \emph{covariant derivative} on a smooth vector bundle $\pi \colon E \to M$ is a map
	\begin{equation}
	\nabla \colon \Gamma \tau M \times \Gamma \pi \to \Gamma \pi, \quad (U,e) \mapsto \nabla_{U}e, \quad (\nabla_Ue)(x) \eqqcolon \nabla_{U(x)}e
	\end{equation}
	which is linear in both arguments and which satisfies the Leibniz rule $\nabla_{U(x)}(f e) = f(x) \nabla_{U(x)}e + (U(x)f)e(x)$ for $f \in C^\infty M$.
\end{defn}
The notation $\nabla_{U(x)}e$ is justified by the fact that the value of the section $\nabla_{U}e$ only depends on the value of $U(x)$ (and on the value of $e$ on any curve whose tangent vector at $x$ is $U(x)$); in general we will denote vectors based at $x \in M$ as $U(x), V(x),\ldots$, reserving $U,V,\ldots$ for vector fields (or just vectors based at an unspecified point). Covariant derivatives will mainly be considered on the tangent bundle $\tau M$: in this case it is automatically extended to the whole of $\bigoplus_{k,l \in \mathbb N}\tau M^{\otimes k} \otimes \tau^* M^{\otimes l}$ given a few compatibility conditions \cite[Lemma 4.6]{L97}. A linear connection on $\tau M$ is equivalently defined by a \emph{Hessian}, i.e.\ an $\bbR$-linear map
\begin{equation}
\nabla^2 \colon C^\infty M \to \Gamma (\tau^*M^{\otimes 2})
\end{equation}
satisfying $\nabla^2(fg) = f \nabla^2 g + g \nabla^2 f + \dif f \otimes \dif g + \dif g \otimes \dif f$ for all $f, g \in C^\infty M$. Covariant derivatives on the tangent bundle and Hessians are equivalent data, and are related by
\begin{equation}\label{covarianthessian}
\langle \nabla^2 f(x), U(x) \otimes V(x) \rangle = (U(x)V - \nabla_{U(x)}V) f
\end{equation}
where $V$ is any extension of $V(x)$ to a local section and $U(x)V$ denotes composition of vector fields, i.e.\ the differential operator whose action on $f \in C^\infty M$ is given by $U(x)(y \mapsto V(y)f)$. Given a chart $\varphi$ we denote with $\partial_k\varphi(x)$ the basis elements of the tangent space $T_xM$ defined by $\varphi$, and we abbreviate $\partial_k \coloneqq \partial_k\varphi$ if there is no risk of ambiguity. Moreover, we denote with $\Gamma_{ij}^k$ the Christoffel symbols of $\nabla$ w.r.t.\ $\varphi$: this means $\nabla_{\partial_i}\partial_j = \Gamma^k_{ij} \partial_k$, and therefore
\begin{equation}\label{eq:nablaVector}
\nabla_{U}V = (U^h \partial_h V^k + U^iV^j \Gamma^k_{ij})\partial_k
\end{equation}
and if $\omega \in \Gamma \tau^*M$
\begin{equation}\label{eq:nablaCovector}
\nabla_{U} \omega = (U^i \partial_i \omega_j - U^i \omega_k \Gamma^k_{ij}) \dif^j
\end{equation}
where $\dif^k \coloneqq \dif_k\varphi$ are the elements of the dual basis of $\{\partial_k\varphi(x)\}_k$. Given two charts $\varphi$, $\overline \varphi$ defined on overlapping domains, the (non-tensorial) transformation rule of the Christoffel symbols is
\begin{equation}\label{eq:chrChange}
\Gamma^{\overline k}_{\overline i \overline j} = \partial_k^{\overline k} \partial_{\overline i}^i \partial_{\overline j}^j \Gamma^k_{ij} + \partial_{\overline i \overline j}^h \partial_h^{\overline k}
\end{equation}
where overlined indices refer to $\overline \varphi$ and simple indices to $\varphi$, and the $\partial$'s refer to the derivatives of the change of chart, e.g.\ $\partial_{\overline i \overline j}^h(x) \coloneqq \partial_{\overline i \overline j} (\varphi^h \circ \overline \varphi^{-1})(x)$). The Hessian can be written in coordinates as
\begin{equation}\label{eq:hessiancoords}
\nabla^2 f = (\partial_{ij} - \Gamma_{ij}^k \partial_k)f \ \dif \varphi^i \otimes \dif \varphi^j	
\end{equation} 
Those connections whose Hessians are valued in $\Gamma(\tau^*M^{ \odot 2})$, or equivalently with vanishing \emph{torsion} tensor $\langle \mathcal T, U \otimes V \rangle \coloneqq \nabla_U V - \nabla_V U - [U,V]$ are called \emph{torsion-free}. In local coordinates $\mathcal T_{ij}^k = \Gamma_{ij}^k - \Gamma_{ji}^k$. We can associate to any connection $\nabla$ a torsion-free one by ${^\odot \!}\nabla_U V \coloneqq \nabla_U V - \frac 12 \langle \mathcal T, U \otimes V \rangle$ or equivalently by projecting its Hessian onto $T^*M \odot T^*M$: the symmetrised connection will then define the same set of (parametrised) geodesics.

Let $\mathscr g$ be a Riemannian metric, i.e.\ a section in $\Gamma (\tau^* M^{\odot 2})$ which is nowhere vanishing and positive-definite at all points (many, but not all, of the considerations made in this paper about Riemannian metrics can be extended to pseudo-Riemannian ones). A connection $\nabla$ is \emph{metric} w.r.t.\ $\mathscr g$ if $\nabla \mathscr g = 0$, or in local coordinates
\begin{equation}\label{eq:nablaMetric}
\scrg_{ij,k} - \scrg_{hj} \Gamma^h_{ki} - \scrg_{ih} \Gamma^h_{kj} = 0
\end{equation}
where indices after the comma denote partial differentiation in the chosen chart, i.e.\ $\scrg_{ij,k} \coloneqq \partial_k \scrg_{ij}$, and $\scrg_{ij}$ the components of the metric in the same chart ($\scrg^{ij}$ will denote the inverse of $\scrg_{ij}$, i.e.\ $\scrg^{ik}\scrg_{kj} = \delta^i_j$). We will also use indices after a semicolon to denote covariant differentiation, e.g.\ $\scrg_{ij;k} \coloneqq (\nabla\scrg)_{ijk}$. There is precisely one such connection which is also torsion-free, called the \emph{Levi-Civita} connection of $\scrg$, which we denote ${^\mathscr g \!}\nabla$, and its Christoffel symbols are given by
\begin{equation}\label{eq:GammaLC}
{^\scrg \!}\Gamma^k_{ij} = \tfrac 12 \scrg^{kh}(\scrg_{hj,i} + \scrg_{ih,j} - \scrg_{ij,h})
\end{equation}
When on a Riemannian manifold we will sometimes use the \emph{musical isomorphisms} ${}^\flat \colon \tau M \to \tau^* M$ with inverse ${}^\sharp$. In coordinates these are given by performing \say{index gymnastics} w.r.t.\ $\scrg$, i.e.\ $V_i \coloneqq V^\flat_i \coloneqq \scrg_{ij} V^j$, $\omega^i \coloneqq (\omega^\sharp)^i = \omega_j \scrg^{ij}$. Similar raising and lowering of indices will be performed with arbitrary tensors.
\begin{rem}\label{rem:contorsion}
	If $\nabla$ is $\scrg$-metric, it is not true in general that ${^\odot\!}\nabla$ is metric. Denoting $\tensor{\mathcal T}{^k_{ij}}$ the components of the torsion tensor, we have that the difference between $\nabla$ and ${^\scrg\!}\nabla$ is quantified by the \emph{contorsion} tensor
	\begin{equation}\label{eq:contorsion}
	\mathscr K^k_{ij} \coloneqq \tfrac 12 (\tensor{\mathcal T}{^k_{ij}} + \tensor{\mathcal T}{_i^k_j} + \tensor{\mathcal T}{_j^k_i}), \quad \Gamma^k_{ij} - {^\scrg\!}\Gamma^k_{ij} = \mathscr K^k_{ij}
	\end{equation}
	which has symmetric part $\frac 12 (\tensor{\mathcal T}{_i^k_j} + \tensor{\mathcal T}{_j^k_i})$.
\end{rem}
The \emph{curvature tensor} associated to a connection $\nabla$ is 
\begin{equation}
\mathscr R(U,V)W \coloneqq \nabla_{U}\nabla_V W - \nabla_{V}\nabla_U W - \nabla_{[U,V]}W
\end{equation}
where $[U,V]$ denotes the Lie bracket of vector fields, which vanishes if the vectors are given by the local basis sections $\partial_k$ defined by a chart. We denote the coefficients
\begin{equation}\label{eq:Rcoords}
\tensor{\mathscr R}{_{ijk}^h} \coloneqq \langle R(\partial_i, \partial_j)\partial_k, \dif^h \rangle = \Gamma^h_{jk,i} - \Gamma^h_{ik,j} + \Gamma^h_{il}\Gamma^l_{jk} - \Gamma^h_{jl}\Gamma^l_{ik} 
\end{equation}
and warn the reader that the ordering of the indices is not standard in the literature (this convention is, for instance, the one followed by \cite{L97,YI73}). The curvature tensor satisfies the symmetry
\begin{equation}\label{eq:Rij}
\tensor{\mathscr R}{_{ijk}^h} = -\tensor{\mathscr R}{_{jik}^h}
\end{equation} 
Moreover, if $\nabla$ is torsion-free 
\begin{equation}\label{eq:Rcyclic}
\tensor{\mathscr R}{_{ijk}^h} + \tensor{\mathscr R}{_{jki}^h} + \tensor{\mathscr R}{_{kij}^h} = 0
\end{equation}
Moreover, if $\nabla$ is $\scrg$-metric (but not necessarily torsion-free)
\begin{equation}\label{eq:Rkh}
\tensor{\mathscr R}{_{ijkh}} = - \tensor{\mathscr R}{_{ijhk}}
\end{equation}
and finally if $\nabla = {^\scrg \!}\nabla$
\begin{equation}\label{eq:ijhk}
\tensor{\mathscr R}{_{ijkh}} = \tensor{\mathscr R}{_{khij}}
\end{equation}
These symmetries are often stated directly in the Levi-Civita case, but hold under the more general hypotheses stated above, as can be seen from a careful reading of their proof \cite[Proposition 7.4]{L97}. We also recall the definition of \emph{Ricci tensor} a symmetric tensor field defined as a contraction of the curvature tensor, and whose components we still denote (without ambiguity, thanks to the different number of indices) with the symbol $\mathscr R$:
\begin{equation}
\tensor{\mathscr R}{_{ij}} = - \tensor{\mathscr R}{_{ki}^k_j} = -\tensor{\mathscr R}{_{hikj}}\mathscr g^{hk}
\end{equation}

Given a smooth fibre bundle $\pi \colon E \to M$ with typical fibre the smooth n-dimensional manifold $R$ (note: this is not $\bbR$, and in general is not even a vector space), its \emph{vertical bundle} $V\pi$ is the subbundle of $\tau E$ with total space $VE \coloneqq \operatorname{ker} (T \pi \colon TE \to TM)$, and we have $V_{e(x)}E = T_{e(x)} E_x$, i.e.\ elements of the total space of $V\pi$ are vectors tangent to the fibres of $\pi$. Recall that for a smooth map of manifolds $f \in C^\infty(P,Q)$ and a fibre bundle $\rho \colon D \to Q$, we define the \emph{pullback} bundle
\begin{equation}
f^*\rho \colon \{(p,d) \in P \times D \mid f(p) = \rho(d)\} \coloneqq f^*Q \to P, \quad (p,d) \mapsto p
\end{equation}
and there is a bundle map $f^*\rho \to \rho$
\begin{equation}
\begin{tikzcd}
f^*Q \arrow[d,"f^*\rho",swap] \arrow[r,"\text{pr}_2"] &D \arrow[d,"\rho"] \\
P \arrow[r,"f"] &Q
\end{tikzcd}
\end{equation}
The \emph{vertical lift} of $\pi$ is defined as the fibre bundle isomorphism
\begin{equation}
\begin{split}
&\pi^*\pi \to V\pi, \quad E_x \times E_x \ni (e(x),U(x)) \mapsto \mathscr v(e(x)) U(x)\\
&\mathscr v(e(x))U(x)(f \in C^\infty E) \coloneqq \frac{\dif }{\dif t}\bigg|_0 f(e(x) + t U(x)) 
\end{split}
\end{equation}
An \emph{Ehresmann connection} is a vector bundle $\eta \colon H \to E$ which is complementary to $V\pi$, i.e.\ $H \oplus VE = TE$. When $\pi$ is a vector bundle, Ehresmann connections and a covariant derivatives are equivalent by further requiring of the former that, denoting the sum and scalar multiplication map by
\begin{equation}\label{eq:sigmalambda}
\Sigma \colon E \oplus E \to E,\qquad \Lambda_a \colon E \to E,\ a \in \mathbb R
\end{equation}
$T\Sigma$ map the subbundle $\{(\alpha(e), \beta(e)) \in H_e \oplus H_e \mid e \in E\} \leq T(E \oplus E)$ to $H$ and that $T\Lambda_a$ map $H$ to itself for all $a \in \mathbb R$. This in particular implies that $H_{0_x} = T_xM$ where we are identifying $M$ with the zero section of $TM$. In order to describe the correspondence we first define the \emph{horizontal lift} (relative to an Ehresmann connection $\eta \colon H \to M$ on the fibre bundle $\pi$) as the fibre bundle isomorphism
\begin{equation}\label{eq:horLift}
\mathscr h \colon \pi^* \tau M \to \eta,\quad E_x \times T_xM \ni (e(x),U(x)) \mapsto \mathscr h(e(x))U(x) \coloneqq T_{e(x)}\pi\big|_{H_{e(x)}}^{-1}(U(x))
\end{equation}
i.e.\ $\mathscr h$ is a splitting of the short exact sequence of vector bundles:
\begin{equation} \label{firstexact}
\begin{tikzcd}
0 \arrow[r] & V\pi \arrow[r] & \tau E \arrow[r,"T\pi",swap] & \pi^*\tau M \arrow[r] \arrow[l, bend right,swap,"\mathscr h"] & 0
\end{tikzcd}
\end{equation}
The Ehresmann connection associated to a covariant derivative (where $\pi$ now is a vector bundle) is given in terms of its horizontal lift as
\begin{equation}\label{eq:horlift}
\mathscr h(e(x))U(x) \coloneqq T_x e(U(x)) - \mathscr v(e(x))\nabla_{U(x)} e
\end{equation}
for any section $e \in \Gamma \pi$ whose value at $x$ is $e(x)$ (the independence on the section $e$ is checked by using the usual characterisation of tensoriality \cite[Lemma 2.4]{L97}, i.e.\ by showing that $\mathscr h(f e(x))U(x) = f(x)\mathscr h(e(x))U(x)$: this is easily done in local coordinates).

If we have a chart $\varphi \colon A \to \bbR^m$ for $A \subseteq M$, a chart $\phi \colon B \to \bbR^n$ for $B \subseteq R$ (the typical fibre of $\pi$, an arbitrary $n$-dimensional manifold) and a trivialisation $\Phi \colon E_A \to A \times R$, the triple $(\varphi, \phi, \Phi)$ defines a chart
\begin{equation}\label{eq:productCoords}
(\varphi \times \phi) \circ \Phi \colon \Phi^{-1}(A \times B) \to \bbR^m \times \bbR^n
\end{equation}
We will call the resulting coordinates \emph{product coordinates}. If $\pi$ is a vector bundle, $R$ can (and always will) be taken equal to $\bbR^n$ and $\phi$ to the identity, and if $\pi = \tau M$ or $\tau^*M$, $\Phi$ can be defined canonically in terms of $\varphi$ as $T\varphi$ or $T^*\varphi^{-1}$. In these cases we will speak of \emph{induced coordinates}.
\begin{conv}\label{conv:indices}
	In what follows we will be working on the manifolds $TM$ (or $T^*M$) and $E$. It will therefore be helpful to establish conventions regarding indexing of the product and induced coordinates. In the absence of other manifolds, ambiguities as to the chart, etc.\ we will denote with Greek indices $\alpha,\beta,\gamma,\ldots = 1,\ldots,m$ the coordinates on $M$, with Latin indices $i,j,k,\ldots = m+1,\ldots,m+n$ the coordinates on $E$ in excess of the aforementioned coordinates of the base space $M$ and with $\widetilde \alpha,\widetilde \beta,\widetilde \gamma,\ldots = m+1,\ldots,2m$ the induced coordinates on $TM$ in excess of those on $M$. More specifically, $\widetilde \gamma \coloneqq m + \gamma$, and we will take this into account when using the Einstein convention, e.g.\ $a_{\alpha\beta}b^{\widetilde \beta \gamma} = \sum_{\beta = 1}^m a_{\alpha\beta}b^{(m+\beta)\gamma}$. Moreover, we will use capital letters $I,J,K,\ldots = 1,\ldots,m+n$ to denote indices that run through all coordinates on $E$, and capital letters $A,B,C,\ldots = 1,\ldots,2m$ to denote indices that run through all the coordinates on $TM$. The following diagrams should help explain this arrangement:
	\begin{equation}
	\begin{split}
	E&: \quad (\underbrace{\overbrace{x^1,\ldots,x^m}^{\alpha,\beta,\gamma,\ldots}, \overbrace{y^1,\ldots, y^n}^{i,j,k,\ldots}}_{I,J,K,\ldots}) \\
	TM &: \quad (\underbrace{\overbrace{x^1,\ldots,x^m}^{\alpha,\beta,\gamma,\ldots}, \overbrace{\widetilde x^{\widetilde 1},\ldots, \widetilde x^{\widetilde m}}^{\widetilde \alpha,\widetilde \beta,\widetilde \gamma,\ldots}}_{A,B,C,\ldots})
	\end{split}
	\end{equation}
	A similar convention is followed when $T^*M$ is replaced with $TM$.
	
	It is important to point out the following potential source of confusion. If $V(x) \in T_xM$ it can be either viewed as a vector in the vector space $T_xM$, with coordinates $V^\gamma(x)$, or as a point in the manifold $TM$, with coordinates
	\begin{equation}
	(V(x)^\gamma, V(x)^{\widetilde \gamma}) = (x^\gamma, V^\gamma(x))
	\end{equation}
	Note the different meaning of $V(x)^\gamma$ and $V^\gamma(x)$; in any case, this ambiguity will be avoided by always considering elements as vectors whenever otherwise mentioned. The use of the twidled indices is seen when considering vectors in $TTM$ and $TT^*M$.
	
	Finally, we mention that the use of Greek/Latin indices will also be used in the separate case in which we are dealing with two different manifolds $M$ and $N$, to distinguish between coordinates on the two manifolds.
\end{conv}
In the case of $\pi$ a vector bundle the change of product coordinates from $\varphi,\Phi$ to $\overline \varphi,\overline\Phi$ can be written as 
\begin{equation}\label{eq:changeProdCoords}
\partial^{\overline K}_K(y) = \begin{pmatrix}
\partial^{\overline \gamma}_\gamma(x) & 0 \\ \partial_\gamma \lambda^{\overline k}_k(x) y^k & \lambda^{\overline k}_k(x)
\end{pmatrix}, \quad  (\overline \Phi \circ \Phi^{-1})(x,y)= \mathbbm (x, \lambda(x)y)
\end{equation}
for $x = \pi(y)$ and $\lambda \in C^\infty (\varphi(A), \mathcal L(\bbR^n,\bbR^n))$. It is worthwhile to specify this to the cases of $\pi = \tau M$ (where $\Phi = T\varphi$) and $\tau^*M$ ($\Phi = T^*\varphi^{-1}$), so $\lambda^{\widetilde {\overline \gamma}}_{\widetilde \gamma} = \partial^{\overline\gamma}_\gamma$ and $\lambda^{\widetilde {\overline \gamma}}_{\widetilde \gamma} = \partial^{\gamma}_{\overline\gamma}$ respectively, and 
\begin{align}
\pi = \tau M:& \quad \partial^{\overline C}_C(y) = \begin{pmatrix}[1.5]
\partial^{\overline \gamma}_\gamma & \partial^{\overline \gamma}_{\widetilde\gamma} \\ \partial^{\widetilde{\overline \gamma}}_\gamma & \partial^{\widetilde{\overline \gamma}}_{\widetilde\gamma}
\end{pmatrix}(y) = \begin{pmatrix}[1.5]
\partial^{\overline \gamma}_\gamma(x) & 0 \\ \partial_{\gamma\alpha}^{\overline \gamma}(x) y^\alpha & \partial^{\overline \gamma}_\gamma(x)
\end{pmatrix}\label{eq:changeTanCoords}\\
\pi = \tau^* M:& \quad \partial^{\overline C}_C(y) = \begin{pmatrix}[1.5]
\partial^{\overline \gamma}_\gamma & \partial^{\overline \gamma}_{\widetilde\gamma} \\ \partial^{\widetilde{\overline \gamma}}_\gamma & \partial^{\widetilde{\overline \gamma}}_{\widetilde\gamma}
\end{pmatrix}(y) = \begin{pmatrix}[1.5]
\partial_{\overline \gamma}^\gamma(x) & 0 \\ \partial_{\overline \beta \overline \gamma}^{\alpha} \partial^{\overline \beta}_\gamma(x) y_\alpha & \partial_{\overline \gamma}^\gamma(x)
\end{pmatrix}\label{eq:changeCotanCoords}
\end{align}

We proceed by providing the expression of the horizontal lift in induced coordinates in the case of $\pi = \tau M$ and $\tau^*M$:\todo{Find suitable reference or include the calculation below (commented out). I think Hsu has it (at least the tangent case)} for the tangent space we have
\begin{equation}\label{eq:horLiftCoords}
(\mathscr h(V)U)^\gamma = U^\gamma, \quad (\mathscr h(V)U)^{\widetilde \gamma} = - \Gamma^\gamma_{\alpha\beta} V^\beta U^\alpha
\end{equation}
and in the case of the cotangent bundle
\begin{equation}
(\mathscr h(\omega) U)^\beta = U^\beta, \quad (\mathscr h(\omega) U)^{\widetilde \beta} = U^\alpha \omega_\gamma \Gamma^\gamma_{\alpha\beta}
\end{equation}
Note that in both cases the coordinates of a horizontal lift are not only linear in the vector being lifted, but in the point in $TM$ (or $T^*M$) at which the lift is based: this is an expression of the linearity of the connection, and will be important to guarantee linearity of parallel transport in \autoref{sec:par}.

It will be helpful to define the \emph{frame bundle} $\phi M \colon FM \to M$, the subbundle of $\tau M^{\oplus m}$ whose fibre at $x \in M$ is given by all $m$-frames (i.e.\ ordered bases) of $T_xM$. Since $FM$ is an open subspace of $TM^{\oplus m}$ it makes sense to use the product coordinates of the latter for the former: these are canonically defined in terms of a chart on $M$ by pairs $(\lambda, \gamma) \in \{1,\ldots,m\}^2$ with the first referring to the copy of $TM$, i.e.\ if $y \in F_xM$ then $y_{\lambda} \coloneqq \text{pr}_{\lambda}(y) \in T_xM$ has coordinates $y^\gamma_{\lambda} = y^{(\lambda,\gamma)}$. If $M$ is Riemannian we may additionally consider the \emph{orthonormal frame bundle} $oM \colon OM \to M$, i.e.\ the subbundle of $\phi M$ with total space consisting of orthonormal frames.

We define the \emph{fundamental horizontal vector fields} $\mathscr H_{\lambda} \in \Gamma \tau FM$, $\lambda = 1,\ldots,m$ by the property $T_y\text{pr}_{\gamma}(\mathscr H_{\lambda}(y)) = \mathscr h(y_{\gamma})y_{\lambda}$, or in coordinates
\begin{equation}\label{eq:funHorCoords}
\mathscr H_{\lambda}^\gamma(y) = y^\gamma, \quad \mathscr H_{\nu}^{(\mu,\gamma)}(y) = -\Gamma^\gamma_{\alpha \beta}(x) y^\beta_{\mu} y^\alpha_{\nu}
\end{equation}
with $y \in F_xM$. If $M$ is Riemannian and $\nabla$ is metric these vector fields restrict to elements of $\Gamma \tau OM$.

To end this subsection, we briefly describe what it means for a smooth map of manifolds to preserve connections. Here we are following \cite{E89}.
\begin{defn}[Affine map]
	Let ${^M\!}\nabla$ (${^N\!}\nabla$) be a linear connection on the tangent bundle of the smooth manifold $M$ ($N$). We will say that $f \in C^\infty(M,N)$ is \emph{affine} if 
	\begin{equation}
	\forall U,V \in \Gamma\tau M \quad	T_xf({^M\!}\nabla_{U(x)}V) = {^N\!}\nabla_{T_xf(U(x))}Tf(V)
	\end{equation}
\end{defn}
Note that the RHS is well-defined, as $Tf(V)$ need only be defined on a curve tangent to $U(x)$ at $x$. The name is justified by the fact that the terminology coincides with the usual notion of affinity for smooth maps of Euclidean spaces. Other examples of affine maps are isometries of Riemannian manifolds (Riemannian isomorphisms that is - local isometries are not affine in general). In terms of the Hessians affinity of $f$ reads
\begin{equation}
T_x^*f^{\otimes 2}({^N\!}\nabla^2 g)(f(x)) = {^M\!}\nabla (g \circ f)(x), \quad g \in C^\infty N
\end{equation}
Symmetrising this identity yields the notion of \emph{symmetric affinity}: this is equivalent to the requirement that $f$ preserve parametrised geodesics, with full affinity holding if $f$ additionally preserves torsion. The most useful characterisation of affinity, however, is the local one
\begin{equation}\label{eq:affLocal}
({^{M,N}\!}\nabla^2 f)^k_{\alpha\beta}(x) \coloneqq \partial_{\alpha\beta}f^k(x) + {^N\!}\Gamma^k_{ij}(f(x))\partial_\alpha f^i\partial_\beta f^j(x) - {^M\!}\Gamma^\gamma_{\alpha\beta}\partial_\gamma f^k(x) = 0
\end{equation}
which symmetrised yields the condition for symmetric affinity:
\begin{equation}\label{eq:sAffLocal}
\partial_{\alpha\beta}f^k(x) = \tfrac 12 ({^M\!}\Gamma^\gamma_{\alpha\beta} + {^M\!}\Gamma^\gamma_{\beta\alpha})\partial_\gamma f^k(x) - \tfrac 12 ({^N\!}\Gamma^k_{ij} + {^N\!}\Gamma^k_{ji})(f(x))\partial_\alpha f^i\partial_\beta f^j(x)
\end{equation}
Of course, there is no difference between the two if both connections are torsion-free. Symmetrised expressions will be of interest to us because of the symmetry of the bracket of a rough path; to lighten the notation we will add $(\alpha\beta)$ in an expression to mean that we are symmetrising it w.r.t.\ to the indices $\alpha,\beta$. For instance, symmetric affinity can be written as $({^{M,N}\!}\nabla^2 f)^k_{\alpha\beta}(x) \stackrel{(\alpha\beta)}{=} 0$ or more succinctly still as $({^{M,N}\!}\nabla^2 f)^k_{(\alpha\beta)}(x) = 0$.

\begin{expl}[Affinity and fibre bundles]\label{expl:twidleAffine}
	It will be important to consider whether the projection map of a fibre bundle $\pi \colon E \to M$ is an affine map w.r.t.\ to chosen linear connections $\widetilde \nabla$ on $E$ and $\nabla$ on $M$. By \autoref{eq:affLocal}, the condition of $\pi$ of being affine reads in coordinates
	\begin{equation}\label{eq:affineTauM}
	\partial_{IJ} \pi^\gamma = \partial_K \pi^\gamma \widetilde \Gamma^K_{IJ} - \Gamma^\gamma_{\alpha\beta} \partial_I \tau M^\alpha \partial_J \tau M^\beta
	\end{equation}
	where the $\widetilde \Gamma$'s denote the Christoffel symbols of $\widetilde \nabla$.
	Keeping in mind that $(\varphi \circ \pi \circ T \varphi^{-1})$ is the map $(x^1,\ldots, x^m,y^1,\ldots,y^n) \mapsto (x^1,\ldots, x^m)$ we compute
	\begin{equation}
	\begin{split}
	\partial_{IJ} \pi^\gamma = 0,\quad
	\partial_\beta \pi^\alpha = \delta^\alpha_\beta, \quad \partial_k \pi^\gamma = 0
	\end{split}
	\end{equation}
	and \autoref{eq:affineTauM} becomes
	\begin{equation}\label{eq:tauAff}
	\widetilde \Gamma^\gamma_{\alpha \beta} = \Gamma^\gamma_{\alpha \beta}, \quad \widetilde \Gamma^\gamma_{\alpha \widetilde \beta} = 0, \quad \widetilde \Gamma^\gamma_{\widetilde \alpha \beta} = 0, \quad \widetilde \Gamma^\gamma_{\widetilde \alpha \widetilde \beta} = 0
	\end{equation}
	It is similarly checked that if $\pi$ is a vector bundle the condition of the inclusions of the fibres (as flat spaces) of being affine reads
	\begin{equation}
	\widetilde \Gamma^K_{\widetilde \alpha \widetilde \beta} = 0
	\end{equation}
	and the condition of the inclusion of $M$ (as the zero section) of being affine reads
	\begin{equation}
	\widetilde \Gamma^\gamma_{\alpha\beta}(0_x) =  \Gamma^\gamma_{\alpha\beta}(x)
	\end{equation}
	Replacing symmetric affinity with affinity results in the above coordinate expressions being symmetrised in the bottom two indices of each Christoffel symbol, e.g.\ the symmetrisation of $\widetilde \Gamma^\gamma_{\alpha \widetilde \beta} = 0$ is $\widetilde \Gamma^\gamma_{\alpha \widetilde \beta}  + \widetilde \Gamma^\gamma_{\widetilde\beta \alpha}= 0$ (not $\widetilde \Gamma^\gamma_{\alpha \widetilde \beta}  + \widetilde \Gamma^\gamma_{\widetilde\alpha \beta}= 0$).
\end{expl}

\subsection{Embedded manifolds}\label{subsec:embedded}
Although we have chosen to write this paper mainly in the framework of intrinsic manifolds and local coordinates, we will relate our work to \cite{CDL15}, in which manifolds are embedded in Euclidean space. In this subsection we revisit some of the notions of the previous subsection, assuming that $M$ is Riemannian and isometrically embedded in $\bbR^d$, $\imath \colon M \hookrightarrow \bbR^d$ (this is always possible by the Nash embedding theorem, for high enough $d$). This means that the connection on $M$ will always be the Levi-Civita connection of the induced metric: this setting is less general than the one considered in the previous subsection, where non-metric connections with torsion were considered.  In order to precisely distinguish between extrinsic and intrinsic formulae, we will always distinguish between objects on $M$ (which will be treated using local coordinates, indexed by Greek letters $\alpha,\beta,\gamma,\ldots$) and their counterparts on $\mathcal M \coloneqq \imath(M)$ (treated using ambient coordinates, indexed by the letters $a,b,c,\ldots$). For instance $T_y \mathcal M$ and $T_y^\bot \mathcal M$ (the normal space) are subspaces of $T_y \bbR^d$, $T_y \bbR^d = T_y \mathcal M \oplus T_y^\bot \mathcal M$, and $T_x \imath \colon T_x M \to T_{\imath(x)} \mathcal M$ is an isomorphism.

The geometry of $\mathcal M$ is characterised by the field of orthogonal projection maps
\begin{equation}
P \in \Gamma \mathcal L(\tau_\mathcal{M} \bbR^d,\tau \mathcal M),\quad \text{i.e.\ } P(y) \in \mathcal L(T_y\bbR^d,T_y \mathcal M), \ y \in \mathcal M
\end{equation}
The use of $\mathcal M$ above means we are availing ourselves of the ambient coordinates to think of $P$ as a $d \times d$ matrix defined smoothly in $y \in \mathcal M$. In practice $\mathcal M$ is often defined (locally) through a Cartesian equation $F(y) = 0$, $F \colon \bbR^d \to \bbR^{d-m}$ with surjective differential, in which case the projection map is given by
\begin{equation}\label{eq:PDF}
P(y) = DF^\intercal \circ (DF \circ DF^\intercal)^{-1} \circ DF(y)
\end{equation}
Note that $DF \circ DF^\intercal$ is invertible thanks to the surjectivity of $DF^\intercal$. Also note that this expression provides a smooth extension of $P$ to a tubular neighbourhood of $\mathcal M$ in $\bbR^d$, although this depends on $F$ (whereas $P$ only depends on $\mathcal M$). We also define
\begin{equation}
Q(y) \coloneqq I_d - P(y) \in \mathcal L(T_y\bbR^d, T^\bot_y \mathcal M),\quad y \in \mathcal M
\end{equation}
the orthogonal projection of $T_y\bbR^d$ onto the normal bundle of $\mathcal M$. Of course we have
\begin{equation}\label{eq:PP}
P_c P^c(y) = P(y),\quad Q_c Q^c(y) = Q(y), \quad P_c Q^c(y) = 0 = Q_c P^c(y)
\end{equation}
for $y \in \mathcal M$.

Another important map is the Riemannian projection, uniquely determined and defined smoothly in a tubular neighbourhood $A$ of $\mathcal M$ in $\bbR^d$ \cite[p.132]{Pet06}
\begin{equation}\label{eq:piPidef}
\pi \colon A \to M,\quad y    \mapsto \underset{x \in M}{\arg\min}\abs[0]{y - \imath(x)}, \quad  \Pi \coloneqq \imath \circ \pi \colon A \to \mathcal M
\end{equation}
The important features of $\pi$ and $\Pi$ are
\begin{equation}\label{eq:idPi}
\pi \circ \imath = \mathbbm 1_M,\quad \Pi \circ \Pi = \Pi, \quad \imath = \Pi \circ \imath 
\end{equation}

We may express the Levi-Civita covariant derivative $\nabla$ of $M$ in ambient coordinates as follows:
\begin{equation}\label{eq:nablaExt}
T_x \imath \nabla_{U(x)}V = P(\imath (x)){^d\!}\nabla_{T_x \imath U(x)} (T\imath V), \quad \nabla_{ U(x)}\omega = {^d\!}\nabla_{T_x \imath U(x)}(\omega \circ T \pi)
\end{equation} 
for $U(x) \in T_xM$, $V \in \Gamma \tau M$, $\omega \in \Gamma \tau^*M$, where we have smoothly extended $T \imath (V)$ to a vector field on a tubular neighbourhood of $\mathcal M$ and ${^d\!}\nabla$ is the canonical covariant derivative on $\bbR^d$ given by taking directional derivatives
\begin{equation}
{^d\!}\nabla_{W(y)}Z \coloneqq W^c \partial_c Z(y), \quad {^d\!}\nabla_{U(y)}\rho \coloneqq U^c \partial_c \rho(y)
\end{equation}
We reiterate that computations are carried out in ambient coordinates, i.e.\ $\partial_c$ is differentiation in the $c$-th variable of $\bbR^d$, and sums go from $1$ to $d$. The Levi-Civita Hessian is the given by
\begin{align}
\nabla^2f = T^* \imath^{\otimes 2} {^d\!}\nabla^2 (f \circ \pi), \quad f \in C^\infty M
\end{align} 
where ${^d\!}\nabla$ is the usual Hessian in $\bbR^d$. 

Differentiating $F \circ \Pi = 0$ (where $F = 0$ is a Cartesian equation defining $\mathcal M$), and the fact that $T^\bot_x \mathcal M$ is spanned by the gradients of the components of $F$, shows that $D\Pi|_{T^\bot \mathcal M} = 0$. Since $\Pi|_\mathcal{M} = \mathbbm 1_\mathcal{M}$, $D\Pi|_{T \mathcal M} = \mathbbm 1_{T\mathcal M}$ we have
\begin{equation}\label{eq:Ppi}
P(y) = D \Pi(y), \quad y \in \mathcal M
\end{equation}
Moreover, if $W(y) \in T_y \mathcal M$, differentiating $\partial_b \Pi(Y_t) = P_b(Y_t)$ at time $0$, where $Y$ is a smooth curve in $M$ with $Y_0 = 0$, $\dot Y_0 = W(y)$ we have $\partial_a P_b W^a(x) = \partial_{ab}  \pi W^a(x)$, or in other words
\begin{equation}\label{eq:partialP}
\partial_a P_b P^a(y) = \partial_{ab}  \Pi P^a(y),\quad y \in \mathcal M
\end{equation}
(and in particular the LHS is independent of the extension of $P$ to a tubular neighbourhood). Another useful fact about the second derivatives of $\Pi$ is the following identity, obtained by differentiating \autoref[second identity]{eq:idPi} twice at $y \in M$ and applying \autoref{eq:Ppi}:
\begin{equation}\label{eq:delP}
\partial_{ce}  \Pi P^c_a P_b^e(y) + P_c \partial_{ab} \Pi^c(y) = \partial_{ab} \Pi(y)
\end{equation}
This implies that for $W(y), Z(y) \in T_y\mathcal M$
\begin{equation}\label{eq:d2piNorm}
\partial_{ab} \Pi W^a Z^b(y) = Q_c \partial_{ab}  \Pi^c W^a Z^b(y) \in T^\bot_y \mathcal M
\end{equation}
This coincides with the second fundamental form of $W(y), Z(y)$, i.e.\ $Q(y) {^d\!}\nabla_{W(y)} Z$, since
\begin{equation}
\begin{split}
Q_c(y) ({^d\!}\nabla_{W(y)} Z)^c &= Q_c \partial_a Z^c W^a(y)\\
&= Q_c \partial_a(\partial_b \Pi^c Z^b) W^a(y) \\
&= Q_c(\partial_{ab} \Pi^c Z^b + P^c_b Z^b)W^a(y) \\
&= Q_c \partial_{ab} \Pi^c W^aZ^b(y)
\end{split}
\end{equation}
and is an extrinsic quantity (it cannot be defined on $M$ without the embedding). Finally, we may express the Christoffel symbols of $\nabla$ (w.r.t.\ some chart on $M$) through $\imath, \pi$ as follows:
\begin{equation}\label{eq:Gammai}
\Gamma^\gamma_{\alpha\beta}(x) = \partial_c \pi^\gamma(\imath(x)) \partial_{\alpha\beta} \imath^c(x)
\end{equation}
To prove this identity, let $\widetilde \partial_\alpha, \widetilde \partial_\beta$ be extensions to a tubular neighbourhood of $\mathcal M$ of $T \imath \partial_\alpha = \partial_\alpha \imath, T \imath \partial_\beta = \partial_\beta \imath$ respectively. \autoref{eq:nablaExt} implies
\begin{equation}
\Gamma^\gamma_{\alpha\beta}(x) = \partial_c \pi^\gamma ({^d\!}\nabla_{\widetilde \partial_\alpha} \widetilde \partial_\beta)^c (\imath(x))
\end{equation}
and
\begin{equation}
\begin{split}
({^d\!}\nabla_{\widetilde \partial_\alpha} \widetilde \partial_\beta) (\imath(x)) &= \partial_e (\widetilde \partial_\beta \imath)(\imath(x)) \partial_\alpha \imath^e(x) \\
&= \partial_e (\partial_\beta \imath \circ \pi)(\imath(x))  \partial_\alpha \imath^e(x) \\
&= \partial_{\gamma\beta} \imath (x) \partial_e \pi^\gamma(\imath(x)) \partial_\alpha \imath^e(x) \\
&= \partial_{\gamma\beta} \imath (x) \delta^\gamma_\alpha \\
&= \partial_{\alpha\beta} \imath (x)
\end{split}
\end{equation}
which concludes the argument.
\todo[inline]{I think I presented things in a bit of a loose and disorderly manner in this subsection, however it does the job. The statements are all short and don't really fit well in lemma environments. I'm not sure how else to structure it.}

\section{Rough paths, rough integration and RDEs on manifolds}\label{sec:rough}
In this section $M$ and $N$ will denote smooth $m$- and $n$-dimensional manifolds respectively. Given a control $\omega$ on $[0,T]$ we say that a continuous path $X \colon [0,T] \to M$ lies in $\mathcal C^p_\omega([0,T],M)$ if for all $f \in C^\infty M$, $f(X) \in \mathcal C^p([0,T],\bbR)$; this agrees with the ordinary definition on vector spaces by \autoref{lem:regSmooth}. Equivalently $X \in \mathcal C^p_\omega([0,T],M)$ if for all charts $\varphi$, $\varphi(X) \in \mathcal C^p([a,b],\bbR^m)$ whenever $X|_{[a,b]}$ is contained in the domain of $\varphi$.

\begin{expl}[Path in a fibre bundle]
	Let $\pi \colon E \to M$ be a smooth fibre bundle with typical fibre $R$. A path $H \in \mathcal C^p_\omega ([0,T], E)$ is characterised as follows: for every local trivialisation $\Phi \colon E_A \to A \times R$ and for every $0\leq a \leq b \leq T$ s.t.\ $H[a,b] \subseteq E_A$, we have $\text{pr}_1 \circ \Phi(H|_{[a,b]}) \in \mathcal C^p_\omega ([a,b],A)$ and $\text{pr}_2 \circ \Phi(H|_{[a,b]}) \in \mathcal C_\omega^p ([a,b],R)$. Examples of such paths are given by smooth sections $\sigma \in \Gamma\pi$ evaluated at $X \in \mathcal C_\omega^p ([a,b],M)$. 
\end{expl}

As a warm-up to the rough path case, we can define the Young integral on a manifold. Let $p \in [1,2)$, $X \in \mathcal C^p_\omega([0,T],M)$, $H \in \mathcal C^p_\omega([0,T],\mathcal L (\tau M,\bbR^e))$ in the fibre of $X$. We can then define the \emph{Young integral}
\begin{equation}\label{eq:YoungMfd}
\int_0^T H^k \dif X \coloneqq \lim_{\abs{\pi} \to 0} \sum_{[s,t] \in \pi} H^k_{\gamma;s} X^\gamma_{st}
\end{equation}
where $X^\gamma \coloneqq \varphi^\gamma(X)$ and $H^k_\gamma \coloneqq (H \circ (T_X\varphi)^{-1})^k_\gamma$ for any chart $\varphi \colon A \to \bbR^m$ with $X_{[s,t]} \subseteq A$. The Riemann summands do not depend on the chart up to $o(\omega(s,t))$, since for another chart $\overline \varphi$, whose components we denote with overlined indices
\begin{align}
\begin{split}
H^k_{\overline \gamma;s} X^{\overline \gamma}_{st} &= (H_s \circ  (T_{X_s} \overline \varphi)^{-1})^k_{\overline \gamma} \cdot  \overline \varphi^{\overline \gamma}(X_{st}) \\
&= (H_s \circ  (T_{X_s} \varphi)^{-1})^k_{\gamma} \cdot  \partial_{\overline\gamma} (\varphi \circ \overline \varphi^{-1})^\gamma(\overline \varphi(X_s))\cdot(\overline \varphi \circ \varphi)^{\overline\gamma}(X)_{st} \\
&\approx (H_s \circ  (T_{X_s} \varphi)^{-1})^k_{\gamma}   \cdot \partial_{\overline\gamma} (\varphi \circ \overline \varphi^{-1})^\gamma(\overline \varphi(X_s))\cdot\partial_{\gamma} (\overline \varphi \circ \varphi^{-1})^{\overline\gamma}( \varphi(X_s)) \cdot X^\gamma_{st}\\
&= H_{\gamma;s} X^\gamma_{st}
\end{split}
\end{align}
where the $\approx$ comes from the fact that the terms of order $2$ and higher in the Taylor expansion are $\approx 0$ since $p<2$. Therefore, the limit is well-defined and converges, since it converges on in every chart. Similarly, given a field of linear homomorphisms $V \in \Gamma\mathcal L(\tau M, \tau N)$ (here $\mathcal L(\tau M, \tau N)$ is the bundle $\mathcal L(TM, TN) \to N \times M$ with fibres $\mathcal L(TM, TN)_{y,x} \coloneqq \mathcal L(T_xM,T_yN)$) we can define the Young differential equation by
\begin{equation}
\dif Y = V(Y,X)\dif X \quad \Longleftrightarrow \quad  \dif Y^k = V^k_\gamma(Y,X)\dif X^\gamma 
\end{equation}
where the coordinates are taken to be w.r.t.\ arbitrary charts on $M$ and $N$. We will give definitions of rough paths, their controlled paths, rough integrals and RDEs in the same spirit, relying on the theory of \autoref{sec:backRps}.
\begin{defn}[Rough path on a manifold]\label{def:rpM}
	Given an atlas $(\varphi \colon A_\varphi \to \bbR^m)_\varphi$ of $M$, an $M$-valued $[2,3) \ni p$-\emph{rough path} controlled by $\omega$ on $[0,T]$, $\uX \in \mathscr C^p([0,T], M)$, consists of a collection of rough paths ${^\varphi \!}\uX = ({^\varphi \!}X, {^\varphi \!}\bbX) \in \mathscr C^\alpha([a_\varphi,b_\varphi],\bbR^m)$, where the intervals $[a_\varphi,b_\varphi]$ are chosen so that their union is $[0,T]$ and no two overlap in a single point, and with the property that for all charts $\varphi,\overline \varphi$ in the atlas s.t.\ $[a_\varphi,b_\varphi] \cap [a_{\overline \varphi},b_{\overline \varphi}] \neq \varnothing$
	\begin{equation}\label{eq:compcond}
	(\overline \varphi \circ \varphi^{-1})_* {^\varphi \!}\uX = {^{\overline \varphi} \!}\uX \in \mathscr C^p_\omega([a_\varphi,b_\varphi] \cap [a_{\overline \varphi},b_{\overline \varphi}],\bbR^m)
	\end{equation}
	The \emph{trace} of $\uX$ is the path $t \mapsto X_t \coloneqq \varphi^{-1} ({^\varphi \!}X_t) \in M$ whenever $t \in [a_\varphi,b_\varphi]$ (independently of $\varphi$), $X \in \mathcal C^p([0,T],M)$. $\uX$ is \emph{geometric}, $\uX \in \mathscr G^p_\omega([0,T],M)$, if ${^\varphi \!}\uX$ is geometric for all $\varphi$.
\end{defn}
To define a rough path on $M$ with trace $X$ we only need as many charts as it takes to cover $X_{[0,T]}$: once the compatibility condition \autoref{eq:compcond} is satisfied for such one such cover, for any further chart $\psi$, ${^\psi \!} \uX$ is automatically defined given this condition; moreover this definition only depends on the smooth structure on $M$ and not on the particular atlas covering the trace, thanks to \autoref[1.]{prop:compat}. The definition of $\mathscr G^p_\omega([0,T],M)$ and of the geometrisation map $\mathscr C^p_\omega([0,T],M) \to \mathscr G^p_\omega([0,T],M)$ is well-defined in charts in the obvious way, thanks to the fact that pushforward commutes with geometrisation \autoref[2.]{prop:compat} (this also guarantees that $\mathscr G^p_\omega([0,T],M)$ is well-defined in the first place). The bracket of a manifold-valued rough path is defined in charts, i.e.\
\begin{equation}\label{eq:bracketTransf}
{^\varphi \!}[\bfX]_{st} \coloneqq [{^\varphi \!}{\bfX}]_{st},\quad {^{\overline \varphi}\!}[\bfX]_{st} \approx \partial_\alpha^{\overline\alpha} \partial_\beta^{\overline \beta}({^\varphi\!}X_s){^\varphi \!}[\bfX]_{st}^{\alpha\beta}
\end{equation} 
for charts $\varphi, \overline \varphi$, thanks to \autoref{prop:KW}; it should be noted that $[\bfX]$ is not an $M$-valued path.

A very similar definition given at the end of \cite{BL15}, where the authors allow for more general than smooth (e.g.\ only Lipschitz) transition maps, focusing on geometric rough paths. Here we will not be concerned with finding the minimal working framework for defining rough paths on manifolds, rather we develop this theory in the familiar context of smooth manifolds (in certain cases endowed with extra structure), keeping in mind that many results can be generalised to the $C^2$ or Lipschitz setting.

Rough paths on manifolds can be pushed forward by smooth maps: if $f \in C^\infty(M,N)$ and $\uX \in \mathscr C^p_\omega([0,T],M)$, $f_* \uX \in \mathscr C^p_\omega([0,T],N)$ is defined by, for a chart $\psi$ on $N$
\begin{equation}
{^\psi\!}(f_* \uX) \coloneqq (\psi \circ f \circ \varphi^{-1})_* {^\varphi\!} \uX 
\end{equation}
independently of the chart $\varphi$ on $M$.
Following \cite{E89} we define an $M$-valued \emph{semimartingale} to be a stochastic process defined on some setup $(\Omega,\mathcal F, P)$ satisfying the usual conditions with the property that $f(X)$ is a real-valued semimartingale for all $f \in C^\infty M$, and denote the set of those defined on the interval $[0,T]$ as $\mathcal S(\Omega, [0,T]; M)$. If $M$ is a finite-dimensional $\bbR$-vector space the two notions of $\mathcal S(\Omega, [0,T];V)$ coincide thanks to It\^o's formula.
\begin{expl}[It\^o and Stratonovich rough paths on $M$]\label{expl:itostratRough}
	Let $X \in \mathcal S(\Omega,[0,T];M)$. We can define its Stratonovich and It\^o lifts respectively by lifting ${^\varphi\!}X$ to ${^\varphi\!}\bbX$ and ${^\varphi\!}\widehat\bbX$ defined in \autoref{eq:ItoStrat} on all stochastic intervals $[a,b]$ s.t.\ $X_{[a,b]} \subseteq A_\varphi$ (the domain of $\varphi$) for $t \in [a,b]$. Crucially, these a.s.\ define $M$-valued stochastic rough paths thanks to \autoref{prop:sound}, and just as in the linear case we have ${_\text{g}\!}\widehat \bfX = \bfX$. These definitions, together with \autoref{rem:stochrough} allow us to restrict all the rough path theory that follows to the semimartingale context and recover the theory of stochastic calculus on manifolds (Stratonovich and It\^o integrals, SDEs, etc.) as presented in \cite{E89}.
\end{expl}

We proceed with the definition of controlled paths, specifically in the case of integrands. While for the definition of rough path we used pushforward to force compatibility, for controlled paths we require it through pullbacks.
\begin{defn}[Controlled integrand]\label{def:contrIntM}
	Let $X \in \mathcal C^p_\omega([0,T],M)$. We define an $\bbR^e$-valued \emph{$X$-controlled integrand} $\bfH = ({^\varphi \!}H, {^\varphi \!}H') \in \mathscr D_X(\mathcal L(\tau M, \bbR^e))$ to be a collection ${^\varphi \!}\bfH \in \mathscr D_{{^\varphi\!}X|_{[a_{\scaleto{\varphi}{2pt}},b_{\scaleto{\varphi}{2pt}}]}}(\bbR^{m\times e})$ with $\varphi, a_\varphi, b_\varphi$ as in \autoref{def:rpM} and 
	\begin{equation}
	(\varphi \circ \overline\varphi^{-1})^* {^\varphi \!}\bfH = {^{\overline\varphi} \!}\bfH, \quad \text{i.e. } H_{\overline\gamma} = H_\gamma \partial^\gamma_{\overline \gamma}, \  H'_{\overline \alpha \overline \beta} = H'_{ \alpha  \beta}\partial^\alpha_{\overline \alpha}\partial^\beta_{\overline \beta} + H_\gamma \partial^\gamma_{\overline \alpha \overline \beta}
	\end{equation}
	The \emph{trace} of $\bfH$ is the path $H \coloneqq {^\varphi\!}H \circ T_X\varphi$, which is valued in the fibre of $X$ of the bundle $\mathcal L (\tau M, \bbR^e) = (\tau^*M)^e$.
\end{defn}

As for $\bbR^d$-valued controlled paths, the most immediate example is given by the evaluation of a 1-form $\sigma \in \Gamma \mathcal L(\tau M, \bbR^e)$: in coordinates this amounts to $\boldsymbol \sigma (X) = (\sigma^k_\gamma(X), \partial_\alpha \sigma_\beta^k(X))$.

A smooth map of manifolds $f \in C^\infty(M,N)$ defines the \emph{pullback} of controlled paths: if $X \in \mathcal C^p_\omega([0,T],M)$ and $\bfH \in \mathscr D_{f(X)}(\mathcal L(\tau N, \bbR^e))$
\begin{equation}
{^\varphi\!}(f^*{\boldsymbol H}) \coloneqq (\psi \circ f \circ \varphi^{-1})^* {^\psi\!}\boldsymbol H
\end{equation}
is defined independently of the chart $\varphi$ on $M$ by \autoref[Pullback]{prop:contrProp} and is checked to be an element of $\mathscr D_X(\mathcal L(\tau M,\bbR^e))$.

\begin{rem}[General controlled paths]\label{rem:genContr}
	More in general, let $\pi \colon E \to M$ be a smooth fibre bundle with typical fibre the $n$-dimensional manifold $R$. We may define $\pi$-valued \emph{$X$-controlled path} as a pair $\bfH = (H,H')$ with $H \in \mathcal C^p_\omega([0,T],E)$ in the fibre of $X$, $H' \in \mathcal C^p_\omega([0,T],\mathcal L(TM, TE))$ in the fibre of $(X,H)$, with $H'_t$ a section of $T_{H_t}\pi$ for all $t \in [0,T]$, i.e.\ $T_{H_t}\pi \circ H'_t = \mathbbm 1_{TM}$. Moreover, we require that for all charts $\varphi \colon A \to \bbR^m$, $\phi \colon B \to \bbR^n$ with $A \subseteq M$, $B \subseteq R$ open and local trivialisations $\Phi \colon E_A \to A \times R$, calling, for all $[a,b] \subseteq [0,T]$ s.t.\ $X_{[a,b]} \subseteq A$, $\Phi(H_{[a,b]}) \subseteq A \times B$, and $s \in [a,b]$
	\begin{align}\label{eq:FFF}
	\begin{split}
	&\mathscr F \coloneqq (\varphi, \phi, \Phi), \quad \Xphi_s \coloneqq \varphi(X_s) \in \bbR^m, \quad (\Xphi_s, \HFi_s) \coloneqq (\varphi \times \phi) \circ  \Phi (H_s) \in \bbR^m \times \bbR^n,\\
	&(\mathbbm 1_{\bbR^m},\HFi'_s) \coloneqq T_{\Phi(H_s)}(\varphi \times \phi) \circ T_{H_s}\Phi \circ H'_s \circ (T_{X_s}\varphi)^{-1} \in \mathcal L( \bbR^m, \bbR^m \oplus \bbR^n)
	\end{split}
	\end{align}
	(the assumptions on $\bfH$ imply that the above coordinate expressions have this form) we have
	\begin{equation}\label{controlledeq}
	{^\mathscr{F} \! }\bfH \in \mathscr{D}_{\Xphi}^p([a,b],\bbR^n), \quad \text{i.e.\ } \HFi_{st} - \HFi_s' \Xphi_{st}  \in O(\omega(s,t)^{2/p}) 
	\end{equation}
	This requirement can be checked to not depend on the choice of $\mathscr F$. Moreover, this definition can be seen to restrict to \autoref{def:contrIntM} for the choice $\pi = \mathcal L(\tau^*M, \bbR^e)$, for which the charts $\phi$ and the trivialisation $\Phi$ can be chosen canonically in terms of $\varphi$. In practice, however, there are no applications of this more general notion of controlled path within our scope, and we will rely solely on the definition of controlled integrand.
\end{rem}

We now give the definition of rough integral on manifolds. Note that this definition already exists for the Stratonovich and It\^o integral of semimartingale \cite[p.93, p.109]{E89}, and thanks to \autoref{rem:stochrough} the notion below extends these when applied to \autoref{expl:itostratRough}. It is easily checked that the na\"ive definition of the rough integral in charts $\int \bfH \dif \bfX \coloneqq \int \bfH_\gamma \dif \bfX^\gamma$ fails to be coordinate-invariant due to the bracket correction in the change of variable formula; to come up with an intrinsic notion we will rely on a connection on $\tau M$.
\begin{defn}\label{def:roughIntM}
	Assume $\tau M$ is endowed with a linear connection $\nabla$ and let $\bfX \in \mathscr C^p_\omega([0,T],M), \bfH \in \mathscr D_X(\mathcal L(\tau M, \bbR^e))$. We define the \emph{rough integral}
	\begin{equation}\label{eq:roughIntMfd}
	\int_0^\cdot \bfH \dif_\nabla \bfX \coloneqq \sum_{[s_\varphi,t_\varphi]} \int_{s_\varphi}^{t_\varphi} \bfH_\gamma \dif \bfX^\gamma + \frac 12 \int_{s_\varphi}^{t_\varphi} H_\gamma \Gamma^\gamma_{\alpha\beta}(X) \dif [\bfX]^{\alpha \beta} \in \mathcal C^p_\omega([0,T],\bbR^e)
	\end{equation}
	where we are summing over a finite partition of $[0,\cdot]$ whose intervals $[s_\varphi,t_\varphi]$ are indexed by charts $\varphi$ with the property that each $[s_\varphi,t_\varphi]$ is contained in the domain of $\varphi$, and the coordinates in the integrals are taken w.r.t.\ to these charts. Moreover, this path can be augmented with the unique second order part $\approx H^i_{\alpha;s} H^j_{\beta;s} \bbX_{st}^{\alpha\beta}$ where the coordinates $\alpha,\beta$ are taken w.r.t.\ to any chart that contains $X_{[s,t]}$, thus defining an element of $\mathscr C^p_\omega([0,T],\bbR^e)$.
\end{defn}

We will often write $\dif_M$ for $\dif_{\nabla}$, especially when more than one manifold is involved. Note how we have defined the rough integral directly as a rough path, without passing through the notion of controlled path. To do so would have required to define what it means for an $\bbR^e$-valued path to be controlled by an $M$-valued one: this is possible by applying the generalised definition of \autoref{rem:genContr} to the trivial vector bundle over $M$ with fibre $\bbR^e$. However, it is much simpler to bypass this step, and we shall do so for solutions of RDEs as well.

\begin{thm}\label{thm:roughIntM}
	\autoref{def:roughIntM} is sound: it depends neither on the partition or on the charts chosen for each interval.
	\begin{proof}
		We begin by dealing with the trace. The first assertion immediately follows from the second by comparing two integrals taken w.r.t.\ two different partitions with that taken w.r.t.\ their common refinement (identity holds by additivity of the rough integral on consecutive time intervals). We then consider two charts $\varphi$, $\overline \varphi$, the latter of whose indices we denote using overlines. Then by \autoref{cor:pushpull} we have
		\begin{equation}
		\int \bfH_{\overline\gamma} \dif \bfX^{\overline\gamma} = \int \bfH_\gamma \dif \bfX^\gamma + \frac 12 \int H_{\overline\gamma} \partial_{\alpha\beta}^{\overline \gamma}(X) \dif [\bfX]^{\alpha\beta}
		\end{equation}
		Moreover, using \autoref{prop:KW} and \autoref{eq:chrChange} we have
		\begin{equation}
		\begin{split}
		\int H_{\overline\gamma} \Gamma^{\overline\gamma}_{\overline\alpha \overline\beta}(X) \dif [\bfX]^{\overline\alpha \overline\beta} &= \int (H_\gamma \partial^\gamma_{\overline \gamma}(X))\cdot(\partial^{\overline \gamma}_\gamma \partial^\alpha_{\overline \alpha}\partial^\beta_{\overline \beta} \Gamma^{\gamma}_{\alpha \beta} + \partial^\delta_{\overline \alpha \overline\beta} \partial^{\overline\gamma}_\delta)(X) \cdot (\partial^{\overline \alpha}_\alpha \partial^{\overline \beta}_\beta(X) \dif[\bfX]^{\alpha\beta}) \\
		&= \int (H_\gamma \Gamma^{\gamma}_{\alpha \beta}(X) + H_\gamma \partial^\gamma_{\overline \alpha \overline\beta}  \partial^{\overline \alpha}_\alpha \partial^{\overline \beta}_\beta(X) )\dif[\bfX]^{\alpha\beta}
		\end{split}
		\end{equation}
		Putting these two identities together, we have
		\begin{equation}
		\begin{split}
		\int \bfH_{\gamma} \dif \bfX^{\gamma} + \int H_{\gamma} \Gamma^{\gamma}_{\alpha \beta}(X) \dif [\bfX]^{\alpha \beta} &= \int \bfH_{\gamma} \dif \bfX^{\gamma} + \int H_{\gamma} \Gamma^{\gamma}_{\alpha \beta}(X) \dif [\bfX]^{\alpha \beta} \\
		&\mathrel{\phantom{=}}+ \int H_\gamma (\partial^\gamma_{\overline \gamma} \partial^{\overline \gamma}_{\alpha\beta} +  \partial^\gamma_{\overline \alpha \overline\beta}  \partial^{\overline \alpha}_\alpha \partial^{\overline \beta}_\beta)(X) \dif[\bfX]^{\alpha\beta}
		\end{split}
		\end{equation}
		But 
		\begin{equation}
		\partial^\gamma_{\overline \gamma} \partial^{\overline \gamma}_{\alpha\beta} + \partial^\gamma_{\overline\gamma} \partial^\delta_{\overline \alpha \overline\beta} \partial^{\overline\gamma}_\delta \partial^{\overline \alpha}_\alpha \partial^{\overline \beta}_\beta = \partial_{\alpha\beta}((\overline \varphi^\gamma \circ \varphi^{-1}) \circ ( \varphi \circ \overline\varphi^{-1})) = \partial_{\alpha\beta}\mathbbm 1^\gamma = 0
		\end{equation}
		which yields the desired identity. As for the second order part, we have
		\begin{equation}
		H^i_{\overline\alpha;s} H^j_{\overline\beta;s} \bbX^{\overline \alpha\overline \beta}_{st} \approx (H^i_{\gamma;s} \partial^{\gamma}_{\overline\alpha}(X_s) H^j_{\delta;s} \partial^{\delta}_{\overline\beta}(X_s)) \cdot(\partial^{\overline\alpha}_\alpha \partial^{\overline\beta}_\beta(X_s) \bbX^{ \alpha \beta}_{st}) = H^i_{\alpha;s} H^j_{\beta;s} \bbX^{\alpha \beta}_{st}
		\end{equation}
		This concludes the proof.
	\end{proof}
\end{thm}

We proceed by proving a few properties of the rough integral on manifolds. The first of these (cf.\ \cite[p.109]{E89} in the case of the It\^o integral) tells us that the definition of the rough integral is indeed the one that yields the correct change of variable formula, i.e.\ in which the second derivative is replaced with the Hessian. Note that \autoref{prop:KW} allows us to define the integral of an element of $K \in \mathcal C^p([0,T],\mathcal L(\tau M^{\otimes 2},\bbR^e))$ above $X$, against $[\bfX]$ in coordinates as 
\begin{equation}
\int K \dif [\bfX] \coloneqq \int K_{\alpha\beta} \dif[\bfX]^{\alpha\beta}
\end{equation}
This is the analogue of \cite[Definition 3.9]{E89} in the rough path context.

\begin{prop}[Properties of the rough integral on manifolds]\label{prop:propRIM}\ \ \ \ \
	\begin{description}
		\item[Exact integrands.] For $f \in C^\infty (M, \bbR^e)$\quad $\displaystyle f(X) - f(X_0) = \int_0^\cdot \emph{\bf{d}}\boldsymbol{f}(X) \edif_\nabla \uX + \frac 12 \int_0^\cdot \nabla^2 f(X) \edif[\uX]$;
		\item[Geometric integrators.] $\int \oH \edif_\nabla \uX$ does not depend on the torsion of $\nabla$, and if $\bfX$ is geometric it is altogether independent of $\nabla$;
		\item[Pushforward-pullback behaviour.] For $\uX \in \mathscr C^p(M)$, $f \in C^\infty (M,N)$, $\oH \in \mathscr D^p_{f(X)}(\mathcal L(\tau N, \bbR^e))$, ${^M\!}\nabla$, ${^N\!}\nabla$ connections on $N$ and $M$ respectively
		\begin{equation}
		\int \oH \edif_N f_*\uX - \int f^* \oH \edif_M \uX = \frac 12 \int H_k ({^{M,N}\!}\nabla^2 f)^k_{\alpha\beta}(X) \edif [\uX]^{\alpha\beta}
		\end{equation}
		where ${^{M,N}\!}\nabla^2 f$ is defined in \autoref{eq:affLocal}. In particular the RHS above vanishes whenever $\bfX$ is geometric or $f$ is symmetrically affine.
	\end{description}
	\begin{proof}
		It suffices to show all three statements in a single chart. The first follows immediately from \autoref{eq:hessiancoords} and \autoref{thm:ito}. The second is evident from the fact that the bracket of a geometric rough path vanishes, and that even when it does not it is a symmetric tensor. The third is handled by using \autoref{cor:pushpull}.
	\end{proof}
\end{prop}

\begin{expl}[Tensorial expansion of the rough integral]
	The Taylor-type approximation
	\begin{equation}
	\int_s^t \bfH \dif_\nabla \bfX \approx H_{\gamma;s} X^\gamma_{st} + H'_{\alpha\beta;s} \bbX^{\alpha\beta}_{st} + \tfrac 12 H_{\gamma;s} \Gamma^\gamma_{\alpha\beta}(X_s) [\bfX]^{\alpha\beta}_{st}
	\end{equation}
	is coordinate-invariant up to $o(\omega(s,t))$, but the single terms in it are not. We may rewrite it as
	\begin{equation}\label{eq:tensExp}
	\int_s^t \bfH \dif_\nabla \bfX \approx H_{\gamma;s} (X^\gamma_{st} + \tfrac 12 \Gamma^\gamma_{\alpha\beta}(X_s)X^\alpha_{st}X^\beta_{st}) + ({^\odot \!}\nabla \bfH)_{\alpha\beta;s} \bbX^{\alpha\beta}_{st}
	\end{equation}
	where for a connection $\nabla$
	\begin{equation}
	(\nabla \bfH)_{\alpha\beta}^k \coloneqq H'^k_{\alpha\beta} - H_{\gamma}^k \Gamma^\gamma_{\alpha\beta}(X) 
	\end{equation}
	(and therefore $({^\odot \!}\nabla \bfH)_{\alpha\beta} = H'^k_{\alpha\beta} - H_{\gamma}^k \Gamma^\gamma_{(\alpha\beta)}(X)$, where we are symmetrising the bottom two indices). $\nabla \bfH$ is defined by analogy with \autoref{eq:nablaCovector}, i.e.\ if $\bfH = \boldsymbol \omega(X)$ for a 1-form $\omega \in \Gamma \mathcal L(\tau M, \bbR^e)$, then $\nabla \bfH = \nabla \omega(X)$. Now all four individual terms $H_{\gamma;s}$, $({^\odot\!}\nabla \bfH)_{\alpha\beta;s}$ (and even $(\nabla \bfH)_{\alpha\beta;s}$), $X_{st} + \frac 12 \Gamma^\gamma_{\alpha\beta}(X_s)X^\alpha X^\beta$ and $\bbX^{\alpha\beta}_{st}$ transform as tensors, in the latter two cases up to an $o(\omega(s,t))$. Note that omitting the symmetrisation in ${^\odot\!}\nabla \bfH$ will result in an incorrect expansion, since the accordingly modified expansion \autoref{eq:tensExp} will not be almost additive, due to the extra term involving the evaluation of the torsion against $X_{su} \wedge X_{ut}$:
	\begin{equation}\label{eq:torsDefect}
	\Gamma^\gamma_{\alpha\beta}(X_s)(X^\alpha_{su} X^\beta_{ut} - X^\beta_{su} X^\alpha_{ut})
	\end{equation}
	If $\bfX$ is geometric the symmetrisation can be omitted by writing the expansion as $H_{\gamma;s}(X^\gamma_{st} + \frac 12 \Gamma^\gamma_{\alpha\beta}(X_s) \bbX^{\alpha\beta}_{st}) + (\nabla \bfH)_{\alpha\beta;s}\bbX^{\alpha\beta}_{st}$ (in this case, of course, the connection is purely auxiliary). All of this seems to suggest that manifolds endowed with non-torsionfree connections are not the correct environment for non-geometric rough integration.
	\todo[inline,backgroundcolor=yellow]{I tried to find an alternative definition of integral which uses the torsion (hoping that there would be some interesting interpretation of the evaluation of the torsion against Levy area terms) but couldn't come up with anything. It would be nice to provide an intuitive explanation of why the functional isn't almost additive in the presence of torsion; this would involve geometric interpretations of the term \autoref{eq:torsDefect} which is what makes almost-additivity fail.} 
\end{expl}

\begin{expl}[It\^o-Stratonovich correction on manifolds]\label{expl:itoStratH}
	We compare the integral of $\bfH \in \mathscr D_X(\mathcal L(\tau M, \bbR^e))$ against $\bfX \in \mathscr C^p_\omega([0,T],M)$ and its geometrisation: by \autoref{expl:difference12} we have, at the path level
	\begin{equation}
	\int \bfH \circ \dif \bfX - \int \bfH \dif_\nabla \bfX = \frac 12 \int \nabla \bfH \dif [\bfX]
	\end{equation}
	This identity is the analogue of \cite[Theorem 5.17]{Dr04} in the context of rough paths. Note how our ability of writing the above It\^o-Stratonovich correction formula in terms of an integral against $\dif [\uX]$ is due to the fact that we are integrating controlled paths. In our context of rough integration this is a necessity, but in stochastic calculus on manifolds one can integrate a much larger class of $\mathcal L(TM,\bbR^e)$-valued processes above $X$, and for these the correction formula will involve the quadratic covariations of the components of $H$ and $X$.
\end{expl}

We will now define RDEs driven by manifold-valued rough paths and with solutions valued in a second manifold. The semimartingale-analogue of the definition below can be found in \cite[p.428]{E90}. A heuristic derivation of the coordinate expression can be derived by writing the \say{intrinsic differential} on a manifold with connection $P$ as $\dif_P\boldsymbol Z^k \coloneqq \dif \boldsymbol Z^c + \Gamma^c_{ab}(Z) \dif [\boldsymbol Z]^{ab}$ and writing the identity $\dif_N \bfY^k = F^k_\gamma(Y,X) \dif_M \bfX^\gamma$:
\begin{equation}\label{eq:intermSDE}
\dif \bfY^k + {^N\!}\Gamma^k_{ij}(Y) \dif [\bfY]^{ij} = F^k_\gamma(Y,X)(\dif \bfX^\gamma + {^M\!}\Gamma^\gamma_{\alpha\beta}(X) \dif \boldsymbol [\bfX]^{\alpha\beta})
\end{equation}
Swapping in $\dif [\bfY]^{ij} = F^i_\alpha F^j_\beta(Y,X) \dif[\bfX]^{\alpha\beta}$ (which follows from \autoref{expl:RDEOps}, given that terms of regularity $p/2$ do not contribute to the bracket) then yields \autoref{eq:RDEM} below. 

\begin{expl}[Riemannian rough integral]\label{expl:rri}
	We may define a $\tau M$-valued $X$-controlled path $\boldsymbol P \in \mathscr D_X(\tau M)$ by a collection of ${^\varphi \!}\boldsymbol P \in \mathscr D_{{^\varphi \!}X}(\bbR^m)$ satisfying the compatibility condition 
	\begin{equation}
	P^{\overline{\gamma}} = \partial^{\overline{\gamma}}_\gamma P^\gamma, \ P'^{\overline \beta}_{\overline\alpha} = \partial^{\overline \beta}_\beta P'^\beta_\alpha \partial_{\overline \alpha}^\alpha + \partial^{\overline \beta}_{\gamma\alpha} P^\gamma \partial^\alpha_{\overline\alpha} 
	\end{equation}
	where we are writing the first index as a superscript since we view $P \in T_XM$ (this definition would fall under the more general \autoref{rem:genContr}). Vector fields evaluated at $X$ (along with their coordinate partial derivatives) are obvious examples. Now, if $\scrg$ is a Riemannian metric on $M$, it is natural to define, for $\boldsymbol{P} \in \mathscr D_X(\tau M)$ the controlled integrand
	\begin{equation}
	\boldsymbol{P}^\flat = (P^\flat_\gamma, P'^\flat_{\alpha\beta}) \coloneqq (\scrg_{\gamma\delta}(X) P^\delta, \scrg_{\beta\delta,\alpha}(X) P^\delta + \scrg_{\beta\delta}(X)P^\delta_\alpha) \in \mathscr D_X(\tau^*M)
	\end{equation}
	and if $\nabla$ is a connection on $M$ (which need not be metric) we may integrate $\boldsymbol{P}$ thanks to the Riemannian metric:
	\begin{equation}
	\int \scrg(\boldsymbol{P}, \dif_\nabla \bfX) \coloneqq \int \boldsymbol{P}^\flat \dif_\nabla \bfX
	\end{equation}
	These definitions can be extended in the multivariate case, i.e.\ when we replace $\tau M$ with $\tau M^{\oplus e}$.
\end{expl}	

\begin{defn}[RDEs on manifolds]\label{def:RDEM}
	Let $V \in \Gamma \mathcal L(\tau M, \tau N)$, $\bfX \in \mathscr C^p_\omega([0,T],M)$ and $y_0 \in N$. We define a \emph{solution} to the RDE
	\begin{equation}\label{eq:RDEM}
	\dif_N \bfY = F(Y,X) \dif_M \bfX, \quad Y_0 = y_0
	\end{equation}
	to mean a $N$-valued rough path $\bfY$ with $Y_0 = y_0$ s.t.\ for any two charts on $M$ and $N$, and on any interval restricted to which $X$ and $Y$ are contained in the respective domains, the following RDE (in the sense of \autoref{def:rde})
	\begin{equation}\label{eq:RDEM}
	\dif \bfY^k = F^k_\gamma(Y,X) \dif \bfX^\gamma + \tfrac 12 (F^k_\gamma(Y,X) {^M\!}\Gamma^\gamma_{\alpha\beta}(X)  - {^N\!}\Gamma^k_{ij}(Y) F^i_\alpha F^j_\beta(Y,X) )\dif [\bfX]^{\alpha\beta}
	\end{equation}
	where coordinates are taken (invariantly) w.r.t.\ the two charts. Note that this implies $\mathbb Y^{ij}_{st} \approx F^i_\alpha F^j_\beta(Y_s,X_s)\mathbb X^{\alpha\beta}_{st}$.
\end{defn}
The coordinate-independence check is analogous to that performed in \autoref{thm:roughIntM} and is therefore omitted. Analogously to the vector space-valued case, notions of global and local solutions can be defined and distinguished, and the smoothness of $F$ ensures local existence and uniqueness of the solution. These results can, as usual, be proved via \say{patching} and applying \autoref{thm:localE}. Also note that, just as for the rough integral, only the connection modulo its torsion is relevant, and is not relevant at all when $\bfX$ is geometric, in which case the usual coordinate expression $\dif \bfY^k = F^k_\gamma(Y,X) \dif \bfX^\gamma$ holds: for this reason we shall omit the $M$ and $N$ subscripts to the differentials in this case.
\begin{rem}[Connections with Schwartz-Meyer theory]\label{rmk:SM}
	We mention that there exists a different approach to the topics of this chapter, one that was used in \cite{E89} for stochastic calculus on manifolds, following the ideas of \cite{Sch82} and \cite{Mey81}. It involves defining a connection-independent \say{second order integral} in which the integrator can be viewed as the combination of the semimartingale $X$ and $\frac 12 [X]$, and in which the integrands are valued in appropriately defined second order cotangent bundles $\mathbb T^*M$. The It\^o and Stratonovich integrals can then be viewed as particular cases of this integral, with different choices of the integrand, the former of which depends on the connection. While this approach could be adapted to the rough path setting (i.e.\ by defining an appropriate class of controlled paths with trace valued in the second order cotangent bundle), we chose the more direct approach of defining the rough integral in coordinates and showing invariance under change of charts. In our context coordinates are, in any case, necessary, since they are used in the definition of rough and controlled paths.
	
	In \cite[Proposition 7.34]{E89} the (scalar) It\^o integral is characterised as being the unique additive map
	\begin{equation}
	\int \ \cdot \ \dif X \colon S(\Omega,[0,T];\tau_X^* M) \to S(\Omega,[0,T];\bbR)
	\end{equation}
	with the property that the chain rule \autoref[Exact integrands]{prop:propRIM} and s.t.\ the associativity axiom
	\begin{equation}
	\int \lambda \dif \big( \textstyle \int H \dif X \big) = \displaystyle \int \lambda H \dif X
	\end{equation}
	for all $\lambda \in \mathcal S(\Omega,[0,T];\bbR)$, hold. The rough path analogue of this characterisation would be cumbersome to formulate, since $\boldsymbol \lambda$ would have to be chosen to be controlled by $X$, and not necessarily by $\int H \dif X$. We can, however, still convince ourselves that \autoref{def:roughIntM} is the \say{correct} definition: indeed, as long as we are searching for a local expression of the form
	\begin{equation}
	\int\bfH \dif_M \bfX \coloneqq \int \bfH_\delta \cdot \boldsymbol f^\delta_\gamma(X) \dif \bfX^\gamma + \int H_\gamma g^\gamma_{\alpha\beta}(X) \dif [\bfX]^{\alpha\beta}
	\end{equation}
	for all $\bfH$ (in particular for $\bfH = \textbf{d}\boldsymbol{\varphi}(X)$), which makes \autoref[Exact integrands]{prop:propRIM} hold, \autoref{thm:doobMey} implies (at least in the non-degenerate case of truly rough $X$) that we must pick $f^\delta_\gamma = \delta^\delta_\gamma$ and $g^\gamma_{\alpha\beta} = \Gamma^\gamma_{\alpha\beta}$. To provide a similarly rigorous justification for \autoref{def:RDEM} we can adapt \cite[Definition 6.35]{E89} and characterise the solution $\bfY$ to \autoref{eq:RDEM} as the unique element of $\mathscr C^p_\omega([0,T],N)$ with the property that for all $\bfH \in \mathscr D_Y(\mathcal L(\tau N,\bbR^e))$ (and it is enough to assume this holds for controlled paths given by scalar exact 1-forms) we have
	\begin{equation}
	\int \bfH \dif_N \bfY = \int \bfK \dif_M \bfX
	\end{equation}
	where $\bfK \in \mathscr D_X(\mathcal L(\tau M,\bbR^e)$ is given in coordinates as $\bfK = (\bfH * \bfY) \cdot \boldsymbol F_*(Y,X)$, i.e.\
	\begin{equation}
	K^c_\gamma = K^c_k F^k_\gamma(Y,X),\quad K'^c_{\alpha\beta} = K'^c_{ij} F^i_\alpha F^j_\beta(Y,X) + K^c_k (\partial_\alpha F^k_\beta + F_\alpha^h \partial_h F^k_\beta)(Y,X)
	\end{equation}
	We will omit this check.
	
	The Schwartz-Meyer framework also applies to SDEs (and in our case would apply to non-geometric RDEs), by considering equations defined by a field of \emph{Schwartz morphisms}, i.e.\ fields of morphisms of appropriately-defined short exact sequences
	\begin{equation}
	\begin{tikzcd}
	0 \arrow[r] & T_xM \arrow[d,"F(y{,}x)"] \arrow[r] &\mathbb T_xM \arrow[r] \arrow[d,"\mathbb{F}(y{,}x)"] &T_xM^{\odot 2} \arrow[r] \arrow[d,"F(y{,}x)^{\odot 2}"] &0 \\
	0 \arrow[r] & T_yN \arrow[r] &\mathbb T_yN \arrow[r] &T_yN^{\odot 2} \arrow[r] &0
	\end{tikzcd}
	\end{equation}
	These are not defined with reference to connections on $\tau M$ and $\tau N$; however, once connections are given, a canonical field of Schwartz morphisms is defined, and the resulting equation is identical to \autoref{eq:RDEM}.\todo{I could elaborate on all of this, but the point of this remark was to give credit w/o taking up too much space. The definitions are purposefully vague for that reason.}
\end{rem}

\begin{rem}[Local existence and uniqueness]\label{rmk:localEM}
	The local existence and uniqueness theorem \autoref{thm:localE} extends verbatim to the case of RDEs on manifolds \autoref{def:RDEM} (where compacts are determined by the manifold topology). This is proven for Schwartz-Meyer SDEs in \cite[Theorem 4]{E90} through an affine embedding argument, but in light of \autoref{rmk:SM} the proof carries over to our rough path setting. In particular, if $N$ is compact \autoref{eq:RDEM} always admits a global solution. The explosive case admits the further distinction of whether the solution admits limit in the Alexandrov compactification of $N$ or not.
\end{rem}

When $M$ is a Euclidean space we may consider autonomous RDEs in which the field $F$ only depends on the solution; this is not possible in the general case since the tangent spaces $T_xM$ are not all canonically identified. The next two examples only deals with manifold-valued semimartingales SDEs, but can be viewed in context of RDEs thanks to \autoref{expl:itostratRough}.
\begin{expl}[Local martingales]
	Recall that if $M$ is endowed with a connection $\nabla$, an $M$-valued local martingales $X$ is an $M$-valued semimartingale s.t.\ for all $f \in C^\infty M$, $f(X) - \tfrac 12 \int \nabla^2 f(X) \dif [X]$ is a local martingale in $\bbR^d$, or in local coordinates
	\begin{equation}\label{eq:locMartCoords}
	\dif_\nabla X = \dif X^\gamma + \tfrac 12 \Gamma^\gamma_{\alpha\beta}(X) \dif [X]^{\alpha\beta}
	\end{equation}
	is the differential of a local martingale in $\bbR^d$. As observed in \cite{E90}, it is easy to see that the martingale-preserving property of It\^o SDEs carries over to the manifold setting: if $X$ is an $M$-valued local martingale and $Y$ is the solution to \autoref{eq:RDEM} (where $\bfX$ is given by the It\^o lift of $X$) is an $N$-valued local martingale.
\end{expl}

\begin{expl}[It\^o diffusions]
	Let $M = \bbR^{1+m}$, $X_t = (t,B_t)$ where $B$ is an $m$-dimensional Brownian motion. Then $F$ can be viewed as a collection of $1+m$ vector fields $F_0,F_1,\ldots,F_m \in \Gamma \tau M$, which we take to not depend on $X$. It is well known that the solution to the Stratonovich SDE $\dif Y = F_0(Y) \dif t + F_\gamma(Y) \circ \dif B^\gamma_t$ is a diffusion with generator $\mathscr L \coloneqq F_0 + \frac 12 \sum_{\gamma = 1}^m F_\gamma^2$ (where $F_\gamma^2$ denotes the differential operator $F_\gamma^2 f(x) \coloneqq F_\gamma(y \mapsto F_\gamma f(y))$): this means that for all $f \in C^\infty M$
	\begin{equation}
	f(Y) - \int_0^\cdot \mathscr L f(Y) \dif t
	\end{equation}
	is a local martingale. It\^o diffusions on manifolds can also be considered: the solution to $\dif_\nabla Y = F_0(Y) \dif t + F(Y) \dif B$ (intended in the same intrinsic sense as \autoref{def:RDEM}) is a diffusion with generator $\mathscr L \coloneqq F_0 + \frac 12 \sum_{\gamma = 1}^m \langle \nabla^2 \cdot , F_\gamma \otimes F_\gamma \rangle$. We may verify this claim in local coordinates:
	\begin{equation}
	\begin{split}
	\dif f(Y) &=  \partial_k f(Y) \dif Y^k + \partial^2_{ij}f(Y) \dif [Y]^{ij} \\
	&= \textstyle\partial_k f(Y) (F_\gamma^k(Y)  \dif B^\gamma_t + (F_0^k - \frac 12 \Gamma^k_{ij} \sum_\gamma F_\gamma^i F_\gamma^j)(Y) \dif t ) + \frac 12 \sum_\gamma \partial_{ij}f  F^i_\gamma F^j_\gamma(Y)\dif t \\
	&= \textstyle\partial_k f(Y) F_\gamma^k(Y) \dif B^\gamma_t + (F_0^k + \frac 12 \sum_\gamma (\partial_{ij}f   - \partial_k f \Gamma^k_{ij})F^i_\gamma F^j_\gamma)(Y) \dif t   
	\end{split}
	\end{equation}
	from which the conclusion follows using \autoref{eq:hessiancoords}. This example carries over to the case in which $F$ depends on $t$ (but is still independent of $B$), in which case $\mathscr L$ will also be time-dependent.
\end{expl}

\begin{expl}[RDEs for which It\^o = Stratonovich]
	Assume we have two RDEs started at the same point
	\begin{equation}\label{eq:twoRDEs}
	\dif_N Y = F (Y,X) \dif_M \bfX,\quad \dif Z = G (Z,X) \circ \dif \bfX,\quad Y_0 = Z_0
	\end{equation}
	(recall that $\circ \dif \bfX \coloneqq \dif {_\text{g}\!}\bfX$, and note that we are only interested in comparing the trace of the solutions). We may then ask what the condition on $F,G$ is that guarantees that these two RDEs are equivalent. Assume $Y = Z$: we then write the two equations in coordinates, using \autoref{expl:rdediff} to write them both as RDEs driven by $(\bfX,[\bfX])$, and apply \autoref{thm:doobMey} (since $\bfX$ is taken to be arbitrary, and the condition must hold in particular in the case of $X$ truly rough) to impose equality of the integrands of $\bfX$ and $[\bfX]$. Since we are searching for conditions that may be stated and are valid  independently of the initial condition and of $F,G$ it makes sense to require these two identities to be valid at all pairs $(y,x) \in N \times M$: therefore we must have (substituting the first identity into the second)
	\begin{align}
	F(y,x) &= G(y,x)\label{eq:itoStratRDE1} \\
	F^k_\gamma(y,x) {^M\!}\Gamma^\gamma_{\alpha\beta}(x)  - {^N\!}\Gamma^k_{ij}(y) F^i_\alpha F^j_\beta(y,x) &\stackrel{(\alpha\beta)}{=} \partial_\alpha F^k_\beta(y,x) + F_\alpha^h \partial_h F^k_\beta(y,x)\label{eq:itoStratRDE2}
	\end{align}	
	(recall that $(\alpha\beta)$ above the equals sign means we are symmetrising the identity in $\alpha,\beta$, which we do for the sake of obtaining a condition that is as sharp as possible). This condition is clearly also sufficient to guarantee that the two RDEs yield the same solution (since they do in each coordinate chart). Having required the first-order condition \autoref{eq:itoStratRDE1}, the second-order one is equivalent to the requirement that for all smooth paths $A$ in $M$ and $B$ in $N$ the following implication holds:
	\begin{equation}\label{eq:geodODE}
	\dot B = F(B,A) \dot A \Longrightarrow \nabla_{\dot B} \dot B = F(B,A) \nabla_{\dot A} \dot A
	\end{equation}
	(we would like to write this condition more succinctly as $\nabla_{FU}FU = F\nabla_{U}U$ for $U \in \Gamma \tau M$, but $FU$ cannot be given a meaning as a vector field on $N$, since $F(y,x)U(x)$ depends on $x$; in \autoref{cor:ptStrat} we will see an important special case in which this can nevertheless be done\todo{Can this be written on the product $M \times N$ using ``partial differentiation w.r.t. $\nabla$'', or otherwise, in a nice intrinsic way?}). Indeed, we have
	\begin{equation}
	\begin{split}
	({^N\!}\nabla_{\dot B} \dot B)^k &= \ddot B^k + {^N\!}\Gamma^k_{ij}(B)\dot Y^i \dot Y^j \\
	&= \frac{\dif}{\dif t} (F^k_\gamma(B,A) \dot A^\gamma) + {^N\!}\Gamma^k_{ij}(B) F^i_\alpha F^j_\beta(B,A) \dot A^\alpha \dot A^\beta \\
	&= \partial_\alpha F^k_\beta(B,A)\dot A^\alpha \dot A^\beta + F^h_\alpha \partial_h F^k_\beta(B,A)\dot A^\alpha \dot A^\beta + F^k_\gamma(B,A) \ddot A^\gamma \\
	&\mathrel{\phantom{=}}+ {^N\!}\Gamma^k_{ij}(B) F^i_\alpha F^j_\beta(B,A)\dot A^\alpha \dot A^\beta \\
	&= F^k_\gamma(B,A)({^M\!}\nabla_{\dot A} \dot A)^\gamma + \big[\partial_\alpha F^k_\beta(B,A) + F^h_\alpha \partial_h F^k_\beta(B,A) \\
	&\mathrel{\phantom{=}}- F^k_\gamma(B,A){^M\!}\Gamma^\gamma_{\alpha\beta}(A) + {^N\!}\Gamma^k_{ij}(B) F^i_\alpha F^j_\beta(B,A)\big]\dot A^\alpha \dot A^\beta
	\end{split}
	\end{equation}	
	Now, the quantity contained in the square bracket is precisely the difference of the LHS and RHS of \autoref{eq:itoStratRDE2}, and it is evaluated at $\dot A \otimes \dot A$. Since $A$ is an arbitrary curve, and by polarisation (for a vector space $V$ the tensors $v \otimes v$ generate $V \odot V$ as $v$ ranges in $V$, as is seen by taking $v = u + w$), this is an arbitrary vector in $TM \odot TM$, concluding the argument. Cf. \cite[Corollary 16]{E90}, where \autoref{eq:geodODE} is specified w.l.o.g.\ to the case in which $A$ is a geodesic: in this case it may be stated by saying that the solution map $A \mapsto B$ maps geodesics to geodesics.
\end{expl}

\begin{expl}[Connections that guarantee It\^o = Stratonovich]
	We now reverse the question, and ask whether there are connections on $M$ and/or $N$ which make the two RDEs in \autoref{eq:twoRDEs} equivalent: the same Doob-Meyer argument still guarantees \autoref{eq:itoStratRDE1} is necessary (again, considering that we are seeking a universal condition), and \autoref{eq:itoStratRDE2} must now be solved for ${^M\!}\Gamma^\gamma_{\alpha\beta}(x)$ and/or ${^N\!}\Gamma^k_{ij}(y)$: because of the symmetrisation, the resulting inhomogeneous linear system therefore has $n m(m+1)/2$ equations and $m(m+1)/2 + n(n+1)/2$ potential unknowns. The connections will then be determined, up to torsion, once it is shown that the Christoffel symbols satisfy the correct change of variables \autoref{eq:chrChange}. Since we see no clear way to accomplish this if the system is underdetermined (and cannot expect to have a general solution if the system is overdetermined) we make the assumption that $n = m$ and that we have fixed one of the two connections: this results in the system being square of order $m^2(m+1)/2$. Moreover, we assume that $F(y,x)$ is an isomorphism $T_xM \cong T_yN$; note that requiring this for all $(y,x)$ is a very strong condition, and in particular implies triviality of $\tau M$, $\tau N$ (since $F(y,x_0) \mathcal B_0$ would define a global frame on $\tau N$ for any $x_0 \in M$ and any frame $\mathcal B_0$ of $T_{x_0}M$, and the same goes for $F(y_0,x)^{-1}$ on $\tau M$). Nevertheless, topologically nontrivial examples of manifolds in which the global condition can hold do exist, e.g.\ the torus $\bbR^2/\mathbb Z^2$ with $F$ defined by translation in $\bbR^2$, which passes to the quotient, and the question is even interesting in the case of Euclidean spaces. Now, if the fixed connection is the one on $N$, the linear system is given by $Ax = b$ where $A$ is the matrix with $m(m+1)/2$ diagonal blocks all equal to the matrix $(F^k_\gamma)$ and zeros elsewhere, $x$ represents the column vector of (symmetrised) Christoffel symbols on $M$ enumerated as ${^M\!}\Gamma^1_{11},\ldots,{^M\!}\Gamma^m_{11},\ldots \ldots, {^M\!}\Gamma^1_{mm},\ldots,{^M\!}\Gamma^m_{mm}$ (according to some enumeration of $11,\ldots,mm$), and $b^k_{\alpha\beta}$ are the constant terms obtained from \autoref{eq:itoStratRDE2} and ordered correspondingly. Letting $\Phi(x,y) = F(y,x)^{-1}$ we then have that the solution must be given by
	\begin{equation}\label{eq:LambdaSol}
	{^M\!}\Gamma^\gamma_{\alpha\beta}(x) \stackrel{(\alpha\beta)}{=} \Phi^\gamma_k(x,y)({^N\!}\Gamma^k_{ij}(y) F^i_\alpha F^j_\beta(y,x) + \partial_\alpha F^k_\beta(y,x) + F_\alpha^h \partial_h F^k_\beta(y,x))
	\end{equation}
	This means, in particular, that the RHS must be independent of $y$ for the solution to be well-defined.\todo{Does anyone understand this condition? It seems extremely specific, and I don't know whether it's a good idea to include the example at all if we cannot explain it satisfactorily.}
	We must now change chart on $M$ and check that these symbols transform in the correct way: we compute
	\begin{equation}
	\begin{split}
	\partial_{\overline \alpha}F^k_{\overline \beta} &= \frac{\partial (F^k_{\overline\beta} \circ \overline\varphi^{-1})}{\partial x^{\overline \alpha}} \\
	&= \frac{\partial (F^k_\beta \circ \varphi^{-1} \circ \varphi \circ \overline\varphi^{-1})}{\partial x^{\overline \alpha}} \\
	&= \frac{\partial (F^k_{\overline \beta} \circ \varphi^{-1})}{\partial x^\alpha} \partial^\alpha_{\overline \alpha} \\
	&= \frac{\partial (F^k_\beta \circ \varphi^{-1} \cdot \partial^\beta_{\overline \beta} \circ \varphi^{-1})}{\partial x^\alpha}\partial^\alpha_{\overline \alpha} \\
	&= \Big(\partial_\alpha F^k_\beta \partial^\beta_{\overline \beta} + F^k_\beta \frac{\partial(\partial^\beta_{\overline \beta} \circ \overline \varphi^{-1} \circ \overline \varphi \circ \varphi^{-1})}{\partial x^\alpha}
	\Big) \partial^\alpha_{\overline \alpha} \\
	&= (\partial_\alpha F^k_\beta \partial^\beta_{\overline \beta} + F^k_\beta \partial^\beta_{\overline \alpha \overline \beta} \partial^{\overline \alpha}_\alpha) \partial^\alpha_{\overline \alpha} \\
	&= \partial_\alpha F^k_\beta \partial^\beta_{\overline \beta}\partial^\alpha_{\overline \alpha} + F^k_\beta \partial^\beta_{\overline \alpha \overline \beta} 
	\end{split}
	\end{equation}
	and thus
	\begin{equation}
	\begin{split}
	{^M\!}\Gamma^{\overline{\gamma}}_{\overline{\alpha} \overline{\beta}} &= \Phi^{\overline\gamma}_k({^N\!}\Gamma^k_{ij} F^i_{\overline\alpha} F^j_{\overline\beta} + \partial_{\overline\alpha} F^k_{\overline\beta} + F_{\overline\alpha}^h \partial_h F^k_{\overline\beta}) \\
	&= \partial^{\overline \gamma}_\gamma \Phi^\gamma_k ({^N\!}\Gamma^k_{ij} (\partial^\alpha_{\overline \alpha} F^i_\alpha) (\partial^\beta_{\overline \beta} F^j_\beta) + (\partial_\alpha F^k_\beta\partial^\beta_{\overline \beta}\partial^\alpha_{\overline \alpha} + F^k_\beta \partial^\beta_{\overline \alpha \overline \beta}) +(\partial^\alpha_{\overline \alpha} F_{\alpha}^h) (\partial^\beta_{\overline \beta}\partial_h F^k_{\beta})) \\
	&= \partial^{\overline \gamma}_\gamma \partial^\alpha_{\overline \alpha} \partial^\beta_{\overline \beta} {^M\!}\Gamma^\gamma_{\alpha\beta} + \partial^{\overline \gamma}_\gamma (\Phi_k^\gamma F_\gamma^k)\partial^\beta_{\overline\alpha \overline \beta}\\
	&= \partial^{\overline \gamma}_\gamma \partial^\alpha_{\overline \alpha} \partial^\beta_{\overline \beta} {^M\!}\Gamma^\gamma_{\alpha\beta} + \partial^{\overline \gamma}_\gamma \partial^\gamma_{\overline\alpha \overline \beta}
	\end{split}
	\end{equation}
	The case in which we fix ${^M\!}\Gamma^k_{ij}$ and solve for ${^N\!}\Gamma^k_{ij}$ is handled similarly: keeping in mind that $(F \odot F)^{-1} = \Phi \odot \Phi$ the solution is given by
	\begin{equation}
	{^N\!}\Gamma^k_{ij}(y) \stackrel{(ij)}{=} \Phi^\alpha_i\Phi^\beta_j(y,x)(F^k_\gamma(y,x) {^M\!}\Gamma^\gamma_{\alpha\beta}(x) - \partial_\alpha F^k_\beta(y,x) + F_\alpha^h \partial_h V^k_\beta(y,x))
	\end{equation}
	provided that its value is independent of $y$.
\end{expl}

\section{The extrinsic viewpoint}\label{sec:extrinsic}
In this paper we have mostly chosen to adopt a local perspective on differential geometry. This choice was motivated by the fact that the most natural definition of rough and controlled paths involve charts, and that therefore the theory stemming from these notions would most easily be handled using local coordinates. While we shall continue with this approach in the next section, one of our objectives is to compare our results with those of the other main paper on this topic, \cite{CDL15}, in which manifolds are handled using an extrinsic approach. To do this, we will revisit the main definitions of the previous section, assuming that all manifolds are smoothly embedded in Euclidean space, and using ambient Euclidean coordinates to express our formulae. We will show that our results do indeed extend those of \cite{CDL15}, in which only geometric rough paths and 1-form integrands are considered. One of the most interesting aspect of this section, however, is that for things to generalise in the correct manner to the case of general controlled integrands, additional non-degeneracy hypotheses will have to be placed on the class of integrands; these are always satisfied if $X$ is truly rough.

We will use the notation introduced in \autoref{subsec:embedded} for embedded manifolds endowed with the Levi-Civita connection of the induced Riemannian metric. We begin by stating when an $\bbR^d$-valued rough path may be considered to lie on $\mathcal M$: this will entail not only the obvious requirement on the trace, but also a condition on the second order part.
\begin{defn}[Constrained rough path]\label{def:consRP}
	Let $\bfX \in \mathscr C^p_\omega([0,T],\bbR^d)$. We will say $\bfX$ is \emph{constrained to $\mathcal M$} if $\Pi_* \bfX = \bfX$, and denote the set of $\mathcal M$-constrained $p$-rough paths controlled by $\omega$ with $\mathscr C^p_\omega([0,T],\mathcal M)$ and its subset of geometric ones with $\mathscr G^p_\omega([0,T],\mathcal M)$.
\end{defn}
$\imath_*$ defines bijections $\mathscr C^p_\omega([0,T],M) \to \mathscr C^p_\omega([0,T],\mathcal M)$ with inverse $\pi_*$, but we still choose to distinguish the two notions, since local coordinates are used in the former case, while $\mathscr C^p_\omega([0,T],\mathcal M) \subseteq \mathscr C^p_\omega([0,T],\bbR^d)$. An equivalent way of stating \autoref{def:consRP} for an $\bbR^d$-valued rough path $\bfX$ is as follows: $X \in \mathcal C^p_\omega([0,T], \mathcal M)$, or equivalently by Taylor's formula
\begin{equation}\label{eq:firstConstr}
X_{st}^c = \Pi^c(X)_{st} \approx P_d^c(X_s)X_{st}^d + \tfrac 12 \partial_{ab} \Pi^c(X_s)X^a_{st}X^b_{st}
\end{equation}
and 
\begin{equation}\label{eq:secondConstr}
\bbX_{st}^{cd} \approx P_a^c P_b^d(X_s)\bbX_{st}^{ab} \Leftrightarrow Q_a^c(X_s)\bbX_{st}^{ab} \approx 0 \Leftrightarrow Q_b^d(X_s)\bbX_{st}^{ab} \approx 0
\end{equation}
Moreover, these imply
\begin{equation}\label{eq:bracketConstr}
[\bfX]^{cd}_{st} \approx P^c_a P^d_b(X_s)[\bfX]^{ab}_{st}
\end{equation}
We note straight away that this definition extends the characterisation \cite[Corollary 3.32 (2)]{CDL15} to the non-geometric setting; the characterisation \cite[Corlollary 3.32 (1)]{CDL15} (which states that $Q_a I_b(X_s)(\bbX_{st}^{ab} - \bbX_{st}^{ba}) \approx 0$) does not hold, however, for non-geometric rough paths, as the symmetric part of their second order part is not determined by their trace (a counterexample is easily found by taking $\bfX \in \mathscr C^p_\omega([0,T],\mathcal M)$ and then adding to $\bbX$ any path $Z \in \mathcal C^{p/2}_\omega([0,T],(\bbR^d)^{\odot 2})$ s.t.\ $Q_b^d(X_s)Z^{ab}_{st} \not \approx 0$).

Instead of defining a notion of \say{constrained controlled path} we directly define a notion of rough integral \say{on $\mathcal M$} which is valid for any path in $\bbR^{e \times d}$ that is controlled by the trace of the integrator. We will then show that, under an additional hypothesis on the integrand, this integral only depends on the restriction of the integrand (and indeed just of its trace) to $T_XM$.
\begin{defn}[Constrained rough integral]\label{def:constRI}
	Let $\bfX \in \mathscr C^p_\omega([0,T],\mathcal M)$, $\bfH \in \mathscr D_X(\bbR^{e \times d})$. We define the \emph{$\mathcal M$-constrained rough integral} of $\bfH$ against $\bfX$ (both as an element of $\mathscr D_X(\bbR^e)$ and as one of $\mathscr C^p_\omega([0,T],\bbR^e)$) as
	\begin{equation}\label{eq:PPiIntId}
	\int \bfH \dif_{\mathcal M} \bfX \coloneqq \int \Pi^* \bfH \dif \bfX = \int (\bfH \cdot \boldsymbol P(X)) \dif \bfX
	\end{equation}
\end{defn}
The identity above is shown by the following simple calculation (we will reuse the letters $e,d$ as indices without the risk of ambiguity)
\begin{equation}\label{eq:equivDefInt}
\begin{split}
\int_s^t \Pi^* \bfH \dif \bfX &\approx H_{d;s} P^d_c(X_s) X^c_{st} + (H'_{ef;s}P^e_a P^f_b(X_s) + H_{d;s} \partial_{ab}\Pi^d(X_s)) \bbX_{st}^{ab} \\
&\approx H_{d;s} P^d_c(X_s) X^c_{st} + (H'_{eh;s} P^h_f(X_s) + H_{d;s} \partial_{eh}\Pi^d P^h_f(X_s)) P_a^e P_b^f(X_s) \bbX_{st}^{ab} \\
&\approx H_{d;s} P^d_c(X_s) X^c_{st} + (H'_{eh;s} P^h_f(X_s) + H_{d;s} \partial_{e} P^d_h P^h_f(X_s) ) P_a^e P_b^f(X_s) \bbX_{st}^{ab} \\
&\approx H_{d;s} P^d_c(X_s) X^c_{st} + (H'_{a h;s} P^h_b(X_s) + H_{d;s} \partial_{a} P^d_b(X_s)) \bbX_{st}^{ab} \\
&\approx \int_s^t (\bfH \cdot \boldsymbol P(X)) \dif \bfX
\end{split}
\end{equation}
where we have used that $\bfX$ is constrained to $M$, \autoref{eq:PP} and \autoref{eq:delP}; at the level of Gubinelli derivatives/second order parts the identity is obvious.

Using \autoref{cor:pushpull} we compute the correction formula for the traces of the ordinary and constrained rough integrals
\begin{equation}\label{eq:MvsRdInt}
\int \bfH \dif \bfX - \int \bfH \dif_{\mathcal M} \bfX = \int \bfH \dif \Pi_*\bfX - \int \Pi^*\bfH \dif \bfX = \frac 12 \int H_c \partial_{ab}\Pi^c(X) \dif [\bfX] 
\end{equation}
while their second-order parts both agree with
\begin{equation}
\bbY_{st}^{ij} \approx H^i_{c;s} H^j_{d;s} P^c_a P^d_b(X_s) \bbX^{ab}_{st} \approx H^i_{a;s} H^j_{b;s} \bbX^{ab}_{st} 
\end{equation}
In particular, if $\bfX$ is geometric
\begin{equation}\label{eq:geomNoP}
\int \bfH \dif_\mathcal{M} \bfX = \int \bfH \dif \bfX
\end{equation}
and hence agrees with \cite[Definition 3.24]{CDL15} when restricted to the case of 1-forms (see \autoref{ex:1formInt} below).

Also note that if $\bfX$ is the It\^o or Stratonovich stochastic rough path associated to a semimartingale, the above definition coincides, thanks to \autoref{rem:stochrough}, with the usual It\^o and Stratonovich integrals, given in extrinsic form in \cite[Definition 5.13]{Dr04}.

Now, it is clear that if $\Pi^* \bfH$ (or equivalently $\pi^* \bfH$, since $\imath_*$ is injective) vanishes, $\int \bfH \dif_\mathcal{M} \bfX$ also vanishes, and we may conclude that the integral depends only on the restriction of $\bfH$ to $\mathcal M$ in the sense that $\Pi^* \bfH = \Pi^*\bfK$ $\Rightarrow$ $\int \bfK \dif_\mathcal{M} \bfX - \int \bfH \dif_\mathcal{M} \bfX$. This, however, falls short of our goal of generalising \cite[Corollary 3.35]{CDL15} (or rather one implication - we will address the second one in \autoref{rem:noConv} below), which states, in our notation, that if $\bfX \in \mathscr G^p_\omega([0,T],\mathcal M)$ then $\int \boldsymbol f(X) \dif \bfX = 0$ for all $f \in \Gamma \mathcal L(\bbR^d,\bbR^e)$ s.t.\ $\imath^*f = 0$. The point is that the requirement is only placed on the trace $f(X)$ of the integrand, not on the whole controlled path. Unfortunately, without further assumptions, the obvious generalisation to the setting of general controlled integrands of this statement fails. The example below exhibits two ways in which this can occur.
\begin{expl}\label{ex:failureWd}
	Take $\mathcal M$ to be the unit circle $S^1$ in $\bbR^2$, so $\Pi$ is given by
	\begin{equation}
	\Pi \colon \bbR^2 \setminus \{(0,0)\} \to \bbR^2, \quad (x,y) \mapsto \frac{(x,y)}{\sqrt{x^2+y^2}} 
	\end{equation}
	Let $\boldsymbol Z \in \mathscr C^p_\omega([0,T],\bbR^2)$ given by
	\begin{equation}
	Z_t \coloneqq (1,0), \quad \mathbb Z_{st} \coloneqq \begin{pmatrix}
	t-s &0 \\ 0 &t-s\end{pmatrix}
	\end{equation}
	which satisfies the Chen identity thanks to the constancy of the trace. Define $\bfX \coloneqq \Pi_* \boldsymbol Z \in \mathscr C^p_\omega([0,T],\mathcal M)$: it is checked that
	\begin{equation}
	X_t = (1,0), \quad \bbX_{st} \approx \begin{pmatrix}
	0 &0 \\ 0 &t-s\end{pmatrix}
	\end{equation}
	Now let 
	\begin{equation}
	H_t = (H_{1;t},H_{2;t}) \coloneqq (1,0), \quad H' \coloneqq 0_{2 \times 2}
	\end{equation}
	Trivially, $(H,H') \eqqcolon \bfH \in \mathscr D_X(\bbR^{1 \times d})$, and we compute
	\begin{equation}
	\begin{split}
	\int_s^t \bfH \dif \bfX &\approx H_{d;s} P^d_c(X_s) X^c_{st} + (H'_{cd;s}P^c_a P^d_b(X_s) + H_{c;s} \partial_{ab}\Pi^c(X_s)) \bbX_{st}^{ab}\\ 
	&= H_{1;s} \partial_{22}\Pi^1(X_s) \bbX_{st}^{22} \\
	&= s-t
	\end{split}
	\end{equation}
	despite the fact that $H|_{T_XM} = 0$ (and even $H'|_{T_XM^{\otimes 2}} = 0$).
	
	Another example is given as follows: let $\mathcal M = \bbR^d$ with $d = 2$ (or embed in $\bbR^3$ if we want non-zero codimension) and let $\bfX$ be the geometric rough path
	\begin{equation}
	X_t \coloneqq (0,0), \quad \bbX_{st} \coloneqq \begin{pmatrix}
	0 &t-s \\ t-s &0\end{pmatrix}
	\end{equation}
	and $\bfH$ be given by
	\begin{equation}
	H \coloneqq 0, \quad H'_{st} \coloneqq \begin{pmatrix}
	0 &1 \\ 0 &0\end{pmatrix}
	\end{equation} 
	Again membership to $\mathscr D_X(\bbR^{2 \times 1})$ is trivially satisfied and proceeding as above we compute
	\begin{equation}
	\int_s^t \bfH \dif \bfX = H'_{cd;s}P^c_a P^d_b(X_s) \bbX_{st}^{ab} = H'_{12} \bbX^{12}_{st} = t-s
	\end{equation}
	
	To summarise, in the first example we were able to have $H|_{T_XM} = 0$, $H'|_{T_XM^{\otimes 2}} = 0$, but the manifold had to be non-flat ($D^2\Pi \neq 0$) and the rough path had to be chosen to be non-geometric (if not $\bbX^{22} = (X_{st}^2)^2 \approx 0$, assuming $X$ is chosen to be in $\mathcal C^{p/2}_\omega([0,T],\bbR^d)$, which is necessary to produce a counterexample). In the second example we were able to choose a geometric rough path, and even a flat manifold, but it was not the case that $H'|_{T_XM^{\otimes 2}} = 0$ (although it still held that $H|_{T_XM} = 0$).
\end{expl}

The fact that in the two examples above the trace of the rough path was chosen to be constant (and in particular an element of $\mathcal C^{p/2}_\omega([0,T],\bbR^d)$) is not accidental. The next proposition places the needed assumptions that rule out this degenerate behaviour of the rough integral.
\begin{thm}\label{prop:wdInt} Let $\bfX \in \mathscr C^p_\omega ([0,T],\bbR^d)$.
	\begin{enumerate}
		\item Let $\bfH \in \mathscr D_X(\bbR^{e \times d})$. If $\bfH = \bfH \cdot \boldsymbol Q(X)$ then $\int \bfH \edif_\mathcal{M} \bfX = 0$ as a rough path;
		\item The above condition on $H'$ is always satisfied if $H_c P^c(X) = 0$ and $X$ is truly rough;
		\item If $\bfH, \bfK \in \mathscr D_X(\bbR^{e \times d})$ with $(\bfK - \bfH) \cdot \boldsymbol Q = 0$, which is always the case if $X$ is truly rough and $H|_{T_XM} = K|_{T_XM}$, then as rough paths
		\begin{equation}
		\int \bfH \edif_\mathcal{M} \bfX = \int \bfK \edif_\mathcal{M} \bfX 
		\end{equation} 
	\end{enumerate}
	\begin{proof}
		The hypothesis $\bfH = \bfH \cdot \boldsymbol Q(X)$ reads
		\begin{equation}\label{eq:wdCond}
		H = H_c Q^c(X), \quad H'_{ab} = H'_{ac} Q^c_b(X) + H_c \partial_a Q^c_b(X) = H'_{ac} Q^c_b(X) - H_c \partial_{ab} \Pi^c(X)
		\end{equation}
		Using this and proceeding as in \autoref{eq:equivDefInt} for the trace we compute
		\begin{equation}
		\begin{split}
		\int_s^t \bfH \dif_\mathcal{M} \bfX &\approx H_{d;s} P^d_c(X_s) X^c_{st} + H'_{ef;s} P^e_a P^f_b(X_s) \bbX^{ab}_{st} + H_{c;s} \partial_{ef} \Pi^c P^e_a P^f_b(X_s) \bbX^{ab}_{st} \\
		&= (H'_{e c;s} Q^c_f - H_{c;s} \partial_{ef} \Pi^c)P^e_a P^f_b(X_s) \bbX^{ab}_{st} + H_{c;s} \partial_{ef} \Pi^c P^e_a P^f_b(X_s) \bbX^{ab}_{st} \\
		&= H'_{e c;s}  P^e_a Q^c_f P^f_b(X_s) \bbX^{ab}_{st} \\
		&= 0
		\end{split}
		\end{equation}
		and the second order part is $\approx H^e_{c;s} H^f_{d;s} Q^c_a Q^d_b(X_s)\bbX^{ab}_{st} \approx 0$. This proves 1.
		
		To prove 2.\ simply observe that $\bfH \cdot \boldsymbol Q(X) \in \mathscr D_X(\bbR^{d \times e})$, so \autoref{thm:gubUniq} implies it can be the only one with trace $K \coloneqq H_c Q^ca(X)$, which is equivalent to $K_c P^c(X) = 0$, implying the statement.
		
		3.\ follows from
		\begin{equation}\label{eq:difference0}
		\int \bfK \dif_\mathcal{M} \bfX - \int \bfH \dif_\mathcal{M} \bfX = \int (\bfK - \bfH) \dif_\mathcal{M} \bfX = 0
		\end{equation}
		If $H|_{T_XM} = K|_{T_XM}$ and $X$ is truly rough we may apply 2.\ to $\bfK - \bfH$.
	\end{proof}
\end{thm}

Of course, when dealing with explicit examples of controlled integrands the true roughness hypothesis is often not needed. To give an example of this, and to relate our results with those of \cite{CDL15} we briefly discuss the example of 1-forms.

\begin{expl}[1-form integrands]\label{ex:1formInt}
	Let $f \in \Gamma \mathcal L(\tau \bbR^d, \bbR^e)$ be a 1-form defined on $\bbR^d$. Assume for the moment that $f(X)|_{T_{X}M} = 0$, then differentiating $(f_c \circ \Pi) P^c (x) = 0$ at $y \in \mathcal M$ we obtain
	\begin{equation}
	\partial_c f_d P^c_a P^d_b (y) + f_c \partial_{ad} \Pi^c P^d_b(y) = 0
	\end{equation}
	which shows that $\int \boldsymbol f (X) \dif_\mathcal{M} \bfX = 0$ by proceeding as in \autoref{eq:equivDefInt}. Note this is a slightly different identity compared to \autoref{eq:wdCond}, which may however be obtained by extending $Q$ outside of $M$, e.g.\ by \autoref{eq:PDF} and assuming (w.l.o.g., as argued below) that $f_c = f_d Q^d_c$ on a tubular neighbourhood of $M$: differentiating this relation then shows $\boldsymbol f(X) = \boldsymbol f(X) \cdot \boldsymbol Q(X)$. Therefore, by arguing as in \autoref{eq:difference0}, for $f \in \Gamma \mathcal L(\tau \bbR^d, \bbR^e)$ the value of $\int \boldsymbol f (X) \dif_\mathcal{M} \bfX$ only depends on $f(X)|_{T_XM}$ (cf.\ \cite[Lemma 3.23]{CDL15} in conjunction with \autoref{eq:geomNoP}). The same conclusion follows if we realise that the formula \cite[Equation 3.17]{CDL15} 
	\begin{equation}
	\int_s^t \boldsymbol f (X) \dif \bfX \approx f_d P^d_c(X_s) X^c_{st} + (\nabla f)_{ef}(X_s) P^e_a P^f_b(X_s) \bbX^{ab}_{st}
	\end{equation}
	extends to the case of non-geometric rough paths and that $\nabla f (x)$ (defined in \autoref{eq:nablaExt}) only depends on $f(x)|_{T_xM}$ for $x \in M$. This is the extrinsic version of \autoref{eq:tensExp} applied to 1-form integrands, and the same expansion would hold for arbitrary controlled paths, by defining
	\begin{equation}
	(\nabla \bfH)_{ab} \coloneqq H'_{ab} + H_c \partial_{ab}\Pi^c(X)
	\end{equation}
	Therefore we have shown that the true roughness assumption is not necessary to integrate ambient 1-forms against constrained rough paths in a manner which is only dependent upon their restriction to $\tau M$, since the first hypothesis in \autoref{prop:wdInt} is automatically satisfied. We have also shown that \autoref{def:constRI} extends \cite[Definition 3.24]{CDL15}.
\end{expl}

\begin{expl}[It\^o-Stratonovich corrections on embedded manifolds]
	By \autoref[2.]{prop:compat} the geometrisation of an $\mathcal M$-constrained rough path is still constrained, and we may use \autoref{expl:difference12} to compute the trace-level difference of the constrained integrals against $\bfX$ and its geometrisation as
	\begin{equation}\label{eq:itoStratExt}
	\int \bfH \circ \dif \bfX - \int \bfH \dif_\mathcal{M} \bfX = \frac 12 \int (\nabla \bfH)\dif [\bfX]
	\end{equation}
	This is the extrinsic version of \autoref{expl:itoStratH}.\todo{In fact it can be deduced from \autoref{thm:exInt} below; decide whether to move it down or to keep it in at all.}
\end{expl}

\begin{rem}\label{rem:noConv}
	In \cite[Corollary 3.20]{CDL15} it is shown that, for $\bfX \in \mathscr G^p_\omega([0,T],\bbR^d)$ with $X$ valued in $\mathcal M$, the condition
	\begin{equation}\label{eq:zeroInt}
	\int \boldsymbol f(X) \dif \bfX = 0 \quad \forall f \in \Gamma \mathcal L(\tau \bbR^d,\bbR^e) \text{ s.t.\ } \Pi^*f = 0
	\end{equation}
	implies \autoref{eq:secondConstr} and thus $\bfX \in \mathscr G^p_\omega([0,T],\mathcal M)$. In order to attempt to generalise this statement to the non-geometric case we must pick which of the integrals in \autoref{eq:PPiIntId} to use; in both cases, however, the statement becomes trivial since, and even replacing the quantifier over 1-forms with one over all controlled integrands $\bfH$, we are dealing with the integral against $\bfX$ of a controlled path with trace $H_c P^c(X) = 0$: if $X$ is truly rough this implies that the whole integrand, and thus the integral, vanishes, regardless of the behaviour of $\bbX$. Note that using the ordinary $\bbR^d$-integral in place of the constrained integral (the two coincide for geometric rough paths by \autoref{eq:geomNoP}) is not meaningful either: indeed, if \autoref{eq:zeroInt} implied $\bfX \in \mathscr C^p_\omega([0,T],\mathcal M)$ for non-geometric $\bfX$, by \autoref{expl:difference12}
	\begin{equation}
	\int \boldsymbol f(X) \circ \dif \bfX = \int \boldsymbol f(X) \dif \bfX + \frac 12 \int Df(X) \dif [\bfX]  =  \frac 12 \int \partial_a f_b(X) \dif [\bfX]^{ab}
	\end{equation}
	which would have to be zero by \autoref{ex:1formInt} and the fact that $\mathscr C^p_\omega([0,T],\mathcal M)$ is closed under geometrisation (or by \autoref{thm:doobMey}). But this is not the case if we pick $\mathcal M$, $\bfX$ as in \autoref[first example]{ex:failureWd}, $e = 1$ and $f(x^1,x^2) = (x^1,x^2)$ (which restricts to $0$ on $T\mathcal M$) we have $[\bfX]^{ab}_{st} = 2\delta^{a2}\delta^{b2}(s-t)$ and therefore $\frac 12 \int \partial_a f_b(X) \dif [\bfX]^{ab} = s-t  \neq 0$, a contradiction.
	
	The only way (that we can think of) to characterise non-geometric rough integrals in terms of ambient ones would be to endow $\bbR^d$ with a connection s.t.\ $\imath$ is symmetrically affine (which can always be done \cite[Lemma 15]{E90}) and replacing the integral in \autoref{eq:zeroInt} with the rough integral in $\bbR^d$ taken w.r.t.\ this connection, in the intrinsic sense of \autoref{def:roughIntM}. This, however, falls short of the goal of characterising constrained rough paths in terms of notions that do not involve manifolds.
\end{rem}

\begin{expl}[Affine subspaces]
	If $M$ is an affine subspace of $\bbR^d$ then $P$ is constant, and 
	\begin{equation}
	\int_s^t \bfH \dif_\mathcal{M} \bfX \approx H_{d;s} P^d_c X^c_{st} + H'_{ef;s} P^e_a P^f_b \bbX^{ab}_{st}
	\end{equation}
	and in particular only depends on $H|_{T_XM}$, $H'|_{T_XM^{\otimes 2}}$. The true roughness hypothesis is still necessary if we want dependence only on $H|_{T_XM}$, as demonstrated by \autoref{ex:failureWd}.
\end{expl}

We still have not related the constrained rough integral with its intrinsic counterpart, defined in \autoref{def:roughIntM}. This is done as follows:
\begin{thm}\label{thm:exInt}
	Let $\bfX \in \mathscr C^p_\omega([0,T],\mathcal M)$, $\bfH \in \mathscr D_X(\bbR^{e \times d})$. Then 
	\begin{equation}
	\int \bfH \edif_{\mathcal M} \bfX = \int \imath^*\bfH \edif_M \pi_* \bfX
	\end{equation}
	\begin{proof}
		Applying \autoref{prop:propRIM} to $\imath$ we obtain
		\begin{equation}
		\begin{split}
		\int \bfH \dif_{\mathcal M} \bfX &= \int \Pi^*\bfH \dif (\imath \circ \pi)_*\bfX \\
		&= \int \imath^* \Pi^* \bfH \dif \pi_*\bfX + \frac 12 \int H_d P^d_c(X) ({^{M,\bbR^d}\!}\nabla^2 \imath)^c_{\alpha\beta}(\pi(X)) \dif [\pi_*\bfX]^{\alpha\beta}
		\end{split}
		\end{equation}
		Now,
		\begin{equation}
		\int \imath^* \Pi^* \bfH \dif \pi_*\bfX = \int (\Pi \circ \imath)^* \bfH \dif \pi_* \bfX = \int \imath^*\bfH \dif \pi_*\bfX
		\end{equation}
		and applying \autoref{eq:affLocal} and \autoref{eq:Gammai}, for $x \in M$, $y \coloneqq \imath(x)$
		\begin{equation}
		\begin{split}
		({^{M,\bbR^d}\!}\nabla^2 \imath)^c_{\alpha\beta}(x) &= \partial_{\alpha\beta} \imath^c(x) - \Gamma^\gamma_{\alpha\beta} \partial_\gamma(y) \imath^c(x) \\
		&= \partial_{\alpha\beta} \imath^c(x) - \partial_e \pi^\gamma(y) \partial_{\alpha\beta} \imath^e \partial_\gamma \imath^c(x) \\
		&= \partial_{\alpha\beta} \imath^c(x) - \partial_e (\imath \circ \pi)^c(y) \partial_{\alpha\beta}\imath^e(x) \\
		&= \partial_{\alpha\beta} \imath^c(x) - P^c_e(x) \partial_{\alpha\beta}\imath^e(y) \\
		&= Q^c_e(y) \partial_{\alpha\beta} \imath^e(x)
		\end{split}
		\end{equation}
		which implies
		\begin{equation}
		H_d P^d_c(X) ({^{M,\bbR^d}\!}\nabla^2 \imath)^c_{\alpha\beta}(\pi(X)) = H_d P^d_c Q^c_e(X) \partial_{\alpha\beta} \imath^e(\pi(X)) = 0
		\end{equation}
		concluding the proof.
	\end{proof}
\end{thm}

We now turn to the extrinsic treatment of RDEs. Let ${^\mathcal{N}\!}\imath \coloneqq N \hookrightarrow \bbR^e$ be a Nash embedding of another Riemannian manifold $N$, ${^\mathcal{N}\!}\imath(N) \eqqcolon \mathcal N$, and ${^\mathcal{N}\!}\pi, {^\mathcal{N}\!}\Pi$ its projections \autoref{eq:piPidef} (we will also left superscripts to denote the inclusion/projections relative to $\mathcal M$ accordingly). Let $F \in \Gamma \mathcal L(\tau \bbR^d, \tau \bbR^e)$ restrict to an element of $\Gamma \mathcal L(\tau \mathcal M, \tau \mathcal N)$ (this means $F(y,x)$ maps $T_x\mathcal M$ to $T_y\mathcal N$ for $x \in M$, $y \in N$): just as in the intrinsic setting the expression $\dif \bfY^k = F^k_c(Y,X) \dif \bfX^c$ is ill-defined, in the extrinsic setting $Y$ will, for $\bfX$ non-geometric, exit $\mathcal M$ when the equation is started on $\mathcal M$. Proceeding heuristically to derive the extrinsic counterpart to the local RDE formula \autoref{eq:RDEM}, with the idea that $\dif_\mathcal{P} \boldsymbol Z = P_c(Z) \dif \boldsymbol Z^c$ for an embedded manifold $\mathcal P$ and $\boldsymbol Z \in \mathscr C^p_\omega([0,T],\mathcal P)$, we interpret 
\begin{equation}
\dif_\mathcal{N} \bfY^k = F^k_c(Y,X) \dif_\mathcal{M} \bfX^c
\end{equation}
as ${^\mathcal{N}\!}P^k_h(Y) \dif \bfY^h = F^k_d(Y,X) {^\mathcal{M}\!}P^d_c(X)\dif \bfX^c$ or in Davie form (using $\bfX \in \mathscr C^p_\omega([0,T],\mathcal M)$, imposing $\bfY \in \mathscr C^p_\omega([0,T],\mathcal N)$ and using \autoref{eq:partialP}) as
\begin{equation}
\begin{split}
{^\mathcal{N}\!}P^k_h(Y_s) Y_{st}^h + \partial_{ij} {^\mathcal{N}\!}\Pi^k(Y_s) \bbY^{ij}_{st} &\approx F_d^k(Y,X) {^\mathcal{M}\!}P^d_c(X) \bfX^c_{st}\\
\bbY^{ij}_{st} &\approx F^i_d F^j_b(Y_s,X_s) \bbX^{ab}_{st}
\end{split}
\end{equation}
Note that we have chosen not to expand the RHS of the first line into first and second-order parts. Now, by \autoref{eq:firstConstr} applied to $Y$, we may rewrite this as 
\begin{equation}
Y^k_{st} - \tfrac 12 \partial_{ij} {^\mathcal{N}\!}\Pi^k(Y_s) Y^i_{st} Y^j_{st} + \partial_{ij} {^\mathcal{N}\!}\Pi^k(Y_s) \bbY^{ij}_{st} \approx F_d^k(Y,X) {^\mathcal{M}\!}P^d_c(X) \bfX^c_{st}
\end{equation}
or as \autoref{eq:extrinsicRDE} in the definition below.
\begin{defn}[Constrained RDE]\label{def:constrRDE}
	Given $\bfX \in \mathscr C^p_\omega([0,T],\mathcal M)$, $y_0 \in \mathcal N$ and $F \in \Gamma \mathcal L(\tau \bbR^d, \tau \bbR^e)$ which restricts to an element of $\Gamma \mathcal L(\tau \mathcal M, \tau \mathcal N)$ we will write
	\begin{equation}
	\dif_\mathcal{N} \bfY^k = F^k_c(Y,X) \dif_\mathcal{M} \bfX^c, \quad Y_0 = y_0
	\end{equation}
	to mean
	\begin{equation}\label{eq:extrinsicRDE}
	\dif \bfY^k = F_d^k(Y,X) {^\mathcal{M}\!}P^d_c(X) \dif \bfX^c + \tfrac 12 \partial_{ij} {^\mathcal{N}\!}\Pi^k(Y) F^i_aF^j_b(Y,X) \dif [\bfX]^{ab},\quad Y_0 = y_0
	\end{equation}
	and say that $\bfY$ solves the $\mathcal N$-\emph{constrained RDE} driven by the $\mathcal M$-constrained rough path $\bfX$.
\end{defn}

The next proposition legitimises this formula.
\begin{thm}\label{thm:exIntRDE}
	Let $\bfX, y_0, F$ be as in \autoref{def:constrRDE}. 
	\begin{enumerate}
		\item The solution to \autoref{eq:extrinsicRDE} only depends on $(F(y,x)|_{T_x\mathcal M})_{x \in \mathcal M, y \in \mathcal N}$ and belongs to $\mathscr C^p_\omega([0,T],\mathcal N)$;
		\item If $\bfX$ is geometric, so is $\bfY$ and the equation can be rewritten as $
		\edif \bfY^k = F^k_c(Y,X) \edif \bfX^c$;
		\item $\bfY \in \mathscr C^p_\omega([0,T],\bbR^e)$ satisfies \autoref{eq:extrinsicRDE} if and only if ${^\mathcal{N}\!}\pi_*\bfY$ solves the RDE driven by $\boldsymbol W \coloneqq {^\mathcal{M}\!}\pi_*\bfX$
		\begin{equation}\label{eq:intExtRDE}
		\edif \boldsymbol Z = (T{^\mathcal{N}\!}\pi \circ F( {^\mathcal{N}\!}\imath(Z), {^\mathcal{M}\!}\imath(W) ) \circ T {^\mathcal{M}\!}\imath) \edif \boldsymbol W
		\end{equation}
		in the sense of \autoref{def:RDEM}.
	\end{enumerate}
	\begin{proof}
		In this proof we will draw on the entirety of the theory of \autoref{subsec:embedded} and the present section, and will therefore omit the precise equations which motivate our computations. The first part of 1.\ will be automatically proved once we show 3.; we therefore proceed to show that $\bfY$ is $\mathcal N$-constrained. We have, omitting all evaluations at $Y_s$ and $X_s$ and relying on indices to distinguish maps referring to $\mathcal M$ and $\mathcal N$ (e.g.\ $P^c_d \coloneqq {^\mathcal{M}\!}P^c_d(X_s)$, $\partial_{ij} \Pi^k \coloneqq \partial_{ij}{^\mathcal{N}\!}\Pi^k (Y_s)$) 
		\begin{equation}\label{eq:RDEM1}
		\begin{split}
		&\mathrel{\phantom{=}}P^k_h Y^h_{st} + \tfrac 12 \partial_{ij} \Pi^k Y^i_{st} Y^j_{st} \\
		&\approx P^k_h \big[F^h_d P^d_c X^c_{st} + (\partial_a F^h_c P^c_b + F^l_a \partial_l F^h_c P^c_b + F^h_c \partial_a P^c_b)\bbX^{ab}_{st} + \tfrac 12 \partial_{ij} \Pi^k F^i_a F^j_b [\bfX]^{ab} \big] \\
		&\mathrel{\phantom{=}} + \tfrac 12 \partial_{ij} \Pi^k F^i_a F^j_b X^a_{st} X^b_{st}
		\end{split}
		\end{equation}
		We calculate
		\begin{equation}
		\begin{split}
		P^k_h F^h_d P^d_c &= F^k_d P^d_c \\
		P^k_h (\partial_a F^h_c P^c_b + F^h_c \partial_a P^c_b) &= \partial_a (F^k_c P^c_b) - Q^k_h \partial_a(F^h_c P^c_b) \\
		&= \partial_a (F^k_c P^c_b) - \partial_a(Q^k_h F^h_c P^c_b) \\
		&= \partial_a (F^k_c P^c_b) \\
		&= \partial_a F^k_c P^c_b + F^k_c \partial_a P^c_b \\
		P^k_h F^l_a \partial_l F^h_c P^c_b &= F^l_a \partial_l(P^k_h F^h_c P^c_b) - F^l_a \partial_l P^k_h F^h_c P^c_b \\
		&= F^l_a \partial_l (F^k_c P^c_b) -  \partial_{ij} \Pi^k F^i_a F^j_b \\
		&= F^l_a \partial_l F^k_c P^c_b -\partial_{ij} \Pi^k F^i_a F^j_b \\
		P^k_h \partial_{ij} \Pi^h F^i_a F^j_b [\bfX]^{ab}_{st} &\approx P^k_h \partial_{ij} \Pi^h F^i_c F^j_d P^c_a P^d_b [\bfX]^{ab}_{st} \\
		&\approx P^k_h \partial_{ij} \Pi^h P^i_l P^j_p F^l_c F^p_d P^c_a P^d_b [\bfX]^{ab}_{st} \\
		&\approx 0
		\end{split}	
		\end{equation}
		Substituting these in \autoref{eq:RDEM1}
		\begin{equation}
		\begin{split}
		&\mathrel{\phantom{=}}P^k_h Y^h_{st} + \tfrac 12 \partial_{ij} \Pi^k Y^i_{st} Y^j_{st} \\
		&\approx F^k_d P^d_c X^c_{st} + ( \partial_a F^k_c P^c_b + F^k_c \partial_a P^c_b + F^l_a \partial_l F^h_c P^c_b ) \bbX^{ab}_{st} + \tfrac 12 \partial_{ij} \Pi^k F^i_a F^j_b(X^a_{st} X^b_{st} - 2 \bbX^{ab}_{st}) \\
		&\approx Y^k_{st}
		\end{split}
		\end{equation}
		
		To prove 2.\ we proceed in a similar fashion: if $\bfX$ is geometric, we have
		\begin{equation}
		\begin{split}
		Y^k_{st} &\approx F^k_d P^d_c X^c_{st} + ( \partial_a F^k_c P^c_b + F^k_c \partial_a P^c_b + F^l_a \partial_l F^h_c P^c_b ) \bbX^{ab}_{st} \\
		&\approx F^k_d(X^d_{st} - \tfrac 12 \partial_{ab} \Pi^d X^a_{st} X^b_{st}) + (\partial_a F^k_b + F^l_a \partial_l F^h_b) \bbX^{ab}_{st} + F^k_d \partial_{ab} \Pi^d (\tfrac 12 X^a_{st} X^b_{st}) \\
		&\approx F^k_d X^d_{st} + (\partial_a F^k_b + F^l_a \partial_l F^h_b) \bbX^{ab}_{st}			\end{split}
		\end{equation}
		and $\bfY$ is geometric because it is the solution to an RDE driven by an $\bbR^d$-valued geometric rough path.
		
		The proof of 3.\ is analogous to that of \autoref{thm:exInt} and therefore omitted.\todo{Is it ok to omit it?}
	\end{proof}
\end{thm}

RDEs can be used to generate elements of $\mathscr C^p_\omega([0,T], \mathcal M)$ starting from any unconstrained rough path (cf.\ \cite[Example 4.12, Proposition 4.13]{CDL15} for the geometric case):
\begin{expl}[Projection construction of constrained rough paths]
	Let $\boldsymbol Z \in \mathscr C^p_\omega([0,T],\mathbb R^d)$. Then the solution $\bfX$ to 
	\begin{equation}
	\dif_\mathcal{M} \bfX^k = P^k_c(X) \dif_{\mathbb R^d} \boldsymbol Z^c, \quad X_0 = x_0 \in M
	\end{equation}
	i.e.\
	\begin{equation}
	\dif \bfX^k = P^k_c(X)\dif \boldsymbol Z^c + \tfrac 12 \partial_{ij} \Pi^k P^i_a P^j_b(X) \dif [\boldsymbol Z]^{ab}
	\end{equation}
	belongs to $\mathscr C^p_\omega([0,T],\mathcal M)$ by \autoref{thm:exIntRDE}. Here $P$ and $\Pi$ refer to the embedded manifold $\mathcal M$. Moreover it is checked, using \autoref{eq:firstConstr} and \autoref{eq:secondConstr} that if $\boldsymbol Z \in \mathscr C^p_\omega([0,T],\mathcal M)$ with $Z_0 = x_0$, then $\bfX = \boldsymbol Z$, i.e.\ this defines a projection $\mathscr C^p_\omega([0,T],\bbR^d) \twoheadrightarrow \mathscr C^p_\omega([0,T],\mathcal M)$.
\end{expl}

\section{Parallel transport and Cartan development}\label{sec:par}
In this section we will discuss parallel transport and Cartan development (or \say{rolling without slipping}) along/of non-geometric rough paths on manifolds. The topic has already been addressed in the geometric case (in the extrinsic setting) in \cite{CDL15}; not assuming geometricity however introduces several complications. The literature on It\^o calculus of semimartingales on manifolds also features similar topics, and we shall reference such instances throughout the chapter; however, because of the adjustments that need to be made for the rough path setting, and because of the greater generality with which the theory is approached (even when restricted to semimartingales), the material presented in this chapter will only depend upon \autoref{sec:backRps}, \autoref{sec:backDG} and \autoref{sec:rough}. We will rely on local coordinates for our computations, and will not explore parallel transport and development in the extrinsic context.

We will tackle parallel transport along non-geometric rough paths by first studying the more general case of RDEs with solutions valued in fibre bundles above the manifold in which the driver is valued; we will progressively restrict our attention to more tractable and interesting cases until we reach the case of the horizontal lift, i.e.\ in which the equation is the natural generalisation of the parallel transport equation to non-geometric rough paths; this will then be used to define Cartan (anti)development. We will see that treating non-geometric rough paths entails adding It\^o-type corrections to the classical formulae, and that the terms appearing in the resulting equations will have to satisfy second-order conditions for properties that are usually taken for granted (well-definedness, linearity, metricity) to hold.

More precisely, we will consider an $m$-dimensional smooth manifold $M$ whose tangent bundle is endowed with a linear connection $\nabla$ which we will think of as fixed throughout this section; given a fibre bundle $\pi \colon E \to M$ and a linear connection $\widetilde \nabla$ on $\tau E$ (note we do not require a connection on the bundle $\pi$), we are interested in equations of the form
\begin{equation}\label{eq:fibreRDE}
\dif_{\widetilde \nabla} \bfY = F(Y) \dif_\nabla \bfX, \quad Y_0 = y_0 \in E_o
\end{equation}
where $F$ is a section of the bundle $\mathcal{L}_E(\tau M, \tau E)$ (the $E$ subscript denotes the base space: this means we are dealing with a bundle over $E$, not $E \times M$, i.e.\ the fibre at $y \in E$ is given by $\mathcal{L}(T_{\pi(y)}M, T_yE)$) s.t.\ and $X_0 = o \in M$ is a basepoint on the manifold which will be fixed throughout this section. The first thing to notice is that such equations are not of the form \autoref{def:RDEM}, since $F$ is not defined for all pairs $(y,x) \in E \times M$; we proceed to introduce the tools that are needed to give this type of equation a meaning. Throughout this section we will use the notation in \autoref{subsec:linearconn}, and in particular  \autoref{conv:indices} for indexing vectors based in the total space of fibre bundles.

The first thing we require of $F$ is that
\begin{equation}\label{eq:Fab}
T_y\pi \circ F(y) = \mathbbm 1_{T_xM} \quad \Longleftrightarrow \quad \delta^\alpha_\beta = (T_y\pi \circ F(y))^\alpha_\beta = (T_y\pi)^\alpha_K F^K_\beta(y) = \delta^\alpha_K F^K_\beta = F^\alpha_\beta(y)
\end{equation}
\textbf{We will assume this condition to hold throughout this section} unless otherwise stated. For $W \in \Gamma \tau M$ we define
\begin{equation}
FW \in \Gamma \tau E, \quad (FW)(y) \coloneqq F(y)W(\pi(y))
\end{equation}
In this section we will understand all expressions as being evaluated at $(x,y)$ with $y \in E_x$ unless otherwise specified. The following definition will be of importance in the study of non-geometric RDEs on fibre bundles:
\begin{defn}
	We define $\widetilde{F} \coloneqq \widetilde{F}(\widetilde \nabla, F)$ by
	\begin{equation}\label{eq:defTwidleF}
	\langle \widetilde{F}, U\otimes V\rangle \coloneqq F \nabla_U V - \widetilde \nabla_{FU}FV \in TE, \quad \text{for } U,V \in \Gamma \tau M
	\end{equation}
\end{defn}

\begin{lem}\label{lem:nablaTwidleF}
	For $U, V \in \Gamma \tau M$ we have
	\begin{equation}
	\begin{split}
	(\widetilde \nabla_{FU}FV)^\gamma &= U^\alpha \partial_\alpha V^\gamma + U^\alpha V^\beta \widetilde \Gamma^\gamma_{\alpha\beta} + F^i_\alpha U^\alpha V^\beta \widetilde \Gamma^\gamma_{i\beta} + U^\alpha F^j_\beta V^\beta \widetilde \Gamma^\gamma_{\alpha j} + F^i_\alpha F^j_\beta U^\alpha V^\beta \widetilde \Gamma^\gamma_{ij} \\
	(\widetilde \nabla_{FU}FV)^k &= U^\alpha(\partial_\alpha F^k_\gamma V^\gamma + F^k_\gamma \partial_\alpha V^\gamma) + F^i_\alpha U^\alpha \partial_i F^k_\gamma V^\gamma \\
	&\mathrel{\phantom{=}} + U^\alpha V^\beta \widetilde \Gamma^k_{\alpha\beta} + F^i_\alpha U^\alpha V^\beta \widetilde \Gamma^k_{i\beta} + U^\alpha F^j_\beta V^\beta \widetilde \Gamma^k_{\alpha j} + F^i_\alpha F^j_\beta U^\alpha V^\beta \widetilde \Gamma^k_{ij} 
	\end{split}
	\end{equation}
	so we have $\widetilde{F} \in \Gamma \mathcal L_E(\tau M^{\otimes 2}, \tau E)$, and
	\begin{align}
	\tensor{\widetilde{F}}{_{\alpha\beta}^\gamma}&= \Gamma^\gamma_{\alpha\beta} - (\widetilde \Gamma^\gamma_{\alpha\beta} + F^i_\alpha  \widetilde \Gamma^\gamma_{i\beta} +  F^j_\beta  \widetilde \Gamma^\gamma_{\alpha j} + F^i_\alpha F^j_\beta  \widetilde \Gamma^\gamma_{ij})\label{eq:FtwidleCoords1} \\
	\tensor{\widetilde{F}}{_{\alpha\beta}^k} &= F^k_\gamma \Gamma^\gamma_{\alpha\beta} - (\partial_\alpha F^k_\beta + F^h_\alpha \partial_h F^k_\beta + \widetilde \Gamma^k_{\alpha\beta} + F^i_\alpha \widetilde \Gamma^k_{i\beta} + F^j_\beta \widetilde \Gamma^k_{\alpha j} + \widetilde \Gamma^k_{ij}F^i_\alpha F^j_\beta) \label{eq:FtwidleCoords2}
	\end{align}
	\begin{proof}
		We compute
		\begin{equation}
		\begin{split}
		\partial_I(FV)^K &= \partial_I (F^K_\gamma (V^\gamma \circ \pi)) \\
		&= \partial_I F^K_\gamma V^\gamma + F^K_\gamma \partial_\beta V^\gamma \partial_I \pi^\beta \\
		&= \begin{cases}
		\partial_\alpha V^\gamma &K= \gamma \leq m, \ I = \alpha \leq m\\
		0 &K= \gamma \leq m, \ I = i > m \\
		\partial_\alpha F^k_\gamma V^\gamma + F^k_\gamma \partial_\alpha V^\gamma & K = k > m, \ I = \alpha \leq m \\
		\partial_i F^k_\gamma V^\gamma & K = k > m, \ I = i > m
		\end{cases}
		\end{split}
		\end{equation}
		Substituting these terms in
		\begin{equation}
		(\widetilde \nabla_{FU}FV)^K = (FU)^I \partial_I (FV)^K + (FU)^I(FV)^J\widetilde\Gamma^K_{IJ}
		\end{equation}
		yields the desired expressions. 
		
		We must now show that $\widetilde{F}$ is bilinear, thus legitimising our use of the notation $\langle \widetilde{F}, U \otimes V \rangle$: this is easily done by computing the RHS of \autoref{eq:defTwidleF} thanks to the previously computed expression, and seeing that the derivatives of $V$ cancel out, leaving us with the desired expressions for $\tensor{\widetilde{F}}{_{\alpha\beta}^\gamma}$ and $\tensor{\widetilde{F}}{_{\alpha\beta}^k}$.
	\end{proof}
\end{lem}
The task is now to extend $F$ to all pairs $(y,x)$ where $y$ does not necessarily lie in $E_x$ (the existence of such extensions is proven in \cite[Lemma 8.16, Proof of Proposition 8.15]{E89}), and to investigate when the resulting \autoref{eq:fibreRDE}, which can now be understood as in \autoref{def:RDEM}, is independent of the extension. To do so we introduce the following condition on $\widetilde \nabla$ and $F$ (with $\nabla$ thought of as fixed).
\begin{cond}\label{cond:nabla}
	Assuming \autoref{eq:Fab} holds, for all $U \in \Gamma \tau M$ we have
	\begin{equation}
	T \pi \widetilde \nabla_{F U} (FU) = \nabla_{U}U
	\end{equation}
	i.e.\ $T\pi \langle \widetilde{F}, U \otimes U \rangle = 0$.
\end{cond}
We have purposefully stated the condition with two copies of the same vector field $U$, instead of $T \pi \widetilde \nabla_{F U} (FV) = \nabla_{U}V$: this is motivated by the fact that we only need the symmetrisation of the latter identity, since these terms will turn out to be the coefficients of the bracket. Requiring the unsymmetrised version of the condition would perhaps have been more natural, but is not as sharp; similar comments hold for \autoref{cond:lin} and \autoref{cond:metric} below.
\begin{lem}\label{lem:condition}
	\autoref{cond:nabla} is equivalent to $\widetilde F$ restricting to an element of $\Gamma \mathcal L_E(\tau M^{\odot 2}, V\pi)$, or in product coordinates to $\tensor{\widetilde F}{_{(\alpha\beta)}^\gamma} = 0$, i.e.\
	\begin{equation}\label{eq:condCoords}
	\widetilde \Gamma^\gamma_{\alpha\beta} + \widetilde \Gamma^\gamma_{i\beta}F^i_\alpha + \widetilde \Gamma^\gamma_{\alpha j} F^j_\beta + \widetilde \Gamma^\gamma_{ij}F^i_\alpha F^j_\beta \stackrel{(\alpha\beta)}{=} \Gamma^\gamma_{\alpha\beta}
	\end{equation}
	Moreover, the condition being satisfied for all choices of $F$ is equivalent to $\tau M$ being symmetrically affine w.r.t.\ $\widetilde \nabla, \nabla$ (with the \say{only if} statement only valid for $m \geq 2$).
	\begin{proof}
		The first characterisation of \autoref{cond:nabla} is obvious, and the expression in local coordinates is a direct consequence of \autoref{lem:nablaTwidleF} and polarisation.
		
		As for the second statement, we must check that the conditions on the symmetrised Christoffel symbols stated in \autoref{expl:twidleAffine} hold. The \say{if} part is immediate. For the converse, first of all reading the identity with $F = 0$ yields $\widetilde \Gamma^\gamma_{\alpha\beta} \stackrel{(\alpha\beta)}{=} \Gamma^\gamma_{\alpha\beta}$, and the identity may be rewritten as
		\begin{equation}
		(\widetilde \Gamma^\gamma_{\alpha j} + \widetilde \Gamma^\gamma_{j\alpha})  F^j_\beta + (\widetilde \Gamma^\gamma_{i \beta} + \widetilde \Gamma^\gamma_{ \beta i}) F^i_\alpha + (\widetilde \Gamma^\gamma_{ij}  + \widetilde \Gamma^\gamma_{ji} ) F^i_\alpha F^j_\beta = 0
		\end{equation}
		Now read the identity for arbitrary but fixed $\alpha \neq \beta$ (which is possible since $m \geq 2$) and $j$, and with $F^k_\gamma \coloneqq \delta^{kj}\delta_{\gamma\beta}$ (this is possible since the coefficients $F^k_\gamma$ are completely arbitrary: we do not even have to argue coordinate-independence, as everything is local and we may take $F$ to be supported in the domain of the chart): this yields $\widetilde \Gamma^\gamma_{\alpha j} \stackrel{(\alpha j)}{=} 0$, reducing our identity to $(\widetilde \Gamma^\gamma_{ij}  + \widetilde \Gamma^\gamma_{ji} ) F^i_\alpha F^j_\beta = 0$. We may then fix arbitrary $i,j$ and pick $F$ exactly as above to conclude $\Gamma^\gamma_{ij} \stackrel{(ij)}{=} 0$.
	\end{proof}
\end{lem}
The next result establishes the link between the condition and well-definedness of RDEs in fibre bundles.
\begin{thm}\label{thm:wellDefRDE}
	If \autoref{cond:nabla} holds \autoref{eq:fibreRDE} is well-defined, i.e.\ it is independent of the extension of $F$ to an element of $\Gamma \mathcal L_{E \times M}(\tau M, \tau E)$. In this case $\pi(Y) = X$ and the coordinate expression of the RDE reduces to its vertical component and is given by
	\begin{equation}\label{eq:wdFibreRDE}
	\begin{split}
	\edif \bfY^k &= F^k_\gamma(Y) \edif \bfX^\gamma \\
	&\mathrel{\phantom{=}}+ \tfrac 12 (F^k_\gamma(Y)\Gamma^\gamma_{\alpha\beta}(X)  - (\widetilde\Gamma^k_{\alpha\beta} + \widetilde\Gamma^k_{\alpha j}F^j_\beta + \widetilde\Gamma^k_{i\beta} F^i_\alpha + \widetilde\Gamma^k_{ij}F^i_\alpha F^j_\beta)(Y) )\edif [\bfX]^{\alpha\beta} 
	\end{split}
	\end{equation}
	which may be written at the trace level as
	\begin{equation}
	\edif Y^k = F^k_\gamma(Y) \circ \edif \bfX^\gamma + \tfrac 12 \tensor{\widetilde{F}}{_{\alpha\beta}^k}(Y) \edif[\bfX]^{\alpha\beta} 
	\end{equation}
	Moreover, if $\bfX$ is geometric the equation is always well-defined and independent of the connections $\nabla$ and $\widetilde \nabla$.
	\begin{proof}
		Denoting still with $F = F(y,x)$ an arbitrary extension of $F$ as a section of the bundle $\mathcal{L}_{E \times M}(\tau M, \tau E)$ (i.e.\ $F(y,\pi(y)) = F(y)$), we have that the first $m$ coordinates of the local form of \autoref{eq:fibreRDE} is given by
		\begin{equation}\label{eq:fibreRDECoords}
		\dif \bfY^\gamma = F^\gamma_\alpha(Y,X) \dif \bfX^\alpha + \tfrac 12 (\Gamma^\gamma_{(\alpha\beta)}(X)  - \widetilde\Gamma^\gamma_{(IJ)}(Y) F^I_\alpha F^J_\beta (Y,X) )\dif [\bfX]^{\alpha\beta}
		\end{equation}
		where we have symmetrised the second order part thanks to the symmetry of the tensor $[\bfX]$. Notice that by \autoref{eq:Fab}, by hypothesis and \autoref{lem:condition} we have that on pairs $(Y,X)$ s.t.\ $Y \in T_XM$ the coefficient of $\dif \bfX^\alpha$ equals $\delta^\alpha_\gamma$ and that of $\dif [\bfX]^{\alpha\beta}$ vanishes. Now consider the RDE defined only locally in the domain of the chart (i.e.\ without the claim that the following is a coordinate-invariant expression)
		\begin{equation}
		\dif \begin{pmatrix}\bfY^\gamma \\ \bfY^k\end{pmatrix} = \begin{pmatrix}\dif\bfX^\gamma \\ F^k_\gamma(Y,X) \dif \bfX^\gamma + \tfrac 12 (F^k_\gamma(Y,X) \Gamma^\gamma_{\alpha\beta}(X)  - \widetilde\Gamma^k_{IJ}(Y) F^I_\alpha F^J_\beta(Y,X) )\dif [\bfX]^{\alpha\beta} \end{pmatrix}
		\end{equation}
		The solution to this RDE stays in the fibre of the trace of the driver $X$. But the solution to this RDE must also solve \autoref{eq:fibreRDECoords}, since it takes its values in the locus in which the coefficients of the two coincide (here we are using the obvious principle that the solution of an RDE does not change if the coefficients are modified away from the solution).
		
		If $\bfX$ is geometric $[\bfX]$ vanishes altogether and we may show well-definedness in the same manner, and is independent of the connections on the source and target manifolds since the driver is geometric.
	\end{proof}
\end{thm}
\textbf{For the remainder of this section we assume \autoref{cond:nabla} is satisfied} unless otherwise stated. An even stronger requirement (which we will instead not assume to hold) would be:
\begin{cond}\label{cond:itoEqStrat}
	$F \nabla_U U = \widetilde \nabla_{FU}FU$ for all $U \in \Gamma \tau M$, i.e.\ $\widetilde F = 0$, or in local coordinates (assuming \autoref{cond:nabla} already holds)
	\begin{equation}
	F^k_\gamma \Gamma^\gamma_{\alpha\beta} \stackrel{(\alpha\beta)}{=} \partial_\alpha F^k_\beta + F^h_\alpha \partial_h F^k_\beta + \widetilde \Gamma^k_{\alpha\beta} + F^i_\alpha \widetilde \Gamma^k_{i\beta} + F^j_\beta \widetilde \Gamma^k_{\alpha j} + \widetilde \Gamma^k_{ij}F^i_\alpha F^j_\beta
	\end{equation}
\end{cond}
The following is immediately inferred through \autoref{thm:wellDefRDE}, although it can also be derived indirectly by extending $F$ to a section of $\mathcal{L}_{E \times M}(\tau M, \tau E)$ and applying \autoref{eq:twoRDEs}.
\begin{cor}\label{cor:ptStrat}
	If \autoref{cond:itoEqStrat} holds then the trace-level of \autoref{eq:fibreRDE} is well-defined and equivalent to
	\begin{equation}
	\edif Y = F(Y) \circ \edif \bfX
	\end{equation}
\end{cor}
We will mostly use the geometrised form of our fibre bundle-valued RDEs, regardless of whether \autoref{cond:itoEqStrat} holds or not, keeping in mind that the second-order part of the solution is given in terms of the original rough path $\bfX$ as $\bbY^{ij}_{st} \approx F^i_a F^j_b(Y_s) \bbX^{ab}_{st}$.

In the following example we show how it is possible to generate integrands by solving cotangent bundle-valued RDEs.
\begin{expl}[Solutions to RDEs as controlled integrands]
	If $\pi = \tau^*M$ we may make the solution to \autoref{eq:fibreRDE} into an element of $\mathscr D_X(\tau^*M) = \mathscr D_X(\mathcal L(\tau M,\bbR))$ (on any interval $[0,S]$ on which the solution is defined) as follows. Working in induced coordinates, it is natural to define $Y'_{\alpha \beta} \coloneqq F^{\widetilde \beta}_\alpha(Y)$: since the resulting $\bfY$ is a controlled path in each coordinate chart, we need only check that it satisfies the transformation rule required in \autoref{def:contrIntM} to conclude it is a controlled integrand. At the path level this is already guaranteed by the soundness of \autoref{def:RDEM}, but for Gubinelli derivatives it must be verified, as we have not defined solutions to RDEs on manifolds as controlled paths: by \autoref{eq:changeCotanCoords} we have
	\begin{equation}
	\begin{split}
	Y'_{\overline \alpha \overline \beta} &= F^{\widetilde{\overline \beta}}_{\overline \alpha}(Y) \\
	&= \partial^{\widetilde{\overline \beta}}_B F^B_\alpha(Y) \partial^\alpha_{\overline \alpha} \\
	&= \partial^{\widetilde{\overline \beta}}_\beta \delta^\beta_\alpha \partial^\alpha_{\overline \alpha} + \partial^{\widetilde{\overline \beta}}_{\widetilde \beta} F^{\widetilde \beta}_\alpha(Y) \partial^\alpha_{\overline \alpha} \\
	&= \partial^\gamma_{\overline \alpha \overline \beta} Y_\gamma +  \partial^\alpha_{\overline \alpha} \partial^\beta_{\overline \beta} Y'_{\alpha\beta}
	\end{split}
	\end{equation}	which is precisely the transformation rule required to be satisfied by Gubinelli derivatives under \autoref[Pullback]{prop:contrProp}. We have thus given a meaning to the scalar rough integral $\int \bfY \dif \bfX$. Examples of connections on the cotangent bundle, analogous to the complete, horizontal and Sasaki lifts, defined below\todo{How are these lifts related to the ones to the tangent bundle? Was not able to understand this.} for the tangent bundle, are given in \cite[p.269, p.286]{YI73} and \cite{Sal11}; it would be interesting to verify that \autoref{cond:nabla} holds for these connections. If we start with the bundle $\mathcal L(\tau M, \bbR^e) = (\tau^*M)^{\oplus e}$, we may similarly define $\bbR^e$-valued rough integrals (here we need a connection on $\tau \mathcal L(T M, \bbR^e)$\todo{I think if we have one on $\tau T^*M$ we should get one on $\tau \mathcal L(T M, \bbR^e)$ for free imposing as many cross terms as possible to be zero. More generally, if we have $n$ different connections on $\tau T^*M$ satisfying some compatibility condition they should give rise to a connection on $\tau \mathcal L(T M, \bbR^e)$. It would be nice to make these comments precise if they are not too difficult to justify.}). It is similarly checked that taking $\pi = \tau M$ (as done below) results in the solution of \autoref{eq:fibreRDE} being a $\tau M$-valued controlled path \autoref{expl:rri} (with Gubinelli derivatives $Y'^\beta_\alpha = F^{\widetilde \beta}_\alpha$), which can then be integrated against $\bfX$ once a Riemannian metric on $\tau M$ is given.
\end{expl}

We proceed with the main theory. \textbf{For the remainder of this section we let $\boldsymbol\pi$ be a vector bundle} unless otherwise stated. We will say that $\widetilde U \in \Gamma \tau E$ is \emph{linear} if
\begin{equation}\label{eq:linearVF}
\widetilde U^\gamma(y) = U^\gamma(x), \quad \widetilde U^k(y) = \widetilde U^k_h(x) y^h
\end{equation}
with $x = \pi(y)$ and for locally defined functions $U^\gamma, \widetilde U^k_h$. We will say that $F$ (which we are assuming satisfies \autoref{eq:Fab} and \autoref{cond:nabla}) is \emph{linear} if $FU$ is a linear vector field for all $U \in \Gamma \tau M$.\todo[backgroundcolor=yellow]{Is there a nicer, more intrinsic way to state this linearity condition (without relying on the connections)? The problem is that the target space of the linear map with coords.\ $F^k_{\gamma h}$ changes with the point ($T_yE$)...}
\begin{lem}\label{lem:linearVF}
	The condition of $U \in \Gamma \tau E$ of being linear is independent of the coordinate system. Therefore that of $F$ of being such is too, and in coordinates it amounts to
	\begin{equation}
	F^k_\gamma(y) = F^k_{\gamma h}(x)y^h
	\end{equation}
	for locally defined functions $F^k_{\gamma h}$.
	\begin{proof}
		With the notations of \autoref{eq:changeProdCoords} we have 
		\begin{equation}
		\begin{split}
		\widetilde U^{\overline \gamma}(y) &= \partial^{\overline \gamma}_K(x) \widetilde U^K(y) = \partial^{\overline \gamma}_\gamma U^\gamma(x) \\
		\widetilde U^{\overline k}(y) &= \partial^{\overline k}_K (x) \widetilde U^K(y) \\
		&= \partial_\gamma \lambda^{\overline k}_k(x) y^k U^\gamma(x) + \lambda^{\overline k}_k(x) \widetilde U^k_h(y)y^h \\
		&= (\partial_\gamma \lambda^{\overline k}_k(\lambda^{-1})^k_{\overline h} U^\gamma + \lambda^{\overline k}_k \widetilde U^k_h (\lambda^{-1})^h_{\overline h})(x) y^{\overline h}
		\end{split}
		\end{equation}
		and we may therefore set $U^{\overline\gamma} = \partial^{\overline \gamma}_\gamma U^\gamma$, $\widetilde U^{\overline k}_{\overline h} = \partial_\gamma \lambda^{\overline k}_k(\lambda^{-1})^k_{\overline h} U^\gamma + \lambda^{\overline k}_k \widetilde U^k_h (\lambda^{-1})^h_{\overline h}$. The linearity condition on $F$ is not stated with reference to a particular coordinate system, and is therefore invariant under change of coordinates because linearity of vector fields is. Finally picking an arbitrary $U \in \Gamma \tau M$, we have $(FU)^\gamma(y) = U^\gamma(x)$ by \autoref{eq:Fab} and for $(FU)^k(y) = F^k_\gamma(y) U^\gamma(x)$ to be of the form $\widetilde U^k_h(x) y^h$ for all $U$ we need $F^k_\gamma(y)$ to be of the form $F^k_{\gamma h}(x) y^h$ (for the \say{only if} implication simply pick $U^\gamma = \delta^\gamma_\beta$ with $\beta = 1,\ldots, m$).
	\end{proof}
\end{lem}
Note that the $k$ index in $F^k_{\gamma h}(x)y^h$ represents a coordinate in $T_yE$, whereas $h$ represents a coordinate in $E_x$; following \autoref{conv:indices} we will not place a twidle on the upper index, as we view $F^k_{\gamma h}$ as the coordinates of a linear map between vector spaces. For the remainder of this section we will assume $F$ is linear, and we will be concerned with the question of whether this implies that the resulting \autoref{eq:fibreRDE} is also linear, i.e.\ that its coordinate expression is linear in the conventional sense. To this end, we introduce the following condition on $F$ and $\widetilde \nabla$ (which does not involve the connection $\nabla$ at all, to the extent that $\widetilde \nabla$ is not defined in terms of it).
\begin{cond}\label{cond:lin}
	Assume $F$ is linear. $\widetilde\nabla_{FU}(FU)$, or equivalently $\langle \widetilde F,U\otimes U\rangle$, is a linear vector field for all $U \in \Gamma \tau M$.
\end{cond}

\begin{lem}\label{lem:coordLin}
	Let $F$ be linear and \autoref{cond:nabla} hold. Then \autoref{cond:lin} is equivalent to $\widetilde F|_{TM^{\odot 2}}$ lying in the image of the map
	\begin{equation}\label{eq:imageMap}
	\begin{split}
	\Gamma (\tau^*M^{\odot 2} \otimes \pi^* \otimes \pi) = \Gamma\mathcal L_M(\tau M^{\odot 2} \otimes \pi, \pi) &\to \Gamma \mathcal L_E (\tau M^{\odot 2},V\pi) \\
	G &\mapsto \big(e \mapsto (U \odot V \mapsto \mathscr v(e)\langle G, U \odot V \otimes e \rangle ) \big)
	\end{split}
	\end{equation}
	where $\mathscr v(e) \colon E \to V_e\pi$ denotes the vertical lift isomorphism based at $e$. In other words, we may write its coordinates (symmetrising in the first two indices) as $\tensor{\widetilde F}{_{(\alpha\beta) h}^\gamma} = 0$ and
	\begin{equation}\label{eq:G}
	\begin{split}
	\tensor{\widetilde F}{_{\alpha\beta h}^k} y^h &\stackrel{(\alpha\beta)}{=} F^k_{\gamma h} \Gamma^\gamma_{\alpha\beta} y^h - (\partial_\alpha F^k_{\beta h}y^h + F^l_{\alpha h} F^k_{\beta l} y^h\\
	&\mathrel{\hphantom{\stackrel{(\alpha\beta)}{=}}} + \widetilde \Gamma^k_{ij}F^i_{\alpha h} F^j_{\beta l} y^hy^l + \widetilde \Gamma^k_{\alpha\beta} + F^i_{\alpha h} \widetilde \Gamma^k_{i\beta}y^h + F^j_{\beta h} \widetilde \Gamma^k_{\alpha j} y^h)
	\end{split}
	\end{equation}
	and it follows that an equivalent formulation of the condition is that the expression
	\begin{equation}\label{eq:coordLin}
	\widetilde \Gamma^k_{\alpha\beta} + \widetilde \Gamma^k_{\alpha j} F^j_{\beta h} y^h + \widetilde \Gamma^k_{i\beta} F^i_{\alpha h} y^h + \widetilde \Gamma^k_{ij} F^i_{\alpha h} F^j_{\beta l} y^h y^l, \quad (\alpha\beta)
	\end{equation}
	is linear in the $y$ coordinates.
	
	Moreover, the condition being satisfied for all choices of $F$ as above (without assuming \autoref{cond:nabla} is) is equivalent to the stronger requirement that $\widetilde \nabla_{\widetilde U} \widetilde U$ be linear for all linear $\widetilde U \in \Gamma \tau E$ (with the \say{only if} statement only valid for $m \geq 2$), which in coordinates reads
	\begin{equation}\label{eq:linConstGamma}
	\begin{split}
	&\widetilde \Gamma^\gamma_{(\alpha\beta)} \text{ constant in $y$}, \quad \widetilde \Gamma^\gamma_{(i\beta)} = \widetilde \Gamma^\gamma_{(\alpha j)} = \widetilde \Gamma^\gamma_{(ij)}= 0 \\
	&\widetilde \Gamma^k_{(\alpha\beta)} \text{ linear in $y$}, \quad \widetilde \Gamma^k_{(\alpha j)}, \widetilde \Gamma^k_{(i\beta)} \text{ constant in $y$}, \quad \widetilde \Gamma^k_{(ij)} = 0
	\end{split} 
	\end{equation}
\end{lem}
Note that in \autoref{eq:G} we are not able to provide an expression for $\tensor{\widetilde F}{_{\alpha\beta h}^k}$, since some of the terms on the RHS are nonlinear (recall that the $\Gamma^\gamma_{\alpha\beta}$'s and $F^k_{\gamma h}$ are evaluated at $x$, but the $\widetilde \Gamma^K_{IJ}$ are non-linearly evaluated at $y$, and moreover there are quadratic terms).

\begin{proof}[Proof of \autoref{lem:coordLin}]
	The first characterisation is just a reformulation of the second, which is evident by \autoref{eq:FtwidleCoords2}, \autoref{lem:condition} and the definition of linear vector field. The third follows from the second by subtracting terms that are already linear in $y$.
	
	As for the second statement, we first observe that linearity of $\widetilde \nabla_{FU}FU$ without requiring \autoref{cond:nabla} entails the additional requirement that (by \autoref{eq:FtwidleCoords1}, rewritten to account for the linearity of $F$) the expression 
	\begin{equation}\label{eq:Fconst}
	\Gamma^\gamma_{\alpha\beta} - (\widetilde \Gamma^\gamma_{\alpha\beta} + F^i_{\alpha h}  \widetilde \Gamma^\gamma_{i\beta}y^h +  F^j_{\beta h} \widetilde \Gamma^\gamma_{\alpha j}y^h + F^i_{\alpha h} F^j_{\beta l}  \widetilde \Gamma^\gamma_{ij}y^h y^l), \quad (\alpha\beta)
	\end{equation}
	be constant in $y$. Then by arguing as in the proof of \autoref{lem:condition} by progressively disregarding constant (resp.\ linear) terms in \autoref{eq:Fconst} (resp.\ \autoref{eq:coordLin}) we may conclude that linearity of $\widetilde \nabla_{FU}FU$ for all $F$ and $U$ as above is equivalent to \autoref{eq:linConstGamma}.
	
	Now, writing $(\nabla_{\widetilde U} \widetilde V)^K = \widetilde U^I \partial_I \widetilde V^K + \widetilde U^I \widetilde V^J \widetilde \Gamma^K_{IJ}$ for $\widetilde U, \widetilde V$ linear (with notation as in \autoref{eq:linearVF}) we obtain
	\begin{equation}
	\begin{split}
	(\nabla_{\widetilde U} \widetilde V)^\gamma &= U^\alpha \partial_\alpha V^\gamma + U^\alpha V^\beta \widetilde \Gamma^\gamma_{\alpha\beta} + U^\alpha \widetilde V^j_h \widetilde \Gamma^\gamma_{\alpha j} y^h + \widetilde U^i_h  V^\beta \widetilde \Gamma^\gamma_{i \beta} y^h + \widetilde U^i_h  \widetilde V^j_l \widetilde \Gamma^\gamma_{ij} y^h y^l \\
	(\nabla_{\widetilde U} \widetilde V)^k &= U^\alpha \partial_\alpha \widetilde V^k_h y^h + \widetilde U^i_h \widetilde V^k_i y^h + U^\alpha V^\beta \widetilde \Gamma^k_{\alpha\beta} + U^\alpha \widetilde V^j_h \widetilde \Gamma^k_{\alpha j} y^h + \widetilde U^i_h V^\beta \widetilde \Gamma^k_{i\beta} y^h \\
	&\mathrel{\phantom{=}} + \widetilde U^i_h \widetilde V^j_l \widetilde \Gamma^k_{ij} y^h y^l
	\end{split}
	\end{equation}
	As usual, we rely on the symbols involved to infer whether a function is evaluated at $y \in E$ or at $x = \pi(y)$. We then see, by arbitrarity of $U^\gamma, V^\gamma, \widetilde U^k_h, \widetilde V^k_h \in C^\infty M$, polarisation, and the usual elimination procedure, that linearity of $\nabla_{\widetilde U} \widetilde U$ is equivalent to \autoref{eq:linConstGamma}.
\end{proof}

Assuming \autoref{cond:nabla} is satisfied we may consider the \emph{flow map} associated to $F$ and $\bfX$ at times $0 \leq s \leq t \leq T$
\begin{equation}
\begin{split}
&\Phi_{ts} = \Phi(F,\bfX)_{ts} \colon E_{X_s} \to E_{X_t}, \quad y \mapsto Y_t \\
&\text{where} \quad \dif Y = F(Y) \dif \bfX, \ Y_s = y
\end{split}
\end{equation}
which is defined as long as $Y_t$ is defined, and by uniqueness we have
\begin{equation}\label{eq:PhiComp}
\Phi_{tu} \circ \Phi_{us} = \Phi_{ts}
\end{equation}
for $0 \leq s \leq u \leq t \leq T$ whenever one of the two sides is defined. The following theorem justifies our interest in the linearity condition. 
\begin{thm}\label{thm:linearFibreRDE}
	Let $F$ be linear and satisfy \autoref{cond:nabla} and \autoref{cond:lin}. Then \autoref{eq:fibreRDE} can be written in coordinates as
	\begin{equation}\label{eq:linearRDE}
	\begin{split}
	\edif Y^k &= F^k_{\gamma h}(X)Y^h \circ \edif \bfX^\gamma + \tfrac 12 \tensor{\widetilde F}{_{(\alpha\beta)h} ^k}(X)Y^h \edif [\bfX]^{\alpha\beta}
	\end{split}
	\end{equation}
	and admits a global solution. Moreover, $\Phi_{ts}$ defines linear isomorphisms $E_{X_s} \cong E_{X_t}$ for all $0 \leq s \leq t \leq T$. These statements also hold, independently of $\widetilde F$, if $\bfX$ is geometric.
	\begin{proof}
		The first statement is a restatement of \autoref{eq:wdFibreRDE} to the case in which \autoref{cond:lin} is satisfied. We may argue global existence by \autoref{thm:localE} and \autoref{rmk:localEM}: indeed, assume that there exists $S \leq T$ such that $Y_{[0,S)}$ is not contained in any compact set of $M$. Since $\pi(Y) = X$ on $[0,S)$, we must have that $\lim_{t \to S^-} \pi(Y_t) = X_S$, i.e.\ $Y$ must \say{explode vertically}. This, however, is not possible either, since if we may pick a system of product coordinates which contains $X_S$, this would mean that the coordinate solution to \autoref{eq:wdFibreRDE} must only be defined for $t < S$, which is ruled out by \autoref{lem:globE}.
		
		Standard uniqueness arguments apply charts to show that $\Phi_{ts}$ is a linear monomorphism (and thus an isomorphism, by dimensionality) when $X_s,X_t$ are contained in a single chart, and these can be combined to yield the global statement by \say{patching} $X[0,T]$ with finitely many charts and applying \autoref{eq:PhiComp}.
	\end{proof}
\end{thm}
We will denote $\Phi_{st} \coloneqq \Phi_{ts}^{-1}$ for $0 \leq s \leq t \leq T$. We proceed to study the local dynamics satisfied by $t \mapsto \Phi_{t0}$ and $t \mapsto \Phi_{0t}$. Fix coordinates for the vector space $E_o = E_{X_0}$, which we denote with the symbols $i^\circ,j^\circ,k^\circ\ldots$; we continue to denote with $\alpha,\beta,\gamma\ldots$ and $i,j,k\ldots$ the local coordinates in and above a neighbourhood containing $X_t$; we do not intend for the former indices to bear any relationship with the latter (e.g.\ $k^\circ$ and $k$ appearing in a common expression have nothing to do with each other). 
\begin{prop}\label{prop:Phi}
	The coordinate expressions $\Phi^k_{k^\circ;t0}$ and $\Phi^{k^\circ}_{k;0t}$ respectively solve the RDEs (at the trace level) driven by $({_\emph{g}\!}\bfX,[\bfX])$
	\begin{equation}
	\begin{split}
	\edif \Phi^k_{k^\circ;t0} &= F^k_{\gamma h}(X_t)\Phi^h_{k^\circ;t0} \circ \edif \bfX^\gamma_t + \tfrac 12 \tensor{\widetilde F}{_{(\alpha\beta)h}^k}(X_t)\Phi^h_{k^\circ;t0} \edif [\bfX]^{\alpha\beta}_t \\
	\edif \Phi^{k^\circ}_{k;0t} &= -\Phi^{k^\circ}_{h;0t} F^h_{\gamma k}(X_t) \circ \edif \bfX^\gamma_t - \tfrac 12 \Phi^{k^\circ}_{h;0t} \tensor{\widetilde F}{_{(\alpha\beta)k}^h}(X_t) \edif [\bfX]^{\alpha\beta}_t
	\end{split}
	\end{equation}
	\begin{proof}
		The statement is local, and we may confine ourselves to the domain of a single set of product coordinates containing $X_t$. By \autoref{thm:linearFibreRDE} we have
		\begin{equation}
		\begin{split}
		(\dif \Phi^k_{k^\circ;t0})y &= \dif (\Phi^k_{k^\circ;t0} y) \\
		&= \dif Y_t \\
		&= F^k_{\gamma h}(X_t)Y_t^h \circ \dif \bfX_t^\gamma + \tfrac 12 \tensor{\widetilde F}{_{(\alpha\beta)h} ^k}(X_t)Y_t^h \dif [\bfX]^{\alpha\beta}_t \\
		&= F^k_{\gamma h}(X_t)\Phi^k_{h^\circ;t0} y \circ \dif \bfX_t^\gamma + \tfrac 12 \tensor{\widetilde F}{_{(\alpha\beta)h} ^k}(X_t)\Phi^h_{k^\circ;t0} y \dif [\bfX]^{\alpha\beta}_t \\
		&= (F^k_{\gamma h}(X_t)\Phi^k_{h^\circ;t0} \circ \dif \bfX_t^\gamma + \tfrac 12 \tensor{\widetilde F}{_{(\alpha\beta)h} ^k}(X_t)\Phi^h_{k^\circ;t0} \dif [\bfX]^{\alpha\beta}_t)y
		\end{split}
		\end{equation}
		We may therefore conclude, by arbitrarity of $y \in E_o$, that the first of the two RDEs holds. As for the second, we have 
		\begin{equation}
		\begin{split}
		0 &= \dif \delta^{k^\circ}_{h^\circ} \\
		&= \dif (\Phi^{k^\circ}_{k;0t} \Phi^{k}_{h^\circ;t0}) \\
		&= (\dif \Phi^{k^\circ}_{k;0t}) \Phi^{k}_{h^\circ;t0} +  \Phi^{k^\circ}_{k;0t} (F^k_{\gamma h}(X_t)\Phi^h_{h^\circ;t0} \circ \dif \bfX_t^\gamma + \tfrac 12 \tensor{\widetilde F}{_{(\alpha\beta)h} ^k}(X_t)\Phi^h_{h^\circ;t0} \dif [\bfX]^{\alpha\beta}_t)
		\end{split}
		\end{equation}
		which we rewrite as
		\begin{equation}
		\begin{split}
		\dif \Phi^{k^\circ}_{l;0t} &= -\Phi^{k^\circ}_{k;0t} (F^k_{\gamma h}(X_t)\Phi^h_{h^\circ;t0} \circ \dif \bfX_t^\gamma + \tfrac 12 \tensor{\widetilde F}{_{(\alpha\beta)h} ^k}(X_t)\Phi^h_{h^\circ;t0} \dif [\bfX]^{\alpha\beta}_t)\Phi^{h^\circ}_{l;0t} \\
		&= -\Phi^{k^\circ}_{k;0t} F^k_{\gamma l}(X_t)\circ \dif \bfX_t^\gamma -\tfrac 12 \Phi^{k^\circ}_{k;0t} \tensor{\widetilde F}{_{(\alpha\beta)l} ^k}(X_t) \dif [\bfX]^{\alpha\beta}_t
		\end{split}
		\end{equation}
		thus concluding the proof.\todo{Is the level of rigour sufficient, or should I make things clearer (e.g.\ by autoreffing the specific operations used in \autoref{expl:RDEOps}?} 
	\end{proof}
\end{prop}
\textbf{For the remainder of this section we will let $\boldsymbol{\pi = \tau M}$} unless otherwise stated. An important feature of the equation in this case is that we can integrate the inverse of the flow map to obtain a $T_oM$-valued rough path. Note that, although we have not explicitly defined controlled integrands with values in an arbitrary finite-dimensional vector space $V$, this is done simply by choosing a basis of $V$ and setting $\mathscr D_X(\mathcal L(\tau M,V)) \coloneqq \mathscr D_X(\mathcal L(\tau M,\bbR^m))$ under the corresponding isomorphism $V \cong \bbR^m$ (all the needed constructions are easily seen not to depend on the choice of the basis). The next lemma states the change of coordinate formula satisfied by $F^\gamma_{\alpha\beta}$:
\begin{lem}\label{lem:changeF}
	\begin{equation}
	F^{\widetilde{\overline\gamma}}_{\overline \alpha \overline \beta} = \partial_\gamma^{\overline \gamma} \partial_{\overline \alpha}^\alpha \partial_{\overline \beta}^\beta F^\gamma_{\alpha \beta} -  \partial_\gamma^{\overline \gamma}\partial_{\overline \alpha \overline \beta}^\gamma
	\end{equation}
	\begin{proof}
		By \autoref{eq:changeTanCoords} we have
		\begin{equation}
		\begin{split}
		F^{\widetilde \gamma}_{\alpha \beta} \partial^{\beta}_{\overline \beta} y^{\overline \beta} &= F^{\widetilde \gamma}_{\alpha \beta} y^{\beta} \\
		&= F^{\widetilde \gamma}_{\alpha} \\
		&= \partial^{\widetilde{\gamma}}_{\overline C} F^{\overline C}_{\overline\alpha} \partial^{\overline \alpha}_{\alpha} \\
		&=\partial^{\widetilde{ \gamma}}_{\overline \gamma} \delta^{\overline \gamma}_{\overline \alpha} \partial^{\overline \alpha}_{ \alpha} + \partial^{\widetilde{\gamma}}_{\widetilde {\overline \gamma}} F^{\widetilde {\overline \gamma}}_{\overline \alpha} \partial^{\overline \alpha}_{\alpha} \\
		&= \partial^{ \gamma}_{\overline \gamma \overline  \beta} y^{\overline \beta} \delta^{\overline \gamma}_{\overline \alpha} \partial^{\overline \alpha}_{ \alpha} + \partial^{\gamma}_{\overline \gamma} F^{\widetilde{\overline \gamma}}_{\overline \alpha \overline \beta} y^{\overline \beta} \partial^{\overline \alpha}_{ \alpha} \\
		&= (\partial^{\gamma}_{\overline \alpha \overline \beta} \partial^{\overline \alpha}_{\alpha} + \partial^{ \gamma}_{\overline \gamma} F^{\widetilde{\overline \gamma}}_{\overline \alpha \overline \beta} \partial^{\overline \alpha}_{\alpha} ) y^{\overline \beta}
		\end{split}
		\end{equation}
		from which
		\begin{equation}
		F^{\widetilde \gamma}_{\alpha \beta} \partial^{\beta}_{\overline \beta} = \partial^{\gamma}_{\overline \alpha \overline \beta} \partial^{\overline \alpha}_{\alpha} + \partial^{ \gamma}_{\overline \gamma} F^{\widetilde{\overline \gamma}}_{\overline \alpha \overline \beta} \partial^{\overline \alpha}_{\alpha} 
		\end{equation}
		thanks to the arbitrarity of $y$, and we may conclude.
	\end{proof}
\end{lem}
As there will be we no risk of ambiguity, we shall reassign $F^\gamma_{\alpha\beta} \coloneqq F^{\widetilde \gamma}_{\alpha\beta}$, and since now, in view of \autoref{lem:coordLin}, $\widetilde F$ may be viewed as restricting to an element of $\Gamma(\tau^* M^{\odot 2} \otimes \tau^* M \otimes \tau M)$ it also makes sense to set $\tensor{\widetilde F}{_{(\alpha\beta)\gamma}^\delta} \coloneqq \tensor{\widetilde F}{_{(\alpha\beta)\gamma}^{\widetilde \delta}}$. The tensor field $\widetilde F$ may now be given the following interpretation: its evaluation against $(U \odot V) \otimes W$ consists of taking the (symmetrisation of the) defect in commutativity between covariant derivatives and horizontal lift, $\mathscr h \nabla_U V - \widetilde\nabla_{\mathscr h U}\mathscr h V$ and mapping its vertical part at $W \in TM$ down isomorphically onto $TM$.
\begin{prop}\label{prop:intPhi}
	$\boldsymbol {\Phi}_{0\cdot} \in \mathscr D_X(\mathcal L(\tau M, T_oM))$, where $\Phi'^{\gamma^\circ}_{\alpha\beta;0t} \coloneqq -\Phi^{\gamma^\circ}_{\gamma;0t} F^\gamma_{\alpha\beta}(X_t)$.
	\begin{proof}
		The local condition is satisfied in each coordinate chart thanks to \autoref{prop:Phi}. We must check that the compatibility condition of \autoref{def:contrIntM} is met: again, this is obvious at the trace level, and for Gubinelli derivatives we have, by \autoref{lem:changeF}
		\begin{equation}
		\begin{split}
		\Phi'^{\gamma^\circ}_{\overline\alpha \overline\beta} &= -\Phi^{\gamma^\circ}_{\overline \gamma} F^{\overline\gamma}_{\overline\alpha \overline\beta} \\
		&= -\Phi^{\gamma^\circ}_{\overline \gamma} F^{\overline\gamma}_{\overline\alpha \overline\beta} \\
		&= -\Phi^{\gamma^\circ}_\gamma \partial^\gamma_{\overline\gamma} (\partial_\delta^{\overline \gamma} \partial_{\overline \alpha}^\alpha \partial_{\overline \beta}^\beta F^\delta_{\alpha \beta} - \partial_\delta^{\overline \gamma} \partial_{\overline \alpha \overline \beta}^\delta )\\
		&= -\Phi^{\gamma^\circ}_\gamma (\partial_{\overline \alpha}^\alpha \partial_{\overline \beta}^\beta F^\gamma_{\alpha \beta} - \partial^\gamma_{\overline\alpha \overline\beta}) \\
		&= \Phi'^{\gamma^\circ}_{\alpha \beta}\partial_{\overline \alpha}^\alpha \partial_{\overline \beta}^\beta + \Phi^{\gamma^\circ}_\gamma \partial^\gamma_{\overline\alpha \overline\beta}
		\end{split}
		\end{equation}
		Thus concluding the proof.
	\end{proof}
\end{prop}

We now restrict our attention for the last time: \textbf{from now on we will consider the case in which $\boldsymbol{F}$ is given by the horizontal lift $\boldsymbol{\mathscr h}$} unless otherwise stated. In terms of the connections $\widetilde \nabla$, this means we are interested in differentiating horizontal vector fields w.r.t.\ horizontal directions, with \autoref{cond:nabla} fixing the horizontal part of such covariant derivatives, while \autoref{cond:lin} and the optional \autoref{cond:itoEqStrat} impose limitations on their vertical part. The reader who is versed in sub-Riemannian geometry may spot the link with horizontal connections \cite[Definition 7.4.1]{CC09}, although it should be remarked that our setting is more specific (i.e.\ not all sub-Riemannian manifolds arise as the total space of a vector or even fibre bundle), and the requirements on the connection is somewhat different (on the one hand we are only interested in $\widetilde \nabla_U$ with $U$ horizontal, and on the other also consider the vertical components of such covariant derivatives). In coordinates
\begin{equation}\label{eq:Fptr}
\begin{split}
F^\gamma_{\alpha\beta} \hspace{-0.6em} &\hspace{0.6em}  = -\Gamma^\gamma_{\alpha\beta} \\
\tensor{\widetilde F}{_{\alpha\beta \delta}^\gamma} y^\delta &\stackrel{(\alpha\beta)}{=} -\Gamma^\gamma_{\varepsilon\delta} \Gamma^\varepsilon_{\alpha\beta} y^\delta - (-\Gamma^\gamma_{\beta\delta,\alpha} y^\delta + \Gamma^\varepsilon_{\alpha\delta}\Gamma^\gamma_{\beta\varepsilon} y^\delta \\
&\mathrel{\hphantom{\stackrel{(\alpha\beta)}{=}}} +\widetilde \Gamma^{\widetilde \gamma}_{\widetilde \xi \widetilde \eta} \Gamma^\xi_{\alpha \delta} \Gamma^\eta_{\beta\varepsilon} y^\delta y^\varepsilon + \widetilde \Gamma^{\widetilde \gamma}_{\alpha\beta} - \Gamma^\xi_{\alpha\delta} \widetilde \Gamma^{\widetilde \gamma}_{\widetilde \xi \beta} y^\delta - \Gamma^\eta_{\beta \delta} \widetilde \Gamma^{\widetilde \gamma}_{\alpha \widetilde \eta} y^\delta)
\end{split}
\end{equation}
Note how \autoref{lem:changeF} agrees with \autoref{eq:chrChange}.

We are now in a position to be able to provide the natural generalisation of parallel transport of vectors and Cartan (anti)development to the setting of non-geometric rough paths, with $\tau TM$ endowed with a linear connection. Since the development of a path is not guaranteed to remain in the manifold for all time, it will be helpful to define the following variations of the rough path spaces (note the use of the double closing parenthesis):
\begin{equation}
\begin{split}
\mathscr C^p_\omega([0,T)],M) &\coloneqq \mathscr C^p_\omega([0,T],M) \cup \{ \bfX \in \mathscr C^p_\omega([0,S),M) \text{ for some } S \leq T \\
&\mathrel{\hphantom{\coloneqq}} \text{ and } \nexists \text{ compact } K \subseteq M \text{ s.t.\ }X_{[0,S)} \subseteq K\} \\
\mathscr C^p_\omega([0,\leq T)],M) &\coloneqq \mathscr C^p_\omega([0,T],M) \cup \bigcup_{0 \leq S \leq T} \mathscr C^p_\omega([0,S),M)
\end{split}
\end{equation}
Note that $\mathscr C^p_\omega([0,T)],M) \subseteq \mathscr C^p_\omega([0,\leq T)],M)$, and we also will use these notations when $M$ is a vector space. Moreover, we will add a modifier in the rough path sets to denote those rough paths which are started at a specific point. The following notions are defined whenever \autoref{cond:nabla} and \autoref{cond:lin} are met.
\begin{defn}\label{def:pt}
	Let $\bfX \in \mathscr C^p_\omega([0, T)],M,o)$. We will denote 
	\begin{equation}
	\ptr{\bfX}{}{}{ts} \coloneqq \Phi(\mathscr h, \bfX)_{ts} \colon T_{X_s}M \xrightarrow{\cong} T_{X_t}M
	\end{equation}
	which by \autoref{eq:linearRDE} is well-defined for all $s, t$ at which $X$ is defined, and call it \emph{parallel transport} of vectors along the non-geometric rough path $\bfX$. We will denote $\Ptr_{ts} \coloneqq \Ptr(\bfX)_{ts}$, $\Ptr_t \coloneqq \ptr{\bfX}{}{}{t} \coloneqq \ptr{\bfX}{}{}{t0} \colon T_oM \to T_{X_t}M$ and $\Aptr_t \coloneqq \ptr{\bfX}{}{}{0t} = \Ptr_t^{-1}$ when there is no ambiguity as to the rough path.
\end{defn}

\begin{rem}[There is no alternate notion of \say{backward parallel transport}]
	A rough path $\bfX$ canonically defines a rough path $\overleftarrow{\bfX} = (\overleftarrow{X}, \overleftarrow{\bbX})$ above the inverted path $\overleftarrow{X}_t \coloneqq X_{T-t}$. This is done by imposing the Chen identity to hold for all $0 \leq s,u,t \leq T$ (not just $s \leq u \leq t$), or equivalently by taking \autoref{eq:bbXint} literally, and results in $\overleftarrow{\bbX}_{st} = -\bbX_{T-t,T-s} + X_{T-t,T-s}^{\otimes 2}$ for $0 \leq s \leq t \leq T$. It is shown that if $\bfH \in \mathscr D_X$ then $\overleftarrow{\bfH} \in \mathscr D_{\overleftarrow{X}}$, where $\overleftarrow{\bfH}_t \coloneqq \bfH_{T-t}$, and that
	\begin{equation}
	\int_0^T \overleftarrow{\bfH} \dif \overleftarrow{\bfX} = - \int_0^T \bfH \dif \bfX
	\end{equation}
	at the trace level. It can then be concluded (by a uniqueness argument) that
	\begin{equation}
	\begin{cases}\dif Y = y_0 + \int F(Y) \dif \bfX \\ \dif \overleftarrow{Y} = Y_T + \int F(\overleftarrow{Y}) \dif \overleftarrow{\bfX}\end{cases} \Longrightarrow \overleftarrow{Y}_t = Y_{T-t} 
	\end{equation}
	which implies that, denoting with $\Phi$ the flow map of the RDE defined by $F, \bfX$ and with $\overleftarrow{\Phi}$ the one defined by $F,\overleftarrow{\bfX}$, $\overleftarrow{\Phi} = \Phi^{-1}$. Therefore, once a rough path is fixed, the definition of $\Aptr$ given above and the one obtained by defining the parallel transport RDE w.r.t.\ $\overleftarrow{\bfX}$ coincide.
\end{rem}

\begin{defn}\label{def:aDev}
	Let $\bfX \in \mathscr C^p_\omega([0, T)],M,o)$. Using \autoref{prop:intPhi} we will denote 
	\begin{equation}\label{eq:aDev}
	\adev{\bfX}{}{}{\cdot} \coloneqq \int_0^\cdot \aptr{\bfX}{}{}{s} \dif_\nabla \bfX_s \in \mathscr C^p_\omega([0,\leq T)],T_oM,0_o)
	\end{equation}
	which we call the \emph{antidevelopment} of $\bfX$. If $\boldsymbol{Z} = \adev{\bfX}{}{}{}$ (up to the time at which $\bfX$ is defined) we will denote $\bfX = \dev{\boldsymbol{Z}}{}{}{}$ and call $\bfX$ the \emph{development} of $\boldsymbol Z$.
\end{defn}
In coordinates \autoref{eq:aDev} amounts to
\begin{equation}\label{eq:aDevCoords}
\dif \ade{\bfX}{\gamma^\circ}{} = \Aptr^{\gamma^\circ}_{\gamma}{} \dif \bfX^\gamma_t + \tfrac 12 \Aptr^{\gamma^\circ}_{\gamma} \Gamma^\gamma_{\alpha\beta}\dif[\bfX]^{\alpha\beta}
\end{equation}

For the moment we have only defined development of a rough path which already is the antidevelopment of an $M$-valued one. If we start from an arbitrary $\bfZ \in \mathscr C^p_\omega([0,T],T_oM,0_o)$ with $Z_0 = 0_o$ we would like to invert \autoref{def:aDev} and define its development as the solution to the path-dependent RDE
\begin{equation}\label{eq:pathDepDev}
\dif_\nabla \de{\bfZ}{}{} = \ptr{\de{\bfZ}{}{}}{}{}{} \dif \bfZ, \quad \de{\bfZ}{}{0} = o
\end{equation}
Heuristically, this means that in an infinitesimal time interval $[t_0,t_0+\dif t]$ we are translating the differential $\dif \bfZ_{t_0} \in T_{Z_{t_{\scaleto{0}{3pt}}}}T_oM$ so that it is based at the origin $0_o$, parallel-transporting it along the already-developed portion of the rough path $\bfX_{[0,t_0]} \coloneqq {\de{\bfZ}{}{[0,t_0]}}$ so that it is now based at $X_{t_0}$, and then using it to \say{roll $T_oM$ on $M$ along $\bfZ$ without slipping} for time $\dif t$. The problem, of course, is that we have not defined such (adaptedly) path-dependent RDEs. Moreover, it should be noted that even once this equation is given a meaning, contrary to the case of parallel transport there is no reason why the solution should not explode (see \cite[Corollary 1.36]{Dri18} for general criteria that rule this out for $\bfX$ geometric).

\vspace{-40pt}
\begin{figure}[h]
	\minipage[b]{0.5\textwidth}
	\includegraphics[width=\linewidth]{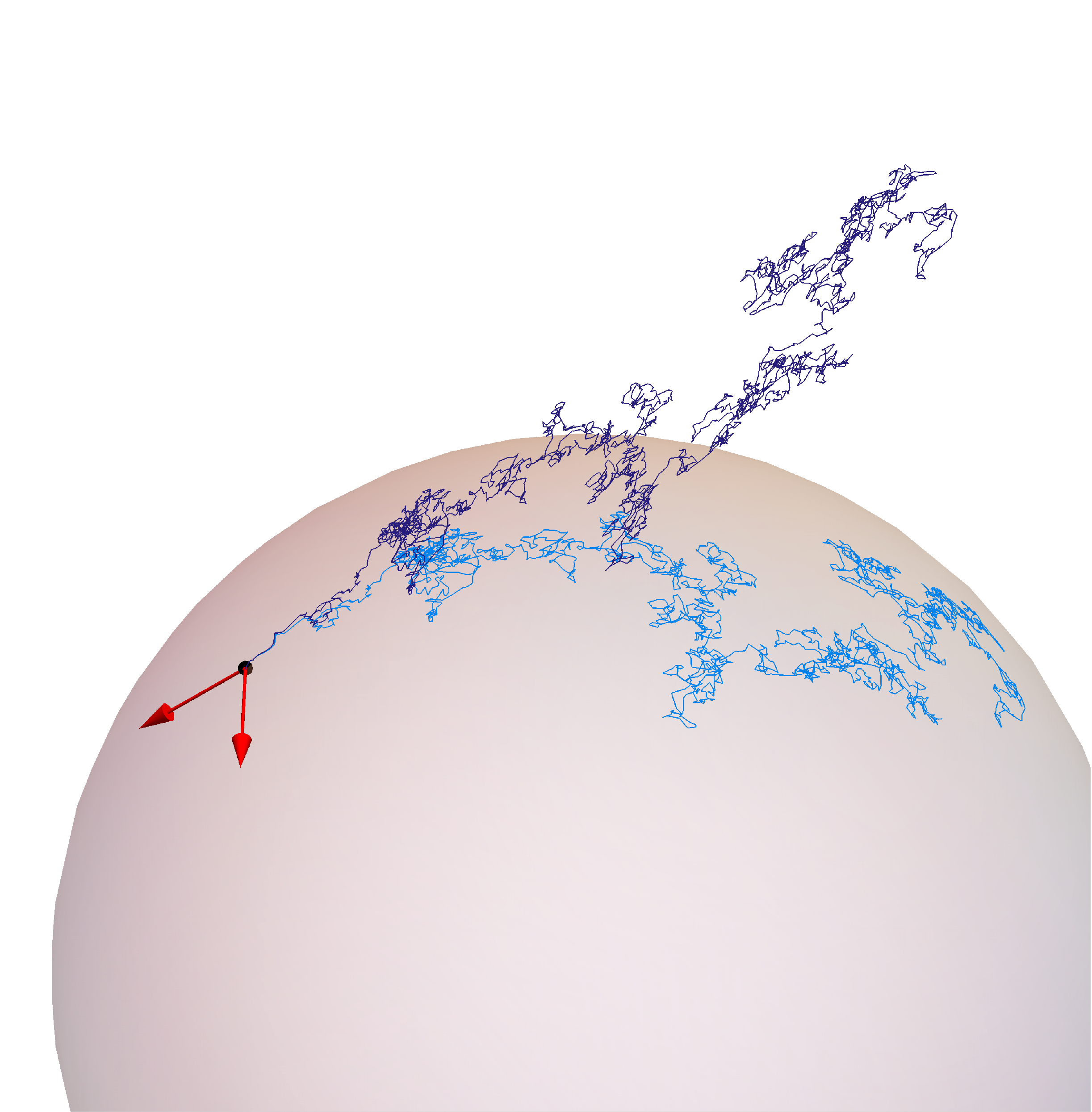}
	\vspace{-20pt}
	\endminipage\hfill
	\minipage[b]{0.5\textwidth}
	\includegraphics[width=\linewidth]{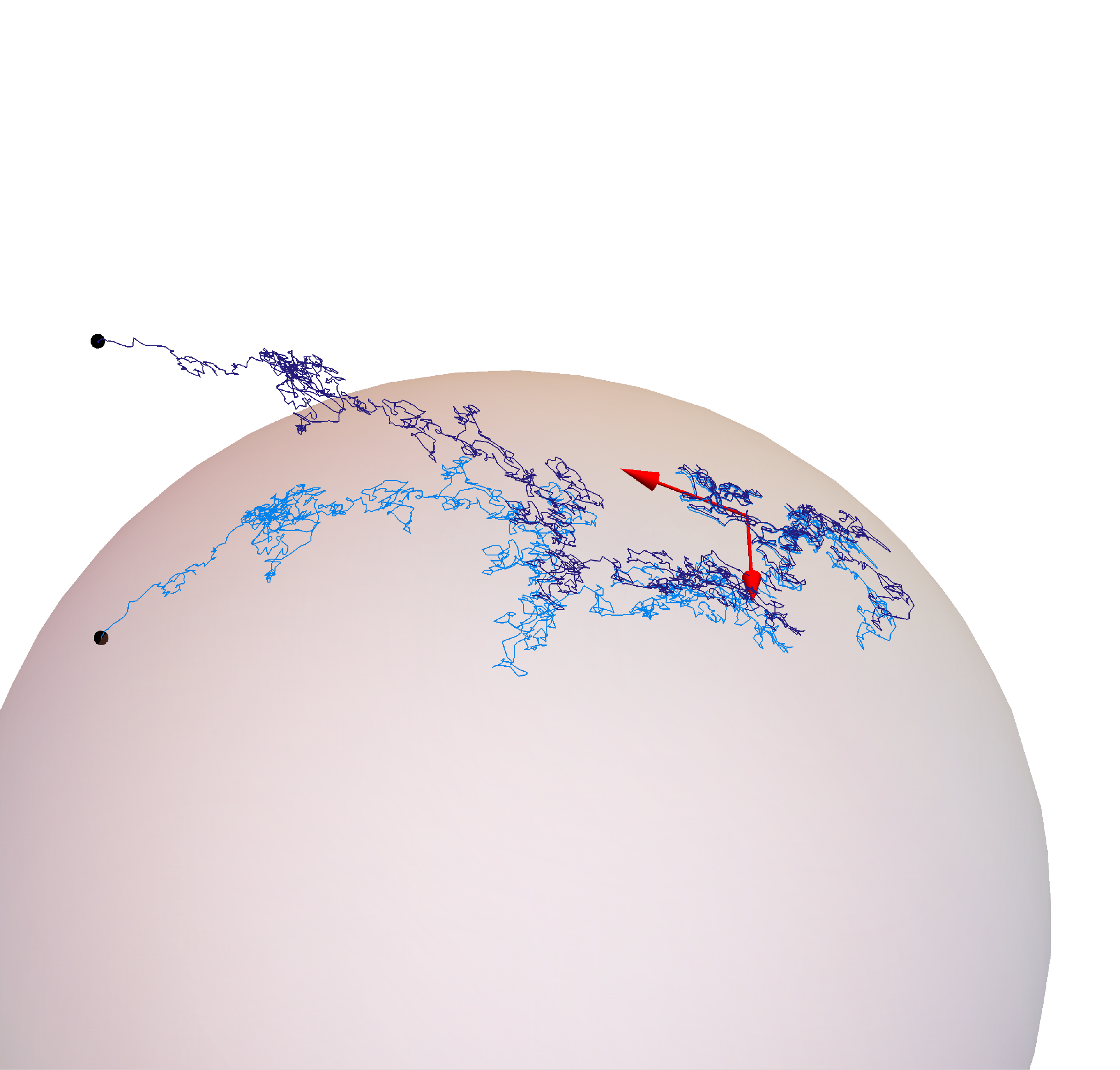}
	\vspace{-20pt}
	\endminipage
	\caption{A 2-dimensional Brownian path in $T_{(1,0,0)}S^2$, plotted in dark blue, and its development onto $S^2$, plotted in light blue. The ``rolling without slipping'' motion is shown at the initial time (on the left) and at a later time (on the right), with the parallel frame based at the point of contact between the tangent space and the manifold.}\label{fig:BM}
\end{figure}
The trick to give \autoref{eq:pathDepDev} a meaning is to consider it jointly with a parallel frame: this transforms the path-dependent RDE into a state-dependent one.
\begin{thm}\label{thm:localDev}
	Let $\bfZ \in \mathscr C^p_\omega([0,T],T_oM)$. Then $\bfX = \dev{\boldsymbol{Z}}{}{}{} \in \mathscr C^p_\omega([0,T],M)$ (possibly up to its exit time from $M$) if and only if $\bfX$ is the unique solution to
	\begin{equation}\label{eq:localDev}
	\begin{split}
	&\begin{dcases}
	\edif \Ptr_{\gamma^\circ}^\gamma &\hspace{-8pt}= \mathscr h^\gamma_{\alpha\beta} \Ptr^\beta_{\gamma^\circ} \Ptr^\alpha_{\delta^\circ} \circ \edif \bfZ^{\delta^\circ} + \tfrac 12 \tensor{\widetilde F(\widetilde \nabla, \mathscr h)}{_{\alpha\beta\delta}^\gamma} \Ptr^\alpha_{\alpha^\circ} \Ptr^\beta_{\beta^\circ} \Ptr^\delta_{\gamma^\circ} \edif[\bfZ]^{\alpha^\circ \beta^\circ} \\
	&\hspace{-8pt}= -\Gamma^\gamma_{\alpha\beta}\Ptr^\beta_{\gamma^\circ} \Ptr^\alpha_{\delta^\circ} \circ \edif \bfZ^{\delta^\circ} \\
	&\hspace{-8pt}\mathrel{\phantom{=}}+ \tfrac 12 \big[-\Gamma^\gamma_{\varepsilon\delta} \Gamma^\varepsilon_{\alpha\beta} \Ptr^\delta_{\gamma^\circ} - (-\Gamma^\gamma_{\beta\delta,\alpha} \Ptr^\delta_{\gamma^\circ} + \Gamma^\varepsilon_{\alpha\delta}\Gamma^\gamma_{\beta\varepsilon} \Ptr^\delta_{\gamma^\circ} \\
	&\hspace{-8pt}\mathrel{\phantom{=}}+\widetilde \Gamma^{\widetilde \gamma}_{\widetilde \xi \widetilde \eta} \Gamma^\xi_{\alpha \delta} \Gamma^\eta_{\beta\varepsilon} \Ptr^\delta_{\gamma^\circ} \Ptr^\varepsilon_{\gamma^\circ} + \widetilde \Gamma^{\widetilde \gamma}_{\alpha\beta} - \Gamma^\xi_{\alpha\delta} \widetilde \Gamma^{\widetilde \gamma}_{\widetilde \xi \beta} \Ptr^\delta_{\gamma^\circ} - \Gamma^\eta_{\beta \delta} \widetilde \Gamma^{\widetilde \gamma}_{\alpha \widetilde \eta} \Ptr^\delta_{\gamma^\circ}) \big] \Ptr^\alpha_{\alpha^\circ} \Ptr^\beta_{\beta^\circ} \edif [\bfZ]^{\alpha^\circ \beta^\circ}
	\\
	\edif \bfX^\gamma &\hspace{-14pt}= \Ptr^\gamma_{\gamma^\circ} \edif \bfZ^{\gamma^\circ} - \tfrac 12 \Gamma^\gamma_{\alpha\beta}\Ptr^\alpha_{\alpha^\circ} \Ptr^\beta_{\beta^\circ} \edif[\bfZ]^{\alpha^\circ \beta^\circ}
	\end{dcases} \\
	&X_0 = o, \quad \{\Ptr_{\gamma^\circ;0}\}_{\gamma^\circ = 1,\ldots,m} \text{ any basis of $T_oM$}
	\end{split} 
	\end{equation}
	with the $\Gamma$'s evaluated at $X_t$ and the $\widetilde \Gamma$'s evaluated at $\Ptr_{\gamma^\circ;t} = \ptr{\bfX}{}{\gamma^\circ}{t}$.
	
	$\Dev$ therefore defines a surjective map $\mathscr C^p_\omega([0, T],T_oM,0_o) \twoheadrightarrow \mathscr C^p_\omega([0,T)],M,o)$ 
	with right inverse $\Adev \colon \mathscr C^p_\omega([0,T)],M,o) \hookrightarrow \mathscr C^p_\omega([0, \leq T)],T_oM,0_o)$ (composed with any map prolonging an element of $\mathscr C^p_\omega([0, S),T_oM,0_o)$ up to time $T$, e.g.\ trivially). In particular, if $M$ is compact $\Dev$ takes values in $\mathscr C^p_\omega([0,T],M,o)$, i.e.\ development exists for all time.
\end{thm}
If $\bfZ$ is geometric this equation may be stated more elegantly as taking values in the frame bundle $\phi M \colon FM \to M$, and defined by the fundamental horizontal vector fields, i.e.\
\begin{equation}
\dif \bfY = \mathscr H_{\lambda^\circ}(Y) \dif \bfZ^{\lambda^\circ}, \ Y_0 \in F_oM \quad \Longrightarrow \quad \dev{\bfZ}{}{}{} = \phi M_* \bfY, \ Y = \ptr{\bfX}{}{}{}
\end{equation}
In this context, compare \autoref{eq:localDev} with \cite[(3.3.9) p.86]{Hsu02}, which is stated in the case of $X$ a Brownian motion, although the formula generalises to more general processes/geometric rough paths. We have decided not to consider frame bundle-valued RDEs in the non-geometric case, since this would require defining a connection on $FM$, which is a delicate matter (some comments to this effect are provided in \cite[p.439]{E90} in the case of the complete lift, though these do not contain an exhaustive description of the connection on $\tau FM$). We have preferred to define development in a coordinate-free manner by simply declaring $\bfX$ to be the development of $\bfZ$ if $\bfZ$ is the antidevelopment of $\bfX$ (as done in \autoref{def:aDev}), and only relying on the local description involving the parallel frame (seen as $m$ vectors which are parallel-transported individually) as an alternative characterisation, useful for explicit computations; in this approach only non-geometric parallel transport of vectors is needed.
\begin{proof}[Proof of \autoref{thm:localDev}]
	By \autoref{eq:aDevCoords} $\bfX = \Dev(\bfZ)$ means
	\begin{equation}
	\begin{split}
	&\dif \bfZ_t^{\gamma^\circ} = \Aptr^{\gamma^\circ}_{\gamma} \dif \bfX^\gamma + \tfrac 12 \Aptr^{\gamma^\circ}_{\gamma} \Gamma^\gamma_{\alpha\beta}\dif[\bfX]^{\alpha\beta}, \quad Z_0 = 0_o \\
	\Rightarrow \quad &\dif[\bfZ]^{\alpha^\circ \beta^\circ} = \Aptr^{\alpha^\circ}_{\alpha}\Aptr^{\beta^\circ}_{\beta} \dif [\bfX]^{\alpha\beta}
	\end{split}
	\end{equation}
	and we have
	\begin{equation}\label{eq:AdZ}
	\begin{split}
	&\mathrel{\phantom{=}}\Ptr^\gamma_{\gamma^\circ} \dif \bfZ^{\gamma^\circ} - \tfrac 12 \Gamma^\gamma_{\alpha\beta}\Ptr^\alpha_{\alpha^\circ} \Ptr^\beta_{\beta^\circ} \dif[\bfZ]^{\alpha^\circ \beta^\circ}_t \\
	&= \Ptr^{\gamma}_{\gamma^\circ} \Aptr^{\gamma^\circ}_{\delta} \dif \bfX^\delta + \tfrac 12 \Ptr^{\gamma}_{\gamma^\circ}\Aptr^{\gamma^\circ}_{\delta} \Gamma^\delta_{\alpha\beta} \dif [\bfX]^{\alpha\beta} - \tfrac 12 \Gamma^\gamma_{\alpha\beta}\Ptr^\alpha_{\alpha^\circ} \Ptr^\beta_{\beta^\circ} \Aptr^{\alpha^\circ}_{\mu}\Aptr^{\beta^\circ}_{\nu} \dif [\bfX]^{\mu\nu}  \\
	&= \dif \bfX^\gamma_t
	\end{split}
	\end{equation}	
	By \autoref{prop:Phi} and \autoref{eq:Fptr} we have
	\begin{equation}\label{eq:dAgamma}
	\begin{split}
	\dif \Ptr^\gamma_{\gamma^\circ} &= F^\gamma_{\varepsilon\delta} \Ptr^\delta_{\gamma^\circ} \circ \dif \bfX^\varepsilon + \tfrac 12 \tensor{\widetilde F}{_{(\alpha\beta)\delta}^\gamma}\Ptr^\delta_{\gamma^\circ} \dif [\bfX]^{\alpha\beta} \\
	&= F^\gamma_{\varepsilon\delta} \Ptr^\delta_{\gamma^\circ} \dif \bfX^\varepsilon + \tfrac 12 \big[ \tensor{\widetilde F}{_{(\alpha\beta)\delta}^\gamma} + \partial_\alpha F^\gamma_{\beta\delta} + F^\varepsilon_{\alpha\delta} F^\gamma_{\beta\varepsilon} \big] \Ptr^\delta_{\gamma^\circ} \dif [\bfX]^{\alpha\beta} \\
	&= -\Gamma^\gamma_{\varepsilon\delta} \Ptr^\delta_{\gamma^\circ} \Ptr^\varepsilon_{\varepsilon^\circ} \dif \bfZ^{\varepsilon^\circ} \\
	&\mathrel{\phantom{=}} + \tfrac 12 \big[ \tensor{\widetilde F}{_{(\alpha\beta)\delta}^\gamma} - \Gamma^\gamma_{\beta\delta,\alpha} + \Gamma^\varepsilon_{\alpha\delta} \Gamma^\gamma_{\beta\varepsilon} + \Gamma^\gamma_{\varepsilon\delta}\Gamma^\varepsilon_{\alpha\beta} \big] \Ptr^\delta_{\gamma^\circ} \Ptr^\alpha_{\alpha^\circ} \Ptr^\beta_{\beta^\circ} \dif [\bfX]^{\alpha\beta} \\
	&= -\Gamma^\gamma_{\varepsilon\delta} \Ptr^\delta_{\gamma^\circ} \Ptr^\varepsilon_{\varepsilon^\circ} \circ \dif \bfZ^{\varepsilon^\circ} + \tfrac 12 \tensor{\widetilde F}{_{(\alpha\beta)\delta}^\gamma} \Ptr^\delta_{\gamma^\circ} \Ptr^\alpha_{\alpha^\circ} \Ptr^\beta_{\beta^\circ} \dif [\bfX]^{\alpha\beta} 
	\end{split}
	\end{equation}
	where in the last step the Gubinelli derivative of $-\Gamma^\gamma_{\beta\delta}(X) \Ptr^\delta_{\gamma^\circ} \Ptr^\varepsilon_{\varepsilon^\circ}$ w.r.t.\ $Z^\alpha$ is computed thanks to the previous step and \autoref{eq:AdZ}. Retracing these steps proves the converse. Note that we do not need to show the coordinate invariance of \autoref{eq:localDev}, as we have shown it is equivalent to $\bfZ = \Adev(\bfX)$, which is defined in \autoref{def:aDev} without reference to a coordinate system.
	
	The map $\Dev$ is then well-defined by uniqueness of RDE solutions applied to the $(m + m^2)$-dimensional system in each coordinate patch, and its right inverse is $\Adev$ by definition. It only remains to show that $\Dev(\bfZ)$ is either defined up to time $T$ or that it is defined up to and excluding some $S \leq T$ with the image of its trace not contained in any compact of $M$. Assume $(X,\Ptr)$ is defined up to time $S$ with $X_{[0,S)}$ contained in a compact $K$ of $M$. Therefore there exists $t_n \searrow S$ s.t.\ $\lim X_{t_n} = \overline x \in K$. We now show that for any neighbourhood $V$ of $\overline x$ there exists $s_0$ s.t.\ $X_{[s_0,S)} \subseteq V$. Consider the image of \autoref{eq:localDev} (defined in $V$) through a change of coordinates $\Phi$ that maps the $\Ptr$ components to a compact, and extend the resulting coefficients smoothly. Now picking a second neighbourhood $U$ of $\overline x$ s.t.\ $\text{Im}(U) \subseteq U$, $\overline U \subseteq V$, an application of \autoref{lem:localUnifE} proves the claim by picking $s_0$ s.t.\ $s_0 < S < s_0+\delta$. (The change of coordinates was necessary because we need to be able to start the equation for $(X,\Ptr)$ at an arbitrary point in $TU^m$.) We may then reason as in \autoref[Proof]{thm:linearFibreRDE} to conclude that $(X,\Ptr)_{[0,S)}$ must also lie in a compact of $TM^m$, and a second application of \autoref{lem:localUnifE} (arguing as in \cite[Theorem 4.2]{CDL15}) then shows that the solution may be prolonged past $S$ (or with its limit if $S = T$). This concludes the proof.\todo{is this level of rigour sufficient? This is quite a delicate argument, given how specifically the theory has been developed.}
\end{proof}
The following result is proven in \cite[Theorem 8.22]{E89} in the case of Stratonovich parallel transport, and interestingly it carries over to the more general case. 
\begin{cor}\label{cor:stratDev}
	At the trace level we have
	\begin{equation}\label{eq:aDevStrat}
	\Adev(\bfX)_\cdot \coloneqq \int_0^\cdot \Aptr(\bfX)_s \circ \edif \bfX_s
	\end{equation}
	and we may replace 
	\begin{equation}
	\edif X^\gamma = \Ptr^\gamma_{\gamma^\circ} \circ \edif \bfZ^{\gamma^\circ}
	\end{equation}
	for the second equation of \autoref{eq:localDev}.
\end{cor}
\begin{rem}
	We emphasise that this does not mean that the (anti)development of a rough path coincides with that of its geometrisation (including at the trace level): in \autoref{eq:aDevStrat} parallel transport is still carried out with reference to the original non-geometric $\bfX$ (and thus depends on the choice of $\widetilde \nabla$), and in the case of development, the first equation of \autoref{eq:localDev} still has the $\dif [\bfZ]$ terms, which are not present when developing ${_\text{g}\!}\bfZ$. Moreover, at the second order level $\bbX_{st}^{\alpha\beta} \approx \Ptr^{\alpha}_{\alpha^\circ;s} \Ptr^{\beta}_{\beta^\circ;s} \mathbb Z^{\alpha^\circ\beta^\circ}_{st}$ locally in terms of the original rough path $\bfZ$.
\end{rem}
\begin{proof}[Proof of \autoref{cor:stratDev}]
	By \autoref{prop:Phi} and \autoref{eq:aDevCoords} we have, at the trace level
	\begin{equation}
	\begin{split}
	\Aptr^{\gamma^\circ}_{\gamma} \circ \dif \bfX &= \Aptr^{\gamma^\circ}_{\gamma} \dif \bfX^\gamma_t + \tfrac 12 \Aptr^{\gamma^\circ}_{\gamma} \Gamma^\gamma_{\alpha\beta} \dif[\bfX]^{\alpha\beta} = \Aptr^{\gamma^\circ}_\gamma \dif_\nabla \bfX^\gamma
	\end{split}
	\end{equation}
	and the second claim is proved analogously by using \autoref{eq:localDev}.
\end{proof}

If $M$ is Riemannian and $\nabla$ is metric we may further ask under what hypotheses the $\Ptr_{ts}$'s are linear isometries $T_{X_s}M \cong T_{X_t}M$. The following condition does not actually require $F$ to be given by horizontal lift, although we will only apply it in that case.
\begin{cond}\label{cond:metric}
	Let $\mathscr g$ be a Riemannian metric on $M$ and $\nabla$ be $\scrg$-metric: $\langle \widetilde F, U \otimes U \otimes V \rangle \in \Gamma \tau M$ is $\scrg$-orthogonal to $V$ for all $U,V \in \Gamma \tau M$. 
\end{cond}\todo{This condition can also be formulated in the general vector bundle case with general linear $F$, but in that case I don't know how to describe the first-order condition. Note that the condition $\nabla \scrg$ appears because $\mathscr h$ is formulated in terms of $\nabla$, and in general it would be formulated in terms of $F$. Also, in that case, we could try and find the condition on $\nabla$ that makes the condition hold for all $F$, as done for the other two properties.}
Note that we are not requiring $\widetilde \nabla$ to be metric w.r.t.\ a Riemannian metric on the manifold $TM$. The statement of this condition in coordinates is given in the following lemma, whose proof is immediate by polarisation.
\begin{lem}
	In coordinates \autoref{cond:metric} corresponds to
	\begin{equation}
	\tensor{\widetilde F}{_{(\alpha\beta)(\gamma \delta)}} = 0
	\end{equation}
\end{lem}
\begin{thm}\label{thm:metric}
	If \autoref{cond:metric} holds, or if $\bfX$ is geometric, $\Ptr(\bfX)_{ts}$ is a linear isometry for all $0 \leq s, t \leq T$.
\end{thm}
The following pattern has emerged: for each property (well-definedness, linearity, and metricity, each required at the level of generality considered) we have a first-order condition (respectively \autoref{eq:Fab}, $F$ linear, and $\nabla$ $\scrg$-metric - as shall be seen in the proof below) and a second order condition (respectively \autoref{cond:nabla}, \autoref{cond:lin}, \autoref{cond:metric}). The first-order conditions are necessary when considering the geometric (or even ODE) case, whereas the second-order conditions become relevant once the driving rough path is no longer geometric. Note how all three conditions are automatically satisfied when \autoref{cond:itoEqStrat} holds.
\begin{proof}[Proof of \autoref{thm:metric}]
	We may assume $s = 0$; then for $y,z \in T_oM$ by \autoref{prop:Phi} we have
	\begin{equation}
	\begin{split}
	&\mathrel{\phantom{=}}\dif \langle \mathscr g(X), \Ptr_{\alpha^\circ} \otimes \Ptr_{\beta^\circ} \rangle \\
	&= \dif (\scrg_{\alpha \beta} \Ptr_{\alpha^\circ}^\alpha \Ptr_{\beta^\circ}^\beta ) \\
	&= \scrg_{\alpha\beta,\gamma}  \Ptr_{\alpha^\circ}^\alpha \Ptr_{\beta^\circ}^\beta \circ \dif \bfX^\gamma  \\
	&\mathrel{\phantom{=}}+ \scrg_{\alpha\beta} \Ptr_{\alpha^\circ}^\delta \Ptr_{\beta^\circ}^\beta (-\Gamma^\alpha_{\gamma \delta} \circ \dif \bfX^\gamma  + \tfrac 12 \tensor{\widetilde F}{_{( \xi\eta)\delta}^\alpha} \dif [\bfX]^{ \xi\eta} ) \\
	&\mathrel{\phantom{=}}+ \scrg_{\alpha\beta} \Ptr_{\alpha^\circ}^\alpha \Ptr_{\beta^\circ}^\delta (-\Gamma^\beta_{\gamma \delta} \circ \dif \bfX^\gamma  + \tfrac 12 \tensor{\widetilde F}{_{( \xi\eta)\delta}^\beta} \dif [\bfX]^{ \xi\eta} ) \\
	&= \Ptr_{\alpha^\circ}^\alpha \Ptr_{\beta^\circ}^\beta \big[(\scrg_{\alpha\beta,\gamma} - \scrg_{\alpha \delta}\Gamma^\delta_{\gamma\beta} - \scrg_{\delta\beta}\Gamma^\delta_{\gamma\alpha}) \circ \dif \bfX^\gamma  + \tensor{\widetilde F}{_{(\xi\eta)(\alpha\beta)}} \dif[\bfX]^{\xi\eta} \big]
	\end{split}
	\end{equation}
	which vanishes by metricity of $\nabla$, \autoref{eq:nablaMetric} and \autoref{cond:metric} or by vanishing of the bracket in the case of $\bfX$ geometric (note how by \autoref{thm:doobMey} the hypotheses are sharp in the case of $X$ truly rough).
\end{proof}

We will now provide three examples of connections $\widetilde \nabla$ on $\tau M$ which it makes sense to consider. The first two, for which we refer to \cite{YI73}, can be viewed as \say{lifts} of the connection $\nabla$ (which is not assumed to be metric or torsion-free), while the third consists of assuming $M$ is Riemannian, defining a Riemannian metric on the manifold $TM$, and taking its Levi-Civita connection.
\begin{expl}[The complete lift of $\nabla$]
	Assume $\nabla$ is a linear connection on $\tau M$, which we do not assume to be torsion-free or metric. The \emph{complete lift} $\widetilde \nabla$ of $\nabla$ is the linear connection on $\tau TM$ whose Christoffel symbols in induced coordinates (w.r.t.\ to any chart $\varphi$ on $M$) are given as functions of the Christoffel symbols $\Gamma^k_{ij}$ of $\nabla$ w.r.t.\ to $\varphi$ as follows:
	\begin{align}\label{eq:comGamma}
	\begin{split}
	&\widetilde \Gamma^\gamma_{\alpha \beta }(x,y) = \Gamma^\gamma_{\alpha \beta }(x), \quad \widetilde \Gamma^\gamma_{\alpha\widetilde \beta }(x,y) = \widetilde \Gamma^\gamma_{\widetilde \alpha  \beta }(x,y) = \widetilde \Gamma^\gamma_{\widetilde \alpha  \widetilde \beta }(x,y) = 0 \\
	&\widetilde \Gamma^{\widetilde \gamma}_{\alpha \beta }(x,y) = \partial_\lambda \Gamma^\gamma_{\alpha \beta }(x)y^\lambda, \quad \widetilde \Gamma^{\widetilde \gamma}_{\alpha  \widetilde \beta }(x,y) = \Gamma^\gamma_{\alpha \beta }(x), \quad \widetilde \Gamma^{\widetilde \gamma}_{\widetilde \alpha  \beta }(x,y) = \Gamma^\gamma_{\alpha \beta }(x), \quad \widetilde \Gamma^{\widetilde \gamma}_{\widetilde \alpha  \widetilde \beta }(x,y) = 0
	\end{split}	
	\end{align}
	From these and \autoref{expl:twidleAffine} it follows that $\tau M$ is an affine map w.r.t.\ $\widetilde \nabla$, $\nabla$. This connection admits the following simple description: having defined the complete lift of $V \in \Gamma \tau M$ as $\widetilde V \in \Gamma \tau TM$ given in induced coordinates by
	\begin{equation}
	\widetilde V^\gamma(x,y) \coloneqq V^\gamma(x), \quad \widetilde V^{\widetilde \gamma}(x,y) \coloneqq y^\lambda \partial_\lambda V^\gamma(x)
	\end{equation}
	(this is checked to be a sound definition; note that no further connection is needed to perform this lift) $\widetilde \nabla$ is characterised by the condition
	\begin{equation}
	\widetilde\nabla_{\widetilde U} \widetilde V = \widetilde{\nabla_UV}, \quad U,V \in \Gamma \tau M
	\end{equation}
	We will only need the local description of $\widetilde\nabla$. However, we remark that the complete lift can be extended to tensor fields, and in particular to Riemannian metrics $\mathscr g$, thus yielding a pseudo-Riemannian metric $\widetilde {\mathscr g}$ on $TM$ (with metric signature $(m,m)$) whose components are given by
	\begin{equation}\label{eq:gc}
	\begin{pmatrix}
	\widetilde {\mathscr g}_{\alpha \beta } & \mathscr g_{\alpha  \widetilde \beta } \\ \mathscr g_{\widetilde \alpha  \beta } & \mathscr g_{\widetilde \alpha  \widetilde \beta } 
	\end{pmatrix}(x,y) = \begin{pmatrix} \partial_\lambda \mathscr g_{\alpha \beta }(x) y^\lambda &  \mathscr g_{\alpha \beta }(x) \\ \mathscr g_{\alpha \beta }(x) & 0
	\end{pmatrix}
	\end{equation}
	If $\nabla$ is $\mathscr g$-metric, then $\widetilde \nabla$ is $\widetilde {\mathscr g}$-metric, and if $\nabla$ is torsion-free then so is $\widetilde \nabla$; therefore $\widetilde{{^\mathscr g \!}\nabla} = {^{\widetilde{\mathscr g}} \!} \nabla$. In general, $\widetilde \nabla$ has the property that its geodesics are given by the Jacobi fields of $\nabla$.
	
	It is easily checked using the theory in this section that \autoref{cond:nabla} and \autoref{cond:lin} are satisfied for all $F$ in the case of the complete lift, and in the case of parallel transport with $\nabla$ torsion-free we have
	\begin{equation}
	\tensor{\widetilde F}{_{\alpha\beta\delta}^\gamma} = \tensor{\mathscr R}{_{\alpha\delta\beta}^{\gamma}}
	\end{equation}
	\autoref{cond:metric}, however, is not satisfied even when $\nabla$ is Levi-Civita, since
	\begin{equation}
	\begin{split}
	\tensor{\widetilde F}{_{(\alpha\beta)(\delta\gamma)}} &= \tfrac 14 (\tensor{\mathscr R}{_{\alpha\delta\beta\gamma}} + \tensor{\mathscr R}{_{\beta\delta\alpha\gamma}} + \tensor{\mathscr R}{_{\alpha\gamma\beta\delta}} + \tensor{\mathscr R}{_{\beta\gamma\alpha\delta}}) \\
	&= \tfrac 12 (\tensor{\mathscr R}{_{\alpha\delta\beta\gamma}} + \tensor{\mathscr R}{_{\alpha\gamma\beta\delta}})
	\end{split}
	\end{equation}
	which does not vanish in general. The resulting parallel transport equation was first studied, for semimartingales, in \cite{DoGu78} and subsequently in \cite[(27)]{Mey82} (we caution the reader that the convention regarding the indices of the curvature tensor differ from the ones used in \autoref{eq:Rcoords}), and it was realised in \cite[p.437]{E90} that this type of parallel transport fits into the framework of SDEs of the type defined in \autoref{def:RDEM}.
\end{expl}

\begin{expl}[The horizontal lift of $\nabla$]
	The second lift of a connection which we examine is the \emph{horizontal lift} of $\nabla$, which we also denote $\widetilde \nabla$ (equivocation will easily be avoided, since we will always use each connection separately). Its Christoffel symbols in induced coordinates are similar to those of the complete lift, with one important difference:
	\begin{align}\label{eq:horGamma}
	\begin{split}
	&\widetilde \Gamma^\gamma _{\alpha \beta }(x,y) = \Gamma^\gamma _{\alpha \beta }(x), \quad \widetilde \Gamma^\gamma _{\alpha \widetilde \beta }(x,y) =  \widetilde \Gamma^\gamma _{\widetilde \alpha  \beta }(x,y) =  \widetilde \Gamma^\gamma _{\widetilde \alpha  \widetilde \beta }(x,y) = \widetilde \Gamma^{\widetilde \gamma }_{\widetilde \alpha  \widetilde \beta }(x,y) = 0 \\
	&\widetilde \Gamma^{\widetilde \gamma }_{\alpha \beta }(x,y) = (\partial_\lambda \Gamma^\gamma _{\alpha \beta } - \tensor{\mathscr R}{_{\lambda \alpha \beta }^\gamma })(x)y^\lambda, \ \widetilde \Gamma^{\widetilde \gamma }_{\alpha  \widetilde \beta }(x,y) = \Gamma^\gamma _{\alpha \beta }(x), \ \widetilde \Gamma^{\widetilde \gamma }_{\widetilde \alpha  \beta }(x,y) = \Gamma^\gamma _{\alpha \beta }(x)
	\end{split}	
	\end{align}
	As for the complete lift, $\tau M$ is an affine map w.r.t.\ $\widetilde \nabla$, $\nabla$; the extra term appearing in $\widetilde \Gamma^{\widetilde \gamma}_{\alpha\beta}$, however, causes $\widetilde \nabla$ to have nonvanishing torsion in general even if $\nabla$ is torsion-free. Just as for the complete lift, the horizontal lift of a connection is motivated by a broader construction which involves lifting other objects defined on $M$, such as vector fields. However, unlike the case of the complete lift, these lifts require a connection on $\tau M$ to begin with, and are performed in a way which is related to \autoref{eq:horlift}; we do not provide more details here. If $\nabla$ is $\mathscr g$-metric, then $\widetilde \nabla$ is $\widetilde{\mathscr g}$-metric, where $\widetilde{\mathscr g}$ is the pseudo-Riemannian metric \autoref{eq:gc} (although, unlike the case of the complete lift, $\widetilde{{^\mathscr g \!}\nabla} \neq {^{\widetilde{\mathscr g}} \!} \nabla$ because the former has torsion in general). The characterisation of geodesics of the horizontal lift of a connection is more complicated than that of its complete lift, but it still holds that $\tau M$ maps $\widetilde \nabla$-geodesics to $\nabla$-geodesics. Moreover, it holds that horizontal lifts of geodesics (namely curves in $TM$ above geodesics whose tangent vectors are horizontal, i.e.\ parallel transports above geodesics) define geodesics w.r.t.\ the horizontal lift: this is seen from \cite[Equation 9.4, p.115]{YI73}.
	
	Like the complete lift, the horizontal lift results in \autoref{cond:nabla} and \autoref{cond:lin} being satisfied for all $F$, but in the case of $F$ given by horizontal lift it additionally satisfies \autoref{cond:itoEqStrat}. Therefore the resulting parallel transport is, at the trace level, the same as geometric/Stratonovich parallel transport, a conclusion which is also noted in \cite{Mey82}, \cite{E90}.
\end{expl}

\begin{expl}[The Sasaki metric]
	Let $\mathscr g$ be a Riemannian metric on $M$. We can lift $\mathscr g$ to a Riemannian metric $\widetilde{\mathscr g}$ on $TM$, called the \emph{Sasaki metric} in the following way: recalling the notations introduced in \autoref{subsec:linearconn} for vertical and horizontal bundles, we declare for all $U(x) \in T_xM$
	\begin{equation}
	V_{U(x)} \tau M \ \bot \ H_{U(x)}, \quad \widetilde{\mathscr g}|_{V_{U(x)} \tau M} \coloneqq \mathscr (\mathscr v(U(x))^{-1})^*\mathscr g, \quad \widetilde{\mathscr g}|_{H_{U(x)}} \coloneqq \mathscr (\mathscr h(U(x))^{-1})^*\mathscr g
	\end{equation}
	In induced coordinates, this is amounts to
	\begin{equation}
	\begin{pmatrix}
	\widetilde{\mathscr g}_{\alpha \beta } & \widetilde\scrg_{\alpha \widetilde \beta } \\
	\widetilde \scrg_{\widetilde \alpha  \beta } &  \widetilde \scrg_{\widetilde \alpha  \widetilde \beta }
	\end{pmatrix}(x,y) = \begin{pmatrix}
	\scrg_{\alpha \beta }(x) + \scrg_{\delta  \varepsilon}\Gamma^\delta_{\mu \alpha } \Gamma^\varepsilon_{\nu \beta }(x)y^\mu y^\nu & \Gamma^\gamma _{\alpha \lambda}\scrg_{\gamma \beta }(x) y^\lambda \\
	\Gamma^\gamma _{\lambda \beta }\scrg_{\alpha \gamma }(x) y^\lambda & \scrg_{\alpha \beta }(x)
	\end{pmatrix}
	\end{equation}
	and
	\begin{equation}
	\begin{pmatrix}
	\widetilde{\mathscr g}^{\alpha \beta } & \widetilde\scrg^{\alpha \widetilde \beta } \\
	\widetilde \scrg^{\widetilde \alpha  \beta } &  \widetilde \scrg^{\widetilde \alpha  \widetilde \beta }
	\end{pmatrix}(x,y) = \begin{pmatrix}
	\scrg^{\alpha \beta }(x) & -\Gamma^\beta _{\lambda \gamma }\scrg^{\alpha \gamma }(x) y^\lambda \\
	-\Gamma^\alpha _{\gamma  \lambda}\scrg^{\gamma \beta }(x) y^\lambda & \scrg^{\alpha \beta }(x) + g^{\delta \varepsilon} \Gamma^\alpha _{\delta\mu}\Gamma^\beta_{\varepsilon\nu}(x) y^\mu y^\nu
	\end{pmatrix}
	\end{equation}
	where the $\Gamma$'s are the Christoffel symbols of ${^\scrg\!}\nabla$. The horizontal lift of ${^\scrg\!}\nabla$ is $\widetilde \scrg$-metric, but does not coincide with ${^{\widetilde \scrg}\!}\nabla$ due to torsion. We will call ${^{\widetilde \scrg}\!}\nabla$ the \emph{Sasaki lift} of ${^\mathscr{g}\!}\nabla$ (even though, strictly speaking, it is the metric that we are lifting). The Christoffel symbols of ${^{\widetilde \scrg}\!}\nabla$ in induced coordinates have more complex expressions than the ones for the other two lifts of connections, and are given as functions of the Christoffel symbols of ${^\mathscr{g}\!}\nabla$ and of the components of its curvature tensor by 
	\begin{align}\label{eq:sasGamma}
	\begin{split}
	&\widetilde \Gamma^\gamma_{\alpha \beta }(x,y) = \Gamma^\gamma_{\alpha \beta }(x) + \tfrac 12 (\tensor{\mathscr R}{_{\mu \delta \alpha }^\gamma} \Gamma^\delta_{\lambda \beta } + \tensor{\mathscr R}{_{\mu \delta \beta }^\gamma} \Gamma^\delta_{\alpha  \lambda} )(x)y^\lambda y^\mu \\
	&\widetilde \Gamma^\gamma_{\widetilde \alpha  \beta }(x,y) = \tfrac 12 \tensor{\mathscr R}{_{\lambda \alpha \beta }^\gamma}(x) y^\lambda, \quad \widetilde \Gamma^\gamma_{\alpha  \widetilde \beta }(x,y) = \tfrac 12 \tensor{\mathscr R}{_{\lambda \beta \alpha }^\gamma}(x) y^\lambda, \quad \widetilde \Gamma^\gamma_{\widetilde \alpha  \widetilde \beta }(x,y) = \Gamma^{\widetilde \gamma}_{\widetilde \alpha  \widetilde \beta }(x,y) = 0 \\
	&\widetilde \Gamma^{\widetilde \gamma}_{\alpha  \beta }(x,y) =  \tfrac 12 (\tensor{\mathscr R}{_{\alpha \lambda \beta }^\gamma} + \tensor{\mathscr R}{_{\beta \lambda \alpha }^\gamma} + 2 \partial_\lambda \Gamma^\gamma_{\alpha \beta })(x) y^\lambda \\
	&\mathrel{\phantom{\widetilde \Gamma^{\widetilde \gamma}_{\alpha  \beta }(x,y) =}}+ \tfrac 12 \Gamma^\gamma_{\nu \delta}(\tensor{\mathscr R}{_{\varepsilon\mu \alpha }^\delta} \Gamma^\varepsilon_{\lambda \beta } + \tensor{\mathscr R}{_{\varepsilon\mu \beta }^\delta} \Gamma^\varepsilon_{\alpha \lambda})(x)y^\lambda y^\mu y^\nu \\
	&\widetilde \Gamma^{\widetilde \gamma}_{\widetilde \alpha  \beta }(x,y) = \Gamma^\gamma_{\alpha \beta }(x) - \tfrac 12 \Gamma^\gamma_{\mu \delta} \tensor{\mathscr R}{_{\lambda \alpha \beta }^\delta}(x) y^\lambda y^\mu, \\
	&\widetilde \Gamma^{\widetilde \gamma}_{\alpha  \widetilde \beta }(x,y) = \Gamma^\gamma_{\alpha \beta }(x) - \tfrac 12 \Gamma^\gamma_{\mu \delta} \tensor{\mathscr R}{_{\lambda \beta \alpha }^\delta}(x) y^\lambda y^\mu
	\end{split}	
	\end{align}\todo{Note these symbols do not depend directly on the metric $\scrg$, rather only on its Levi-Civita connection: this raises the question of whether this connection can be defined if we just started from an arbitrary connection on $M$ (and what if it is metric but has torsion, or pseudometric...). This could be answered by performing the very long check of transforming the above Christoffel symbols to another induced coordinate system, and seeing whether there is any condition (e.g.\ metricity, torsion-freeness of $\nabla$) on the above expressions of being independent of the induced coordinate system. Strange that this has not been raised anywhere in the literature.}
	These symbols are taken from \cite{Sa58} with one important caveat: the $\tensor{\mathscr R}{_{\alpha\beta\gamma}^\delta}$'s therein have been transcribed into $\tensor{\mathscr R}{_{\gamma \beta \alpha}^\delta}$'s here. This is because the author follows a different ordering in the coordinate expression of the curvature tensor. This convention is not stated in the paper, but it can be deduced by computing any one of the Christoffel symbols involving a curvature term. This check can be performed by using the fact that the horizontal lift of $\nabla$ is $\widehat \scrg$-metric and \autoref{eq:contorsion}. Let $\widetilde{\mathscr T}$ denote the torsion tensor of the horizontal lift of $\nabla$: its only nonzero component is given by
	\begin{equation}
	\tensor{\widetilde{\mathscr T}}{^{\widetilde \gamma }_{\alpha \beta }}(x,y) = (\tensor{\mathscr R}{_{\lambda \beta \alpha }^\gamma } - \tensor{\mathscr R}{_{\lambda \alpha \beta }^\gamma })(x)y^\lambda
	\end{equation}
	Thus $\tensor{\widetilde{\mathscr T}}{_{\alpha}^{\gamma}_{\widetilde \beta}}(x,y) = 0$ and, performing index gymnastics w.r.t.\ $\widetilde \scrg$ and using \autoref{eq:Rij}, \autoref{eq:Rcyclic}, \autoref{eq:Rkh} and \autoref{eq:ijhk} we obtain
	\begin{equation}
	\begin{split}
	\tensor{\widetilde{\mathscr T}}{_{\widetilde \alpha}^{\gamma }_\beta }(x,y) &= \widehat \scrg_{\widetilde \alpha \widetilde \varepsilon } \widehat \scrg^{\gamma \delta }(x,y) (\tensor{\mathscr R}{_{\lambda \beta \delta }^\varepsilon } - \tensor{\mathscr R}{_{\lambda \delta \beta }^\varepsilon })(x)y^\lambda \\
	&= \scrg_{\alpha \varepsilon }  \scrg^{\gamma \delta } (\tensor{\mathscr R}{_{\lambda \beta \delta }^\varepsilon } - \tensor{\mathscr R}{_{\lambda \delta \beta }^\varepsilon })(x)y^\lambda \\
	&= \scrg_{\alpha \varepsilon }  \scrg^{\gamma \delta } (\tensor{\mathscr R}{_{\lambda \beta \delta }^\varepsilon } + \tensor{\mathscr R}{_{\delta \lambda \beta }^\varepsilon })(x)y^\lambda \\
	&= -\scrg_{\alpha \varepsilon }  \scrg^{\gamma \delta } \tensor{\mathscr R}{_{\beta \delta \lambda}^\varepsilon }(x)y^\lambda \\
	&= -\scrg^{\gamma \delta } \tensor{\mathscr R}{_{\beta \delta \lambda \alpha}}(x)y^\lambda \\
	&= -\scrg^{\gamma \delta } \tensor{\mathscr R}{_{\lambda \alpha \beta \delta }}(x)y^\lambda \\
	&= -\tensor{\mathscr R}{_{\lambda \alpha \beta }^\gamma }(x)y^\lambda \\
	\end{split}
	\end{equation}
	Then, since $\tensor{\widetilde{\mathscr K}}{^\gamma_{\widetilde \alpha\beta}} = \frac 12 \tensor{\widetilde{\mathscr T}}{_{\widetilde \alpha}^{\gamma}_\beta}$ (where $\mathscr K$ is the contorsion tensor defined in \autoref{eq:contorsion}), \autoref{eq:horGamma} yield the value of $\widetilde \Gamma^\gamma_{\widetilde \alpha \beta}$ in \autoref{eq:sasGamma}. Similarly to the case of the horizontal lift of a connection, the horizontal lift of a Riemannian geodesic is a geodesic w.r.t.\ Sasaki metric.
	
	\say{Sasaki parallel transport} has not, to our knowledge, been considered in the literature. Like for the horizontal lift, \autoref{cond:itoEqStrat} is satisfied w.r.t.\ to the Sasaki lift of $\nabla$ when $F$ is given by horizontal lift, and the resulting definition of parallel transport is therefore equivalent to the geometrised one. However, unlike the complete and horizontal lifts, \autoref{cond:nabla} cannot be expected to hold for general $F$, since the Sasaki lift does not make $\tau M$ symmetrically affine: this can be seen, again with reference to \autoref{expl:twidleAffine}, by noting that, for instance $\widetilde \Gamma^\gamma_{\widetilde \alpha \beta}(x,y) = \tfrac 12 \tensor{\mathscr R}{_{\lambda \alpha\beta}^\gamma}(x) y^\lambda$ is not, in general, antisymmetric in $\alpha,\beta$. This means we may not, in general, define arbitrary equations \autoref{eq:fibreRDE} when $E = TM$ is given the Sasaki lift of $\nabla$.
\end{expl}

\begin{expl}[Local martingales and Brownian motion]\label{expl:lmBm}
	It is well known that Stratonovich (anti)development preserves local martingales and if $M$ is Riemannian it preserves Brownian motions. In our setting the first statement always holds in all cases (assuming \autoref{cond:nabla} and \autoref{cond:lin} hold, so that (anti)development is defined), as can be easily seen from the local characterisation of manifold-valued martingales \autoref{eq:locMartCoords}, and the local expressions \autoref{eq:aDevCoords}, \autoref{eq:localDev}. This is also observed (by a different argument) in \cite[p.440]{E90}.
	
	As for the preservation of Brownian motion, we first recall that the Levy criterion on manifolds \cite[Proposition 5.18]{E89} immediately implies the following local characterisation of Brownian motion of a Riemannian manifold $(M,\scrg)$: $X$ is a Brownian motion on $M$ if and only if it is a local martingale and
	\begin{equation}
	\dif[X]^{\alpha\beta} = \scrg^{\alpha\beta}(X) \dif t
	\end{equation}
	If $Z$ is a $T_oM$-valued Brownian motion, then if \autoref{cond:metric} holds, we have for $X = \dif \Dev(Z)$
	\begin{equation}
	\dif [X]^{\alpha\beta} = \Ptr^\alpha_{\alpha^\circ} \Ptr^\beta_{\beta^\circ} \dif[Z]^{\alpha^\circ \beta^\circ} = \Ptr^\alpha_{\alpha^\circ} \Ptr^\beta_{\beta^\circ} \delta^{\alpha^\circ \beta^\circ} \dif t = \mathscr g^{\alpha\beta}(X) \dif t
	\end{equation}
	where the last identity holds thanks to the fact that $\scrg^{\alpha^\circ \beta^\circ}(o) = \delta^{\alpha^\circ \beta^\circ}$ and \autoref{thm:metric}. That antidevelopment maps Brownian motions to Brownian motions under the same hypotheses is checked analogously.
	
	We may therefore conclude that (anti)development defined w.r.t.\ the complete, horizontal and Sasaki lifts to preserve local martingales, but only that defined w.r.t.\ the horizontal and Sasaki lifts to preserve Brownian motion.
	
	We also note that we should expect (anti)development taken w.r.t.\ two different $\widetilde \nabla$'s to be different pathwise, even when both satisfy the linearity and metricity conditions. For Brownian motion this might mean that the law of the (anti)developments coincide (i.e.\ they are both Brownian motions), despite the paths defined by the same state $\omega \in \Omega$ being different. Another way of generating pathwise-distinct Brownian motions through (anti)development of the same Brownian motion is by adding a contorsion term (see \autoref{rem:contorsion}) to the Levi-Civita connection $\nabla$ and taking Stratonovich development. In general, by the It\^o isometry the cross-covariance matrix of the It\^o antidevelopments ${^1\!}\Adev(X)$ and ${^2\!}\Adev(X)$ of the same $M$-valued semimartingale $X$ taken w.r.t.\ ${^1\!}\nabla$, ${^1\!}\widetilde\nabla$ on the one hand and ${^2\!}\nabla$, ${^2\!}\widetilde\nabla$ on the other is given by $E[{^1\!}\adev{X}{\alpha^\circ}{}{}{^2\!}\adev{X}{\beta^\circ}{}{}] = E[\int {^1\!}\Aptr^{\alpha^\circ}_\alpha {^2\!}\Aptr^{\beta^\circ}_\beta \dif [X]^{\alpha\beta}]$, with ${^k\!}\Aptr$, $k = 1,2$ denoting the respective parallel transports above $X$.
\end{expl}

\begin{expl}[$\Ptr$ along Brownian motion on Einstein manifolds w.r.t.\ the complete lift]
	We assume $(M,\scrg)$ is an Einstein manifold, i.e.\ a Riemannian manifold whose Ricci tensor is proportional to the metric tensor, $\mathscr R_{\alpha\beta} = \lambda \scrg_{\alpha\beta}$ with $\lambda \in \bbR$ (the best known such example is the sphere, in all dimensions). Let $Z$ be a Brownian motion on $T_oM$ and $X$ its Stratonovich development, an $M$-valued Brownian motion, and we compare the behaviour of Stratonovich parallel transport $\Ptr(X)$ with parallel transport defined w.r.t.\ to the complete lift $\widetilde \nabla$ of the Levi-Civita connection $\nabla$, which we denote $\widetilde{\Ptr}(X)$. By proceeding as in \autoref[Proof]{thm:metric} and \autoref{expl:lmBm} we compute
	\begin{equation}
	\begin{split}
	\dif \scrg(\Ptr_{\alpha^\circ},\widetilde\Ptr_{\beta^\circ}) &= \tfrac 12 \scrg_{\alpha\beta} \Ptr^\alpha_{\alpha^\circ} \tensor{\mathscr R}{_{\xi \gamma \eta}^\beta} \widetilde \Ptr^\gamma_{\beta^\circ} \scrg^{\xi \eta} \dif t \\
	&= -\tfrac 12 \tensor{\mathscr R}{_{\alpha \beta}} \Ptr^\alpha_{\alpha^\circ} \widetilde\Ptr^\beta_{\beta^\circ} \dif t \\
	&= -\tfrac{\lambda}{2} \scrg(\Ptr_{\alpha^\circ},\widetilde\Ptr_{\beta^\circ}) \dif t
	\end{split}
	\end{equation}
	which implies $\scrg(\Ptr_{\alpha^\circ},\widetilde\Ptr_{\beta^\circ}) = \exp(-\lambda t/2)\delta^{\alpha^\circ \beta^\circ}$, and similarly $\scrg(\widetilde \Ptr_{\alpha^\circ},\widetilde\Ptr_{\beta^\circ}) = \exp(-\lambda t)\delta^{\alpha^\circ \beta^\circ}$. In other words $\widetilde \Ptr(X)$ preserves orthogonality, but not orthonormality, since it consists of a scaling by the above exponential factor. Note that this behaviour of $\widetilde \Ptr$ can only be expected to hold along the Brownian motion $X$, and not along $\widetilde X \coloneqq \widetilde{\Dev}(Z)$, the development of $Z$ taken according to the complete lift of $\widetilde \nabla$, which is not in general a Brownian motion (even given the Einstein assumption): this can be seen by writing $\dif \widetilde X = \sum_{\gamma^\circ} \widetilde \Ptr^\alpha_{\gamma^\circ}(\widetilde X) \widetilde\Ptr^\beta_{\gamma^\circ}(\widetilde X) \dif t$ and by showing that the SDE satisfied by $\sum_{\gamma^\circ} \widetilde \Ptr^\alpha_{\gamma^\circ}(\widetilde X) \widetilde\Ptr^\beta_{\gamma^\circ}(\widetilde X)$ has an extra drift term when compared to that satisfied by $\scrg^{\alpha\beta}(\widetilde X)$.
\end{expl}

\begin{expl}[Linearising rough integrals and rewriting Driver's integration by parts formula]\label{expl:Driver}
	Antidevelopment can be used to write rough integrals against $M$-valued rough paths as ones against $T_oM$-valued ones. Let $\bfX \in \mathscr C^p_\omega([0,T],M,o)$ and $\bfH \in \mathscr D_X(\mathcal L(\tau M, \bbR^e))$. Then it is checked that
	\begin{equation}
	\begin{split}
	\Ptr^* \bfH &= ((\Ptr^* H)^c_{\gamma^\circ},(\Ptr^* H)'^c_{\alpha^\circ \beta^\circ}) \coloneqq ( H^c_\gamma \Ptr^\gamma_{\gamma^\circ}, H'^c_{\alpha\beta} \Ptr^{\alpha}_{\alpha^\circ} \Ptr^{\beta}_{\beta^\circ}  - H^c_\gamma \Gamma^\gamma_{\alpha\beta} \Ptr^{\alpha}_{\alpha^\circ} \Ptr^{\beta}_{\beta^\circ} ) \\
	&\in \mathscr D_{Z}(\mathcal L(T_oM, \bbR^e))
	\end{split}
	\end{equation}
	with $\bfZ \coloneqq \Adev(\bfX)$ (independently of the chart used for the coordinates $\alpha,\beta,\gamma$) and that
	\begin{equation}
	\int \bfH \dif_\nabla \bfX = \int \Ptr^* \bfH \dif \Adev(\bfX)
	\end{equation}
	Note how, in particular, this is independent of the connection $\widetilde \nabla$ on $\tau TM$ used to define $\Ptr$ and $\Adev$. Now, assume that $M$ is Riemannian, $\nabla$ is metric and \autoref{cond:metric} holds. Then for $\boldsymbol P \in \mathscr D_X(\tau M^{\oplus e})$
	\begin{equation}
	\begin{split}
	\int \scrg(\boldsymbol P,\dif_\nabla \bfX) &= \int \boldsymbol P^\flat \dif_\nabla \bfX \\
	&= \int \Ptr^* \boldsymbol P^\flat \dif \Adev(\bfX) \\
	&= \int (\Aptr \boldsymbol P)^\flat \dif \Adev(\bfX) \\
	&= \int \Aptr \boldsymbol P \cdot \dif \Adev(\bfX)
	\end{split}
	\end{equation}
	where the dot product denotes the metric at $o$ and
	\begin{equation}
	\Aptr \boldsymbol P = ((\Aptr P)^{\gamma^\circ}, (\Aptr P)'^{\beta^\circ}_{\alpha^\circ}) \coloneqq (\Aptr_\gamma^{\gamma^\circ} P^\gamma, \textstyle\sum_{\beta,\eta} (\Gamma^\beta_{\xi \eta} \Aptr^{\beta^\circ}_{\eta} \Ptr^\xi_{\alpha^\circ}  P^\beta + \Aptr^{\beta^\circ}_\beta P^\beta_\alpha \Ptr^\alpha_{\alpha^\circ}))
	\end{equation}
	The converses of these statements, i.e.\ 
	\begin{equation}\label{eq:flatToCurved}
	\begin{split}
	&\int \bfK \dif \bfZ = \int \Aptr^*\bfK\dif \Dev(\bfX), \quad \int \boldsymbol Q \cdot \dif \bfZ= \int \scrg(\Ptr \boldsymbol Q, \dif \Dev(\bfZ))\\
	&\text{with } \Aptr^* \boldsymbol K = ((\Aptr^* K)_\gamma,(\Aptr^* K)_{\alpha\beta}) \coloneqq (K_{\gamma^\circ} \Aptr^{\gamma^\circ}_\gamma, K_{\alpha^\circ \beta^\circ}\Aptr^{\alpha^\circ}_\alpha \Aptr^{\beta^\circ}_\beta + K_{\beta^\circ} \textstyle\sum_\eta\Gamma^\beta_{\alpha \eta} \Ptr^{\beta^\circ}_\eta  ) \\
	&\text{and } \Ptr \boldsymbol Q = ((\Ptr Q)^\gamma,(\Ptr Q)^\beta_\alpha) \coloneqq (\Ptr^\gamma_{\gamma^\circ} Q^{\gamma^\circ}, \Gamma^\beta_{\alpha\eta} \Ptr^\eta_{\beta^\circ} Q^{\beta^\circ} + \Ptr^\beta_{\beta^\circ} Q^{\beta^\circ}_{\alpha^\circ} \Aptr^{\alpha^\circ}_{\alpha}) \\
	&\text{for } \bfK \in \mathscr D_Z(\mathcal L(T_oM, \bbR^e)), \ \boldsymbol Q \in \mathscr D_Z(T_oM^{\oplus e}) 
	\end{split}
	\end{equation}
	are similarly shown to hold.
	
	As an application of the latter, we show how the integration by parts formula \cite[Theorem 7.32]{Dr04} can be rewritten as an It\^o integral on $M$. Let $Z$ be a $T_oM$-valued Brownian motion, $X$ its Stratonovich development, $H$ a Cameron-Martin process above $X$, $h \coloneqq \Aptr H$ with $h = \int u \dif t$, $U \coloneqq \Ptr u$ (for the precise terminology pertaining to curved Wiener space see the above reference). Then we may write the integration by parts formula, i.e.\ a formula for the adjoint of the gradient operator
	\begin{equation}\label{eq:Driver}
	\begin{split}
	\mathscr D^* H &= \int \Big(u^\cdot + \frac 12 \Aptr_\gamma^\cdot \tensor{\mathscr R}{^\gamma_\beta}(X) \Ptr^\beta_{\beta^\circ} h^{\beta^\circ} \Big) \cdot \dif Z \\
	&= \int \scrg\Big( U + \frac 12 \tensor{\mathscr R}{^\cdot_\gamma}(X) H^\gamma, \dif_\nabla X \Big)
	\end{split}
	\end{equation}
	as an It\^o integral on $M$. Moreover, if $u$ admits a Gubinelli derivative w.r.t.\ $Z$, so does $U$ w.r.t.\ $X$, and \autoref{expl:itoStratH}, \autoref{expl:rri} may be combined to yield the expression of this as a Stratonovich integral on $M$, plus a correction term involving the covariant derivative of Ricci tensor.
\end{expl}

\todo[inline,backgroundcolor=yellow]{I would have liked to provide the link with development on a Lie group (and with the hyperbolic development used in the signature papers) but was not able to. I think the two notions are fairly disconnected, since I think the Lie group development uses Lie parallel transport (which does not require a connection). It could still hold that e.g.\ if the connection is left-invariant then the two notions agree. It would be good to make this remark, and particularly to provide the link with hyperbolic development, which figures on many papers on the signature.}

\begin{figure}[h]
	\minipage[b]{0.5\textwidth}
	\includegraphics[width=\linewidth]{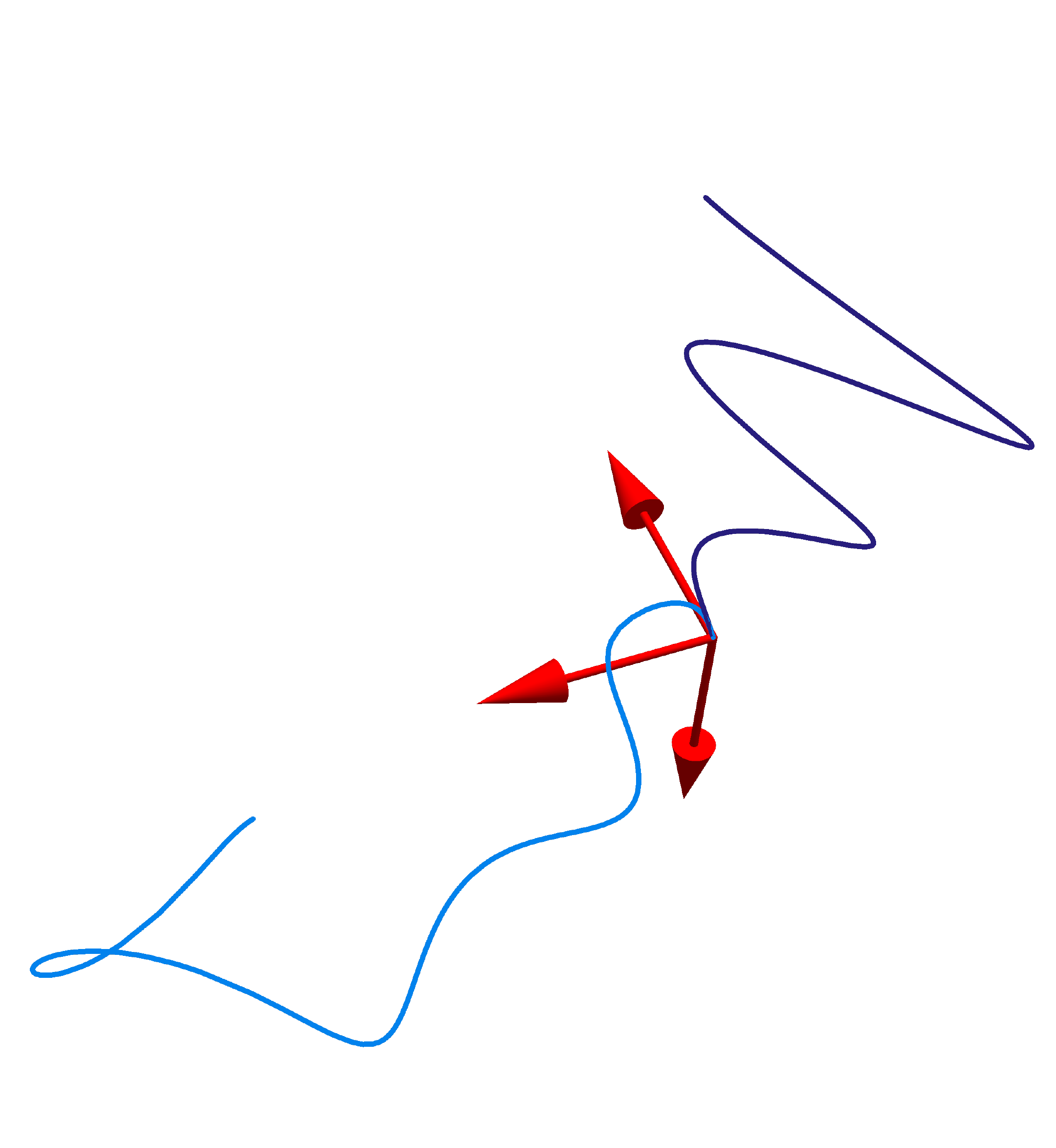}
	\vspace{-30pt}
	\endminipage\hfill
	\minipage[b]{0.5\textwidth}
	\includegraphics[width=\linewidth]{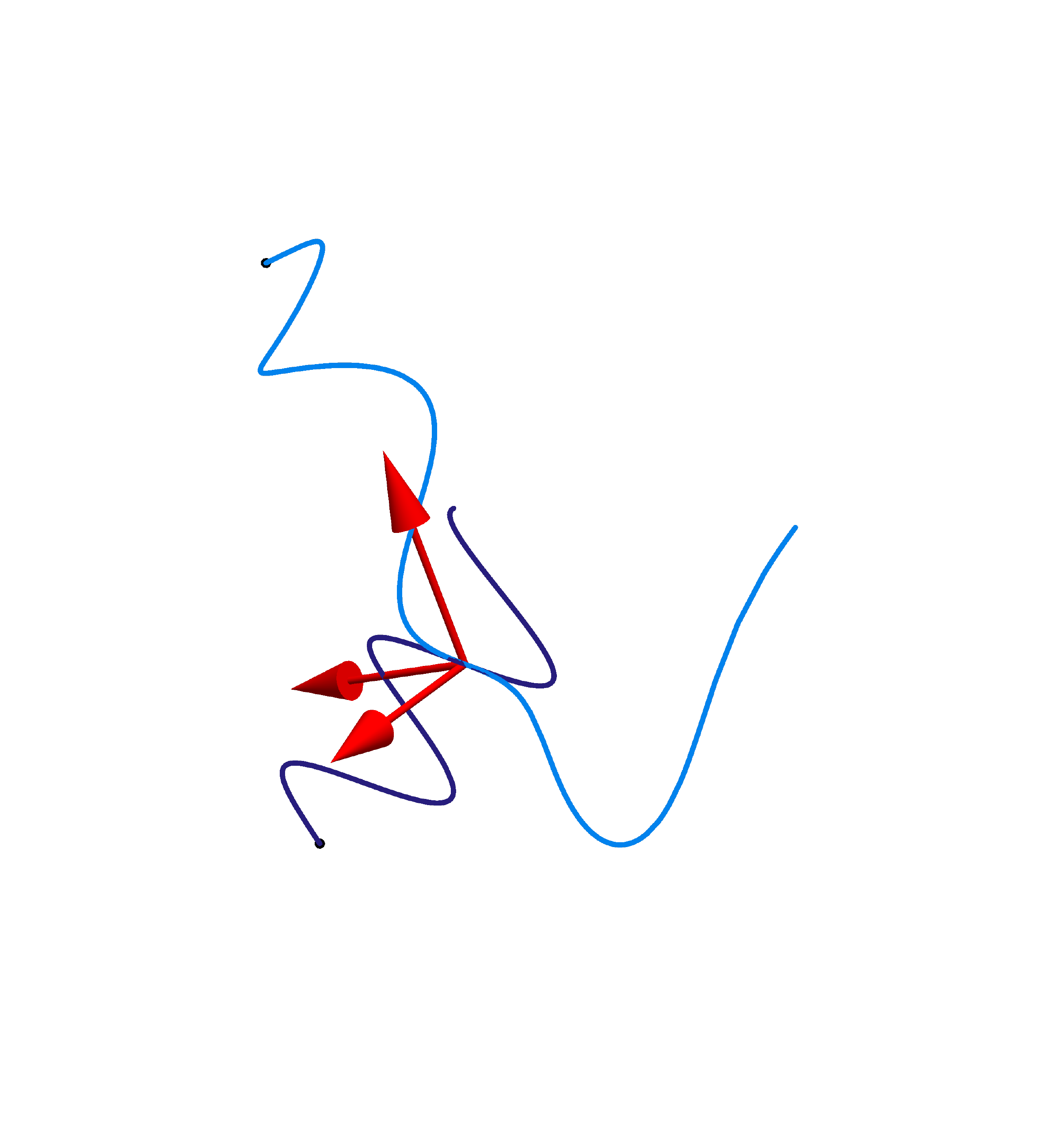}
	\vspace{-30pt}
	\endminipage
	\caption{This figure relates to \autoref{expl:torsion}. The two plots are analogous to \autoref{fig:BM}, with the manifold in question being $\bbR^3$, endowed with the connection defined above, and the path being developed is the parametrised smooth curve $X_t \coloneqq (2t \cos(t), 10\sin(t), 3t)$. The two copies $M$ and $T_0M$ of $\bbR^3$ are superimposed, with coinciding axes in the first plot. We observe how the two curves are not identical, which would be the case if the connection on $\bbR^3$ were the Euclidean one. Also note how $X$ and $\Dev(X)$ are tangent curves at their point of contact.}
	\label{fig:torsionRolling}
\end{figure}

\begin{cutout}{2}{0.8\linewidth}{0pt}{8}%
	\begin{expl}[Torsion]\label{expl:torsion}
		In general, RDEs of the form \autoref{def:RDEM} are independent of the torsion of the connections on the source and target manifolds. For parallel transport, however, torsion of $\nabla$ directly affects the field $F = \mathscr h$ that defines the RDE, and to that extent it influences the definition of $\Ptr$ and therefore that of $\Dev$ and $\Adev$ (both for the trace and second order levels of the rough paths considered). The torsion of $\widetilde \nabla$, instead, plays no role.
		To exhibit the relevance of torsion for parallel transport and development we need only focus on smooth paths. Take $M = \bbR^3$ with its canonical coordinates, and $\nabla$ with constant Christoffel symbols $\Gamma^1_{23} = 1 = -\Gamma^1_{32}$, $\Gamma^2_{31} = 1 = -\Gamma^2_{13}$, $\Gamma^3_{12} = 1 = -\Gamma^3_{21}$ and $\Gamma^k_{ij} = 0$ otherwise. This connection has the same geodesics as the Euclidean connection (straight lines), but, as described in \cite{MO20510}, parallel transport along geodesics looks like a spinning rugby ball, as illustrated in \autoref{fig:rugby} for an orthonormal frame. While the Euclidean connection and $\nabla$ agree on geodesics, they define different notions of developments: identifying $T_0\bbR^3 = \bbR^3$ we have $\Dev = \mathbbm 1$ according to the former, while this is not the case for the latter, as shown in \autoref{fig:torsionRolling}.
	\end{expl}
\end{cutout}
\vspace{-50pt}

\newpage
\begin{expl}[The dimensionality preservation question for Cartan development]\label{expl:hormander}\todo{This example is very sketchy and somewhat orthogonal to the main theme, and we should decide whether to leave it in. I think the idea behind it (and the plots) are interesting though.}
	In this example we confine ourselves to geometric/Stratonovich development, and we consider the question of whether, given a sub-vector space $P \subseteq T_oM$, there exists a submanifold $N$ of $M$, at least defined in a neighbourhood of $o$, with the property that for all $P$-valued (rough) paths $\bfZ$ taking values in $P$ and starting at $0_o$, $\Dev(\bfZ)$ is valued in $N$. Since this must hold when $Z$ is a straight line, and since straight lines develop to geodesics, if such $N$ exists it must (at least locally) be given by $\exp(P)$. Moreover, considering the case of $Z$ a piecewise linear path leads to the conclusion that, at least in the torsion-free case, $\exp(P)$ must be an affine submanifold of $M$, since it must contain every piecewise geodesic path started at $o$.
	
	At the other extreme, we may be interested in showing that, when $Z$ is an $S$-valued Brownian motion, $\Dev(Z)$ admits a density w.r.t.\ to a (hence any) Lebesgue measure on $M$. If $P = T_oM$ we should expect this to hold in view of the fact that the vectors $T_y \phi M(\mathscr H_{\lambda^\circ}(y))$ (with $\phi M$ the projection map of frame bundle) span $T_xM$ for any $y \in F_xM$. Now let $P$ be of dimension $k < m$. We can show that the first two orders of the iterated Lie brackets of the fundamental horizontal vector fields, projected down onto $TM$, are given respectively by torsion and curvature (details are omitted):
	\begin{equation}\label{eq:hormTR}
	\begin{split}
	T_y \phi M^\gamma [\mathscr H_{\mu^\circ},\mathscr H_{\nu^\circ}] &= \mathcal T^\gamma_{\beta\alpha}(x) y^{(\alpha,\mu^\circ)} y^{(\beta,\nu^\circ)} \\
	T_y \phi M^\gamma [\mathscr H_{\lambda^\circ},[\mathscr H_{\mu^\circ},\mathscr H_{\nu^\circ}]] &=  \tensor{\mathscr R}{_{\alpha\beta\delta}^\gamma}(x) y^{(\alpha,\mu^\circ)} y^{(\beta,\nu^\circ)} y^{(\delta,\lambda^\circ)}	\end{split}
	\end{equation}
	What is needed, in concrete cases, to conclude that $\Dev(Z)$ admits a smooth density is a version of the H\"ormander condition which, stated in coordinates, only applies to the first $e_1$ components of an $(e_1+e_2)$-dimensional SDE $\dif Y = F_\gamma(Y) \dif Z^\gamma$, and correspondingly only requires that $\text{Lie}(F_\gamma \colon \gamma = 1,\ldots,d)(y)$ span $\bbR^{e_1}$. This is because we are only interested in the existence of the density of $\Dev(Z)$, and not of the parallel frame above it, with which the development SDE is jointly written. See \autoref{fig:hormander} for an example of what this condition holding at different orders or not holding looks like. In view of \cite{CF10} we can expect these considerations, once made rigorously, to carry over to the case of Gaussian RDEs.
\end{expl}
\begin{figure}[h]
	\minipage[b]{0.5\textwidth}
	\includegraphics[width=\linewidth]{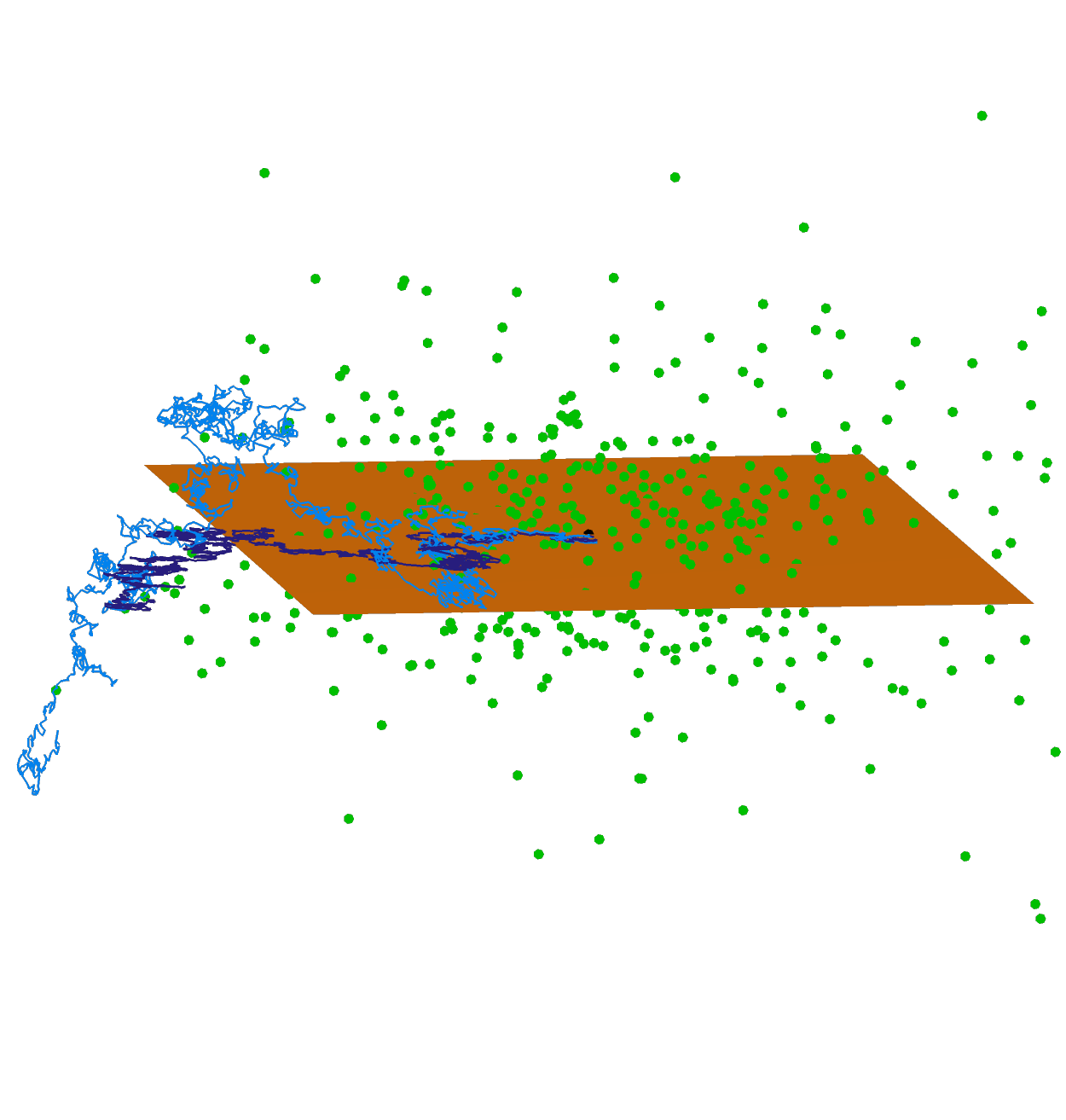}
	\vspace{-50pt}
	\endminipage\hfill
	\minipage[b]{0.5\textwidth}
	\includegraphics[width=\linewidth]{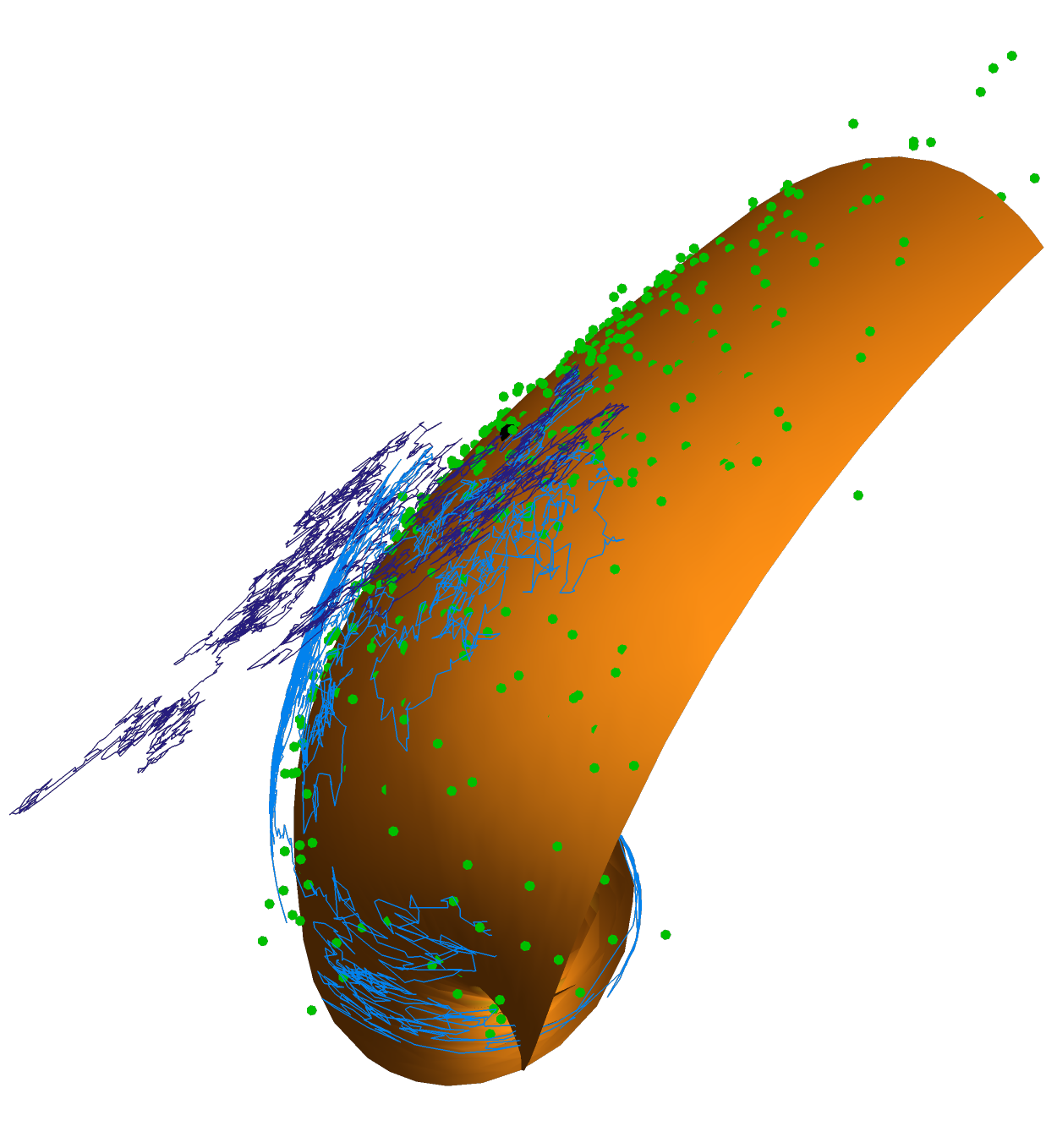}
	\vspace{-30pt}
	\endminipage\\
	\minipage[b]{0.5\textwidth}
	\includegraphics[width=\linewidth]{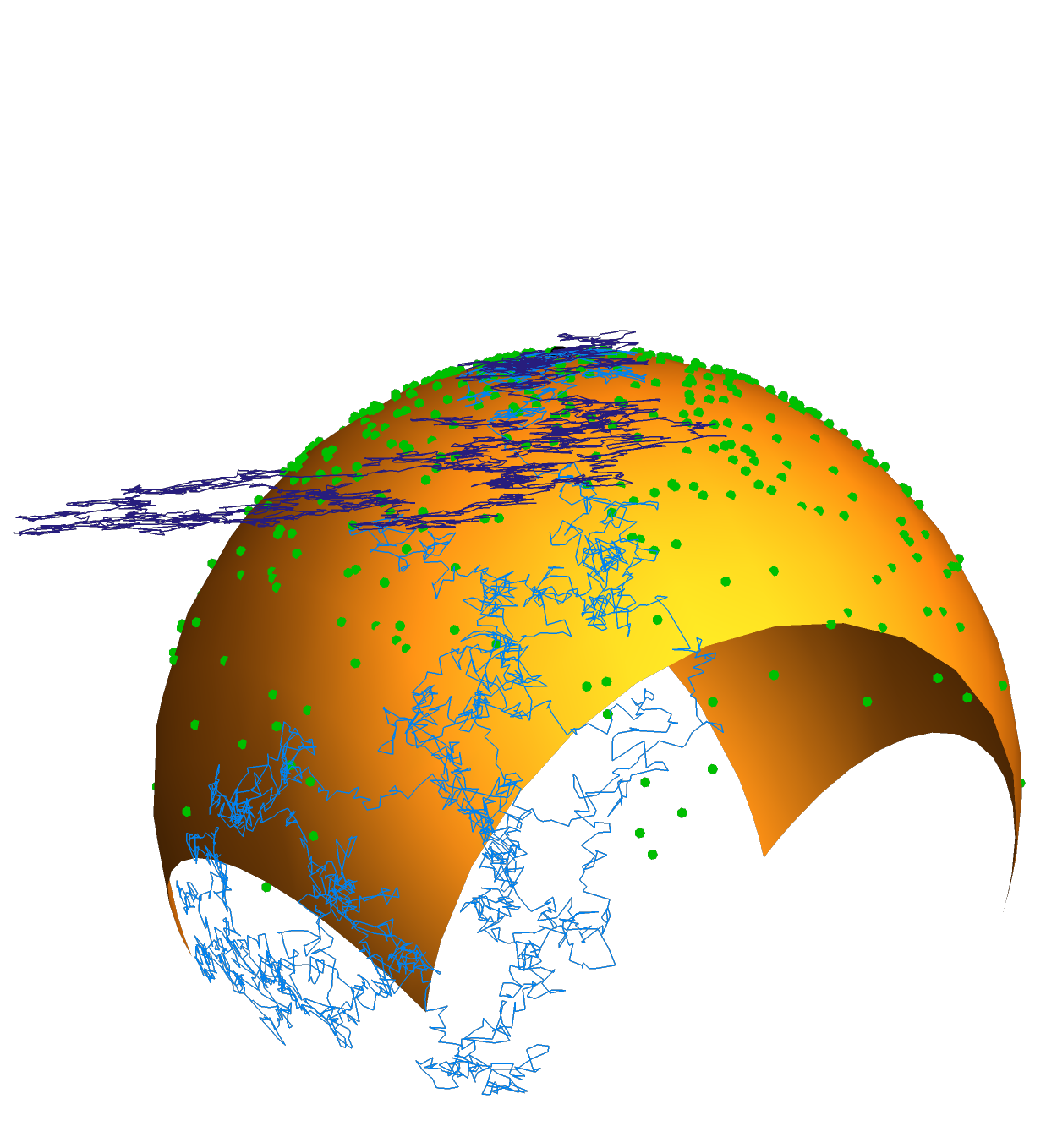}
	\vspace{-30pt}
	\endminipage\hfill
	\minipage[b]{0.5\textwidth}
	\includegraphics[width=\linewidth]{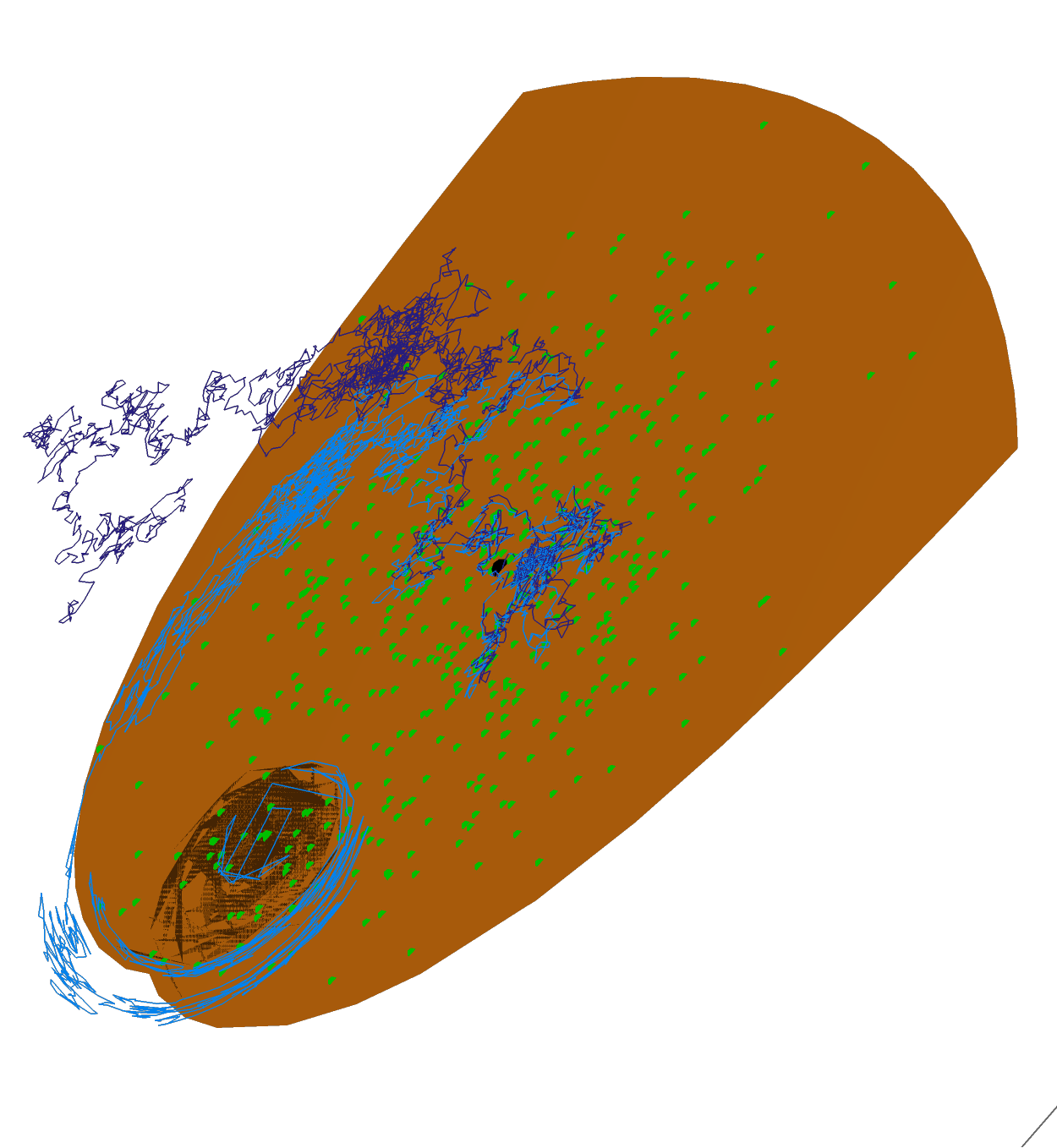}
	\vspace{-30pt}
	\endminipage
	\label{fig:hormander}
\end{figure}
\newpage

\captionof{figure}{\label{fig:hormander}The purpose of these plots is the study of different behaviours of the development of a 2-dimensional Brownian motion within a 3-dimensional manifold $M$, specifically with regards to the question of whether it remains constrained to a 2 dimensional submanifold, or whether it admits a density. In the upper left we are considering $M = \bbR^3$ with the antisymmetric connection of \autoref{expl:torsion}, while in the other three cases $M = \bbR^3 \setminus \{0\}$ with the connection whose Christoffel symbols are given by taking those for the Euclidean metric, written in spherical coordinates, and setting the $\Gamma^r$'s to zero (the geodesics in this connection are therefore circles centred at the origin and rays out of the origin, and $M$ admits the foliation into the affine submanifolds given by concentric 2-spheres centred at the origin). In each case we have plotted one 2-dimensional Brownian path valued in some subplane $P$ of $T_oM$ (in dark blue), its development onto the $3$-manifold in question (in light blue), a cloud of $\sim 5000$ points consisting of the developments of the Brownian motion at terminal time (in green), and the locally nondegenerate surface parametrised by the exponential map applied to $P$. In the first case $o$ is the origin and $P = \exp(P)$ is just the $xy$-plane, and we see how the point cloud (as well as the developed path) takes up three-dimensional space well: this is consistent with the H\"ormander condition being satisfied at order 1, thanks to torsion \autoref{eq:hormTR}. In the other three plots $o = (1,0,0)$ and $P$ is a plane intersecting the $xy$-plane in the line $(0,t,0)$ and with different inclinations w.r.t.\ the $z$-axis: $\pi/4$ in the upper right, $0$ in the bottom left and $\pi/2$ in the bottom right (this means in last two cases $P$ is respectively tangent to the unit sphere, and coincides with the $xy$-plane). Note that we have rotated the plots for improved visibility of all the components. In the plot on the upper right we see how the point cloud and the developed path do not quite adhere to $\exp(P)$, consistent with the H\"ormander condition being satisfied at order 2, thanks to curvature, but not at order 1, since the connection is torsion-free. In the other two cases $\exp(P)$ is an affine submanifold of $M$ (in the first case a sphere with the Levi-Civita connection - a leaf in the aforementioned foliation - and in the second case the punctured $xy$-plane with a non-Euclidean connection), and therefore development remains constrained to the surface plotted in orange, and does not admit a density w.r.t.\ 3-dimensional Lebesgue measure.}

\section*{Conclusions and further directions}\label{sec:concl}

In this paper we have provided a theoretical framework for the study of manifold-valued rough paths of bounded $3>p$-variation which are not required to be geometric. We have dealt with the resulting notions of rough integral and RDE, with the question of how these definitions dovetail with the previous ones for geometric rough paths and It\^o calculus on manifolds, and we have shown how to construct a theory of parallel transport and Cartan development. To conclude, we mention a few directions in which we believe these ideas can be taken further.

An obvious generalisation would consist of abandoning the regularity constraint $p<3$. A paper on geometric rough paths on manifolds of arbitrarily low regularity, of which the last two named authors are co-authors, is forthcoming. In the future we plan on transposing the theory of branched rough paths, defined originally in \cite{Gub10} and related to geometric rough paths in \cite{HK15}, to the curved setting. 

A number of topics of the last section deserve closer attention. It would be interesting to explore additional examples of connections on $\tau TM$ which result in notions of parallel transport different from the geometric/Stratonovich one: the Levi-Civita connection of the Cheeger-Gromoll metric \cite{MuTr88} and the connections defined in \cite[\S 5]{ArTh03} could be worthwhile to test. It would also be interesting to rewrite the section from the point of view of principal bundles. The challenge here is to define a correspondence between connections on $\tau TM$ and on $\tau FM$ (and optionally on $\tau T^*M$), so that the development equation can be viewed as an RDE defined by the fundamental horizontal vector fields, just as in the geometric case. Finally, we mention two ideas, only involving geometric rough paths, which were hinted at in two examples. The first would consist of extending Driver's integration by parts formula to more general Gaussian rough paths, using the notation in \autoref{expl:Driver}. The second is study the laws of developments of lower-dimensional Brownian motion, as discussed in \autoref{expl:hormander}, and to formulate the version of H\"ormander's theorem necessary to show (in certain cases) the existence of the density.

\bibliographystyle{alpha} \addcontentsline{toc}{section}{References}
{\footnotesize\bibliography{mainRefs}}

\end{document}